\magnification=\magstephalf

\baselineskip=1.1\baselineskip

\def\R{{\bf R}}

\def\eps{\varepsilon}
 \def\frac#1#2{{#1\over #2}}
 \def\ol{\overline}
\def\dist{{\bf d}}
 \def\diam{\hbox {\rm diam}}
 \def\bP{{\bf P}}
\def\bE{{\bf E}}

\def\df{{\mathop {\ =\ }\limits^{\rm{df}}}}
\def\bone{{\bf 1}}
\def\wt{\widetilde}
\def\wh{\widehat}
\def\n{{\bf n}}
\def\prt{\partial}

\def\B{{\cal B}}
\def\F{{\cal F}}
\def\S{{\cal S}}
\def\T{{\cal T}}
\def\K{{\cal K}}

\def\<{\langle}
\def\>{\rangle}
\def\[{\lfloor}
\def\]{\rfloor}

\def\thalf{{\textstyle {1 \over 2}}}
\font\tenmsb=msbm10
\font\sevenmsb=msbm7
\font\fivemsb=msbm5
\newfam\msbfam
\textfont\msbfam=\tenmsb
\scriptfont\msbfam=\sevenmsb
\scriptscriptfont\msbfam=\fivemsb
\def\Bbb#1{{\fam\msbfam\relax#1}}

\def\qed{{\hfill $\square$ \bigskip}}
\def\sqr#1#2{{\vcenter{\vbox{\hrule height.#2pt
        \hbox{\vrule width.#2pt height#1pt \kern#1pt
           \vrule width.#2pt}
        \hrule height.#2pt}}}}
\def\square{\mathchoice\sqr56\sqr56\sqr{2.1}3\sqr{1.5}3}

 \centerline{\bf SYNCHRONOUS COUPLINGS }
 \centerline{ \bf OF REFLECTED BROWNIAN
  MOTIONS IN SMOOTH DOMAINS}

\footnote{$\empty$}{\rm Research partially supported by NSF grant
DMS-0303310.} \vskip0.3truein

\centerline{\bf Krzysztof Burdzy{\rm,} \quad Zhen-Qing Chen \quad
{\rm and} \quad Peter Jones }

\vskip1truein

{\narrower\narrower \noindent{\bf Abstract}. For every bounded
planar domain $D$ with a smooth boundary, we define a ``Lyapunov
exponent'' $\Lambda(D)$ using a fairly explicit formula. We
consider two reflected Brownian motions in $D$, driven by the same
Brownian motion (i.e., a ``synchronous coupling''). If
$\Lambda(D)>0$ then the distance between the two Brownian
particles goes to $0$ exponentially fast with rate $\Lambda
(D)/(2|D|)$ as time goes to infinity. The exponent $\Lambda(D)$ is
strictly positive if the domain has at most one hole. It is an
open problem whether there exists a domain with $\Lambda(D)<0$.

 }

\bigskip
\noindent{\bf 1. Introduction and main results}.

Suppose $D\subset\R^2$ is an open connected bounded set with
$C^4$-smooth boundary, not necessarily simply connected. Let
$\n(x)$ denote the unit inward normal vector at $x\in\prt D$. Let
$B$ be standard planar Brownian motion and consider the following
Skorokhod equations,
 $$\eqalignno{
 X_t &= x_0 + B_t + \int_0^t  \n(X_s) dL^X_s
 \qquad \hbox{for } t\geq 0, &(1.1)\cr
 Y_t &= y_0 + B_t + \int_0^t  \n(Y_s) dL^Y_s
 \qquad \hbox{for } t\geq 0. &(1.2) }
$$
Here $L^X$ is the local time of $X$ on $\prt D$. In other words,
$L^X$ is a non-decreasing continuous process which does not
increase when $X$ is in $D$, i.e., $\int_0^\infty \bone_{D}(X_t)
dL^X_t = 0$, a.s. Equation (1.1) has a unique pathwise solution
$(X,L^X)$ such that $X_t \in \ol D$ for all $t\geq 0$ (see [LS]).
The reflected Brownian motion $X$ is a strong Markov process. The
same remarks apply to (1.2), so $(X, Y)$ is also strong Markov. We
will call $(X, Y)$ a ``synchronous coupling.'' Note that on any
interval $(s,t)$ such that $X_u \in D$ and $Y_u \in D$ for all $u
\in (s,t)$, we have $X_u - Y_u = X_s - Y_s$ for all $u \in (s,t)$.

Before we state our main results, we will introduce some notation
and make some technical assumptions on $\prt D$. We will assume
that for every point $x\in\prt D$, there exists a neighborhood $U$
of $x$ and an orthonormal system $CS_x$ such that $\n(x) =(0,1)$
and $x=(0,0)$ in $CS_x$, and $\prt D \cap U$ is a part of the
graph of a function $y_2 = \psi_x(y_1)$ satisfying $ \psi_x(y_1) =
(1/2)\nu(x) y_1^2 + O(y_1^3)$. This defines the curvature $\nu(x)$
for $\partial D$ at $x$. We will assume that there is $N_1<\infty$
such that for every unit vector ${\bf m}$, there are at most $N_1$
points $x\in \prt D$ with $\n(x) = {\bf m}$. Recall that $\prt D$
is assumed to be $C^4$-smooth. We will assume that there is only a
finite number of $x\in\prt D$ with $\nu(x) =0$ and that for every
such $x$, we have $\psi_x(y_1)= c_x y_1^3 + O(y_1^4)$ with $c_x\ne
0$. The distance between $x$ and $y$ will be denoted $\dist(x,y)$.

\bigskip
\noindent{\bf Theorem 1.1}. {\sl  If $D$ satisfies the above
assumptions and it has at most one hole then $\dist(X_t, Y_t) \to
0$ as $t\to \infty$, a.s., for every pair of starting points
$(x_0, y_0)\in \ol D \times \ol D $. }

\bigskip

The above theorem complements the results in [BC] where it has
been proved that the distance between $X_t$ and $Y_t$ converges to
0 as $t\to \infty$ for two classes of domains: (i) polygonal
domains, i.e., domains whose boundary consists of a finite number
of closed polygons, and (ii) ``lip domains'', i.e., bounded
Lipschitz domains which lie between graphs of two Lipschitz
functions that have Lipschitz constants  strictly less than 1. The
number of holes plays no role in the case of polygonal domains but
it is an open problem whether it does in the case of smooth
domains (see Section 2).

Earlier research of Cranston and Le Jan ([CLJ1, CLJ2]) on
synchronous couplings of reflected Brownian motions was focused on
convex domains. In that case, it is clear that $t \to \dist(X_t,
Y_t)$ is non-increasing. Cranston and Le Jan proved that for a
large class of convex domains, $\dist(X_t, Y_t) >0$ for all $t\geq
0$, a.s., if $\dist(X_0,Y_0)>0$. The present paper, especially
Theorem 1.2 below, answers a problem posed at the end of [CLJ1]
and improves on the estimate given in the Appendix of [CLJ2].

\bigskip

Next we will present our main technical result on the ``Lyapunov
exponent,'' which is a crucial step in the proof of Theorem 1.1.
We need some more notation. Let $\sigma^X_t = \inf\{s\geq 0: L^X_s
\geq t\}$. For every bounded planar domain $D$ we have $\lim_{t\to
\infty} L^X_t = \infty$ so $\sigma^X_t < \infty$ for all $t\geq
0$, a.s. The arc length measure on $\prt D$ will be denoted
``$dx$'', e.g., we will write $\int_{\prt D} f(x) dx$ to denote
the integral of $f$ with respect to the arc length. For any
$x,y\in \prt D$, we let $\alpha(x,y)$ be the angle formed by the
tangent lines to $\prt D$ at $x$ and $y$, with the convention that
$\alpha(x,y)\in[0,\pi/2]$. For every point $x\in \prt D$, let
$\omega_x (dy)$ be the ``harmonic measure'' on $\prt D$ with the
base point $x$, defined as follows. Let $K(x,y)$, $y\in D$, be the
Martin kernel in $D$ with the pole at $x$, i.e., the only (up to a
multiplicative constant) positive harmonic function in $D$ which
vanishes everywhere on the boundary of $D$ except for a pole at
$x$. Then we let $\omega_x(dy) = a_x \, {\prt K(x, y) \over \prt
 \n(y)} dy$
where the constant $a_x $ is chosen so that $\lim_{y \to x} \pi
\dist(x,y)^2 \omega_x(dy)/dy = 1 $. Let $|D|$ denote the area of
$D$.

\bigskip
\noindent{\bf Theorem 1.2}. {\sl Let
$$\Lambda(D) = \int_{\prt D} \nu(x) dx
 + \int_{\prt D} \int_{\prt D} |\log \cos\alpha(x,y)|\omega_x(dy) dx.
\eqno(1.3)
 $$
If $\Lambda(D)  > 0$, then for any $x_0,y_0\in\ol D$, a.s.,
 $$\lim_{t\to \infty} { \log \dist(X_{t}, Y_{t})
\over t} = - {\Lambda(D)\over  2 |D|}. \eqno(1.4)
 $$

 }

\bigskip

By the Gauss-Bonnet Theorem, the first integral in (1.3), that is,
$\int_{\prt D} \nu(x) dx$, is equal to $2\pi \chi (D)$, where
$\chi(D)$ is the Euler characteristic of $D$. In our case,
$\chi(D)$ is equal to 1 minus the number of holes in $D$. We are
not aware of a simple representation of the second (double)
integral in (1.3). The integral $\int_{\prt D} \nu(x) dx$, which
appeared in [CLJ2], emerges in our arguments as the limit of
$(1/t)\int_0^t \nu(X_s)dL^X_s$ when $t\to\infty$. See [H] for some
results involving $\int_0^t \nu(X_s)dL^X_s$.

It is elementary to check using the definition (1.3) that
$\Lambda(D)$ is invariant under scaling, i.e., for any $a>0$,
$\Lambda(D) = \Lambda(a D)$, where $a D = \{x\in\R^2: x = a y
\hbox{  for some } y\in D\}$.

\bigskip
We will now explain the intuitive content of Theorem 1.2. The disc
with center $x$ and radius $r$ will be denoted ${\cal B}(x,r)$.
Suppose that at some time $t$, $\dist(X_t,Y_t)$ is very small so
that when one of the processes is on the boundary of the domain
then $\prt D$ looks like a very flat parabola inside the disc
${\cal B}(X_t, 2\dist(X_t,Y_t))$. Suppose further that the line
segment $\ol{X_t,Y_t}$ is ``almost'' parallel to $\prt D$. Then
the local time components in (1.1) and (1.2) will be almost
identical over a short time period $[t,t+\Delta t]$, except for a
small difference between the reflection vectors due to the
curvature of $\prt D$. This small difference translates into the
first integral in (1.3). From time to time, $X$ makes large
excursions from $\prt D$, whose endpoints are at a distance
comparable to the diameter of $D$. At the end of any such
excursion, one and only one of the processes $X$ or $Y$ gets a
substantial local time push, until again $\ol{X_t,Y_t}$ is almost
parallel to $\prt D$. This results in the reduction of
$\dist(X_t,Y_t)$ by a factor very close to $\cos \alpha(x,y)$,
where $x$ and $y$ are the endpoints of the excursion. The double
integral on the right hand side of (1.3) represents the change in
$\dist(X_t,Y_t)$ due to large excursions. We find it surprising
and intriguing that the magnitudes of the two phenomena affecting
the distance $\dist(X_t,Y_t)$, described above, are comparable and
give rise to two ``independent'' terms on the right hand side of
(1.3).

\bigskip
We will briefly sketch the idea behind the proof of Theorem 1.2.
First, we prove that the distance between the particles will be
small at least from time to time, so that we can apply methods
appropriate for processes reflecting on very flat parabolas. The
main part of the proof deals with the two phenomena described in
the previous paragraph. When the line segment $\ol{X_t,Y_t}$ is
``almost parallel'' to $\prt D$ and one or both processes reflect
on $\prt D$, the change of $\dist(X_t,Y_t)$ is ``almost''
deterministic in nature and so are our methods. The change in
$\dist(X_t,Y_t)$ due to ``large'' excursions of $X_t$ in $D$ is
much harder to analyze and that part of the proof is very
complicated. We list here several of the challenges. First, it is
conceivable that even a single excursion may result in a reduction
of $\dist(X_t,Y_t)$ to $0$, if the endpoints of the excursion are
at $x,y\in \prt D$ with $\cos \alpha(x,y)=0$. Proving that this is
not the case takes considerable effort. Second, we use excursion
theory and ergodicity of $X_t$ to prove that $\log \dist(X_t,Y_t)$
obeys a strong law of large numbers, in the sense of (1.4). The
problem here is that although $X_t$ is recurrent and ergodic, the
vector process $(X_t,Y_t)$ is neither, and so we have to analyze
the behavior of $Y_t$ by proving that it is ``close'' to that of
$X_t$. Finally, one has to find upper bounds for probabilities of
various ``unusual'' events which clearly cannot happen, from the
intuitive point of view, but which have to be accounted for in a
rigorous argument.

\bigskip
The rest of the paper is organized as follows. Section 2 is
devoted to the discussion of some open problems and examples,
mostly related to Theorem 1.1. It also contains a (very short)
proof of Theorem 1.1. The proof of Theorem 1.2, consisting of many
lemmas, is given in Sections 3 and 4. Most arguments in Section 3
are deterministic or analytic in nature. Section 4 contains
arguments based on the excursion theory.

\bigskip

We are grateful to Greg Lawler, Nick Makarov, Don Marshall and
B\'alint Vir\'ag for
 very useful discussions and advice.

\bigskip

\noindent{\bf 2. Examples and open problems}.

The paper was inspired by the following problem which still
remains open.

\bigskip
\noindent{\bf Problem 2.1}. {\sl (i) Does there exist a bounded
planar domain such that with positive probability,
$$ \limsup_{t\to\infty} \dist(X_t,Y_t) > 0?
$$

(ii) Does there exist a bounded domain $D$ with $\Lambda(D) < 0$?

 }
\bigskip

The two problems are related to each other via the following
conjecture.

\bigskip
\noindent{\bf Conjecture 2.2}. {\sl If $\Lambda(D) < 0$ then with
probability one, $\limsup_{t\to \infty} \dist(X_t,Y_t) > 0$. }
\bigskip

We believe that the above conjecture can be proved using the same
methods as in the proof of Theorem 1.2. Since we do not know
whether any domains with $\Lambda(D) <0$ exist, we have little
incentive to work out the details of the proof for Conjecture 2.2.

A technical problem arises in relation to Problem 2.1 (i)---it is
not obvious how to define a ``synchronous coupling'' of reflected
Brownian motions in an arbitrary domain. It is desirable from both
technical and intuitive point of view to have the strong Markov
property for the process $(X_t,Y_t)$. See [BC] for a discussion of
these points. So far, the existence of synchronous couplings of
reflected Brownian motions with the strong Markov property can be
proved only in these domains where the stochastic Skorokhod
equations (1.1)-(1.2) have a unique strong solution. A recent
paper ([BBC]) shows that this is the case when $D$ is a planar
Lipschitz domain with the Lipschitz constant less than 1.

We will next present some speculative directions of research
related to Problem 2.1 (ii). We start by explaining how Theorem
1.1 follows from Theorem 1.2.

\bigskip
\noindent{\bf Proof of Theorem 1.1}. If $D$ has at most one hole
then the first integral on the right hand side of (1.3) is equal
to $2\pi$ or 0, by the Gauss-Bonnet Theorem. The integrand in the
double integral in (1.3) is non-negative and it is easy to see
that it is strictly positive on a non-negligible set. Hence,
$\Lambda(D)
>0$ and, consequently, (1.4) holds, according to Theorem 1.2.
Thus, Theorem 1.1 follows from Theorem 1.2. \qed

\bigskip

The above proof suggests the following strategy for finding a
domain with $\Lambda(D) < 0$. One should find a domain where the
first integral on the right hand side of (1.3) is significantly
less than zero. This is because the contribution from the second
term is always non-negative. In other words, one has to consider
domains with many holes because, as we have already mentioned in
Section 1, the first term is equal to 1 minus the number of holes,
multiplied by $2\pi$. The obvious problem with this strategy is
that punching holes in a domain may increase the double integral
on the right hand side of (1.3), and this may offset the effect of
holes on the first integral.

Here is a possible avenue of research based on the above idea.
Suppose that $D$ has a large number of small holes. Here ``small''
means that the holes have diameters very small in comparison with
the diameter of the domain. Let us assume that distances between
different holes, and distances between holes and the outside
boundary of $D$ are large in comparison with diameters of holes.
Then it is not hard to see that the right hand side of (1.3) is
very close to the sum of analogous formulas for each connected
component of $\prt D$. In other words, there is little interaction
between different connected components of $\prt D$, if the holes
are small and far apart. If we can find a hole with the shape
which yields $\Lambda (D)< 0$ for a single hole, then
$\Lambda(D)<0$ for a domain $D$ with a large number of small holes
of this shape.

Simple heuristic estimates show that if $\prt D$ has
``approximate'' corners, like a polygonal domain (the corners have
to be ``approximate'' because the domain has to be smooth), then
the double integral on the right hand side of (1.3) is very large.
Domains, or rather holes, with this property will not help us in
our search for a domain $D$ with $\Lambda(D)<0$. The ultimate
domain without corners is a disc. At the moment we are concerned
with ``holes'' so we will find $\Lambda(D)$ for $D$ which is the
exterior of a disc.

Recall that a disc with center $x$ and radius $r$ is denoted
$\B(x,r)$.

\bigskip

\noindent{\bf Proposition 2.3}. {\sl If $D=\B((0,0),1)^c$ then
$\Lambda(D)=0$.}

\bigskip

\noindent{\bf Proof}. Recall that the first integral in (1.3) is
equal to $-2\pi$. We will parametrize $\prt D$ using
$\theta\in[0,2\pi)$ and writing $x=e^{i\theta}$ for $x\in \prt D$.
The formula for the harmonic measure in $D$ is well known and easy
to derive using standard complex analytic methods (conformal
mappings). This easily leads to the following formula for
the ``harmonic measure'' $\omega_x$,
 $${\omega_x(dy)\over dy}
 = {1 \over 4 \pi \sin^2 ({\theta -\theta'\over 2})},$$
where $x=e^{i\theta}$ and $y = e^{i\theta'}$. Thus the double
integral in (1.3) is equal to
 $$\int_0^{2\pi} \int_0^{2\pi}
 {|\log|\cos(\theta-\theta')||
 \over 4 \pi \sin^2 ({\theta -\theta'\over 2})} d\theta d\theta'
 = 2\pi\int_0^{2\pi}
 {|\log|\cos(\theta)||
 \over 4 \pi \sin^2 (\theta / 2)} d\theta
 = \int_0^\pi
 {|\log|\cos(\theta)||
 \over  \sin^2 (\theta / 2)} d\theta .
 $$
We have
 $$\eqalign{&\int_0^{\pi/2}
 {|\log|\cos(\theta)||
 \over  \sin^2 (\theta / 2)} d\theta
 =-\int_0^{\pi/2}
 {\log\cos(\theta)
 \over  \sin^2 (\theta / 2)} d\theta\cr
 &= [2(\theta + \cot(\theta/2) \log \cos \theta
 - \log(\cos(\theta/2) - \sin(\theta/2))
 + \log (\cos(\theta/2) +
 \sin(\theta/2))]\Big|_{\theta=0}^{\theta=\pi/2}\cr
 &= \pi + 2 \log 2,
 }$$
and
 $$\eqalign{\int_{\pi/2}^\pi
 {|\log|\cos(\theta)||
 \over  \sin^2 (\theta / 2)} d\theta
 &=-\int_{\pi/2}^\pi
 {\log(-\cos(\theta))
 \over  \sin^2 (\theta / 2)} d\theta\cr
 &= \Big[2(\theta + \cot({\theta/ 2}) \log(- \cos \theta)
 - \log( \sin({\theta/ 2})-\cos({\theta/ 2}))\cr
 &\qquad + \log (\cos({\theta/ 2}) +
 \sin({\theta/ 2}))\Big]\Bigg|_{\theta=\pi/2}^{\theta=\pi}\cr
 &= \pi - 2 \log 2,
 }$$
so
$$\int_0^{2\pi} \int_0^{2\pi}
 {|\log|\cos(\theta-\theta')||
 \over 4 \pi \sin^2 ({\theta -\theta'\over 2})} d\theta d\theta'
 = \int_0^\pi
 {|\log|\cos(\theta)||
 \over  \sin^2 (\theta / 2)} d\theta = 2\pi,
 $$
and $\Lambda (D) = 0$. \qed

\bigskip

We have proved that $\Lambda(D)=0$ for the exterior of a disc by a
brute force calculation. It is a natural question whether the same
result follows from some elegant symmetry argument---we have not
found one so far.

Since $\Lambda(D) =0$ for the exterior of the disc, discs are not
helpful as holes in the (hypothetical) construction of a domain
$D$ with $\Lambda(D)<0$. Our next observation is that the exterior
of a line segment would be a great candidate for a useful hole.
This is because for any points $x$ and $y$ on a line segment, we
have $\alpha(x,y) =0$ and, therefore, the double integral on the
right hand side of (1.3) vanishes for the exterior of a line
segment. Hence, $\Lambda(D)<0$ for the exterior of a line segment.
Unfortunately, we cannot use line segments as holes because their
boundaries are not smooth. Instead, we can try a domain ``close''
to a line segment but with a smooth boundary. A natural candidate
is a very elongated ellipse. Our preliminary numerical
calculations showed that $\Lambda(D) =0$ for the exterior of any
ellipse. We are grateful to B\'alint Vir\'ag for the following
rigorous proof of this result.

\bigskip

\noindent{\bf Proposition 2.4}. {\sl If $D$ is the exterior of an
ellipse then $\Lambda(D)=0$.}

\bigskip

\noindent{\bf Proof}. We will use complex analysis and complex
notation in this proof. Recall that $\Lambda(D)$ is invariant
under scaling. Hence, we can consider any ellipse with the given
eccentricity. In other words, it is enough to prove that the
proposition holds for any ellipse $D$ that can be represented as
$D= g(U)$, where $U=\B(0,1)$, $g(z) = z + a/z$, and $a$ is any
real number in $ (0,1)$.

We start by proving the following claim. Suppose that $f$ is an
analytic function in $\ol U$, $f(1) $ is purely imaginary, and
$f'(1)$ is real. Then
 $${1\over \pi} \int _{\prt U} {\Re f (x) \over |1-x|^2} |dx|
 = - f'(1). \eqno (2.1)
 $$
Since $\Re f(y)$ is harmonic in $U$ and continuous on $\ol U$,
 $$\Re f(y) = {1\over 2\pi}
 \int_{\prt U} \Re f(x) {1- |y|^2 \over |y-x|^2} |dx|,$$
for $y\in U$. We have assumed that $f(1)$ is imaginary and $f'(1)$
is real, so, by dominated convergence,
 $$\eqalign{
 -f'(1)&= \lim_{r\to 0, r>0} {\Re f(1-r) - \Re f(1)\over r}
 = \lim_{r\to 0, r>0} {1\over 2\pi} \int_{\prt U} \Re f(x)
  {1- |1-r|^2 \over r |1-r-x|^2} |dx|\cr
  &= {1\over \pi} \int _{\prt U} {\Re f (x) \over |1-x|^2} |dx|.
 }$$
We have shown that (2.1) holds.

The first integral in (1.3) is equal to $-2\pi$. It will suffice
to show that the second (double) integral is equal to $2\pi$. The
second integral in (1.3) is equal to
 $${1\over \pi} \int_{\prt U} \int_ {\prt U}
 { \log \Re \left( {x g'(x) \over |g'(x)|}\cdot {|g'(y)| \over y g'(y)}
 \right) \over |x-y|^2 } |dx| |dy|.
 $$
Note that $g'(z) = 1 - a/z^2$ and for $z\in \prt U$, $\ol z =
1/z$. Let $\beta = 1/(2 g'(y) y)$ for some $y\in \prt U$. If we
write $h(x) = \Re\left ( {x g'(x) \over  y g'(y)} \right)$ then
for $x\in \prt U$, we have
 $$h(x)= \Re (2 \beta g'(x) x) = \beta g'(x) x + \ol \beta \ol {g'(x)}
 \ol x
 = \beta (1 - a/x^2 ) x + \ol \beta (1 - a x^2) (1/x).
 $$
We have
 $$\eqalign{
 \log \Re \left( {x g'(x) \over |g'(x)|}\cdot {|g'(y)| \over y g'(y)}
 \right)
 &= \log \left( {|g'(y)| \over |g'(x)|}\cdot
 \Re\left ( {x g'(x) \over  y g'(y)} \right)
 \right)
 = \log \left( {|g'(y)| \over |g'(x)|}
 \right) + \log h(x)\cr
 &= \Re \log \left( {g'(y) \over g'(x)}
 \right) + \Re \log h(x)
 = \Re \log \left(h(x) \cdot {g'(y) \over g'(x)}
 \right).
 }$$
Hence,
 $${1\over \pi} \int_{\prt U}
 { \log \Re \left( {x g'(x) \over |g'(x)|}\cdot {|g'(y)| \over y g'(y)}
 \right) \over |x-y|^2 } |dx|
 ={1\over \pi} \int_{\prt U}
 { \Re \log \left(h(x) \cdot {g'(y) \over g'(x)}
 \right) \over |x-y|^2 } |dx|.
 $$
Recall that $|y|=1$. Substituting $x = y/v$, we see that the last
integral is equal to
 $${1\over \pi} \int_{\prt U}
 { \Re \log \left(h(y/v) \cdot {g'(y) \over g'(y/v)}
 \right) \over |x-y|^2 }{|x|^2\over |y|} |dv|
 = {1\over \pi} \int_{\prt U}
 { \Re \log \left(h(y/v) \cdot {g'(y) \over g'(y/v)}
 \right) \over |1-v|^2 } |dv|.
 $$
We have
 $$\eqalign{
 h(y/v) \cdot {g'(y) \over g'(y/v)}
 &= {(\beta (1 - a/(y/v)^2 ) (y/v) +
 \ol \beta (1 - a (y/v)^2) (1/(y/v))) g'(y)
 \over 1 - a/(y/v)^2 }\cr
 &= g'(y) {y\over v} {v^2 (\ol \beta - a \beta)
 + y^2 (\beta - a \ol \beta)
 \over y^2 - av^2 }.
 }$$
Let
 $$k(v) = v
 \left(h(y/v) \cdot {g'(y) \over g'(y/v)}\right ) =
 g'(y) y {v^2 (\ol \beta - a \beta)
 + y^2 (\beta - a \ol \beta)
 \over y^2 - av^2 },
 $$
and note that, since $|v|=1$,
 $$ {1\over \pi} \int_{\prt U}
 { \Re \log \left(h(y/v) \cdot {g'(y) \over g'(y/v)}
 \right) \over |1-v|^2 } |dv|
 = {1\over \pi} \int_{\prt U}
 { \Re \log k(v) \over |1-v|^2 } |dv|.
 $$

Next we will verify that (2.1) can be applied to $f(x) = \log
k(x)$. We have for $y \in \prt U$,
 $$\beta = {1 \over 2 g'(y) y} =
 { 1 \over 2 (1 - a/y^2) y } = { y \over 2(y^2 - a)},
 $$
 $$\ol \beta = { \ol y \over 2(\ol y^2 - a)}
 = { 1/y \over 2 ( (1/y)^2 - a)} = { y \over 2(1 - a y^2)},
 $$
 $${\ol \beta \over \beta} = {y^2 -a \over 1 - a y^2},$$
 $${\beta - a \ol \beta\over
 \ol \beta - a \beta}
 = {{ y \over 2(y^2 - a)} - a { y \over 2(1 - a y^2)}\over
 { y \over 2(1 - a y^2)} - a { y \over 2(y^2 - a)}}
 ={ 1+a^2 -2a y^2 \over y^2 - 2a + a^2 y^2},
 $$
 $$k(v) =
 g'(y) y {v^2 (\ol \beta - a \beta)
 + y^2 (\beta - a \ol \beta)
 \over y^2 - av^2 }
 = {1 \over 2 \beta}
 {v^2 (\ol \beta - a \beta)
 + y^2 (\beta - a \ol \beta)
 \over y^2 - av^2 },
 $$
 $$\eqalign{
 k(1) &= {1 \over 2 \beta}
 { (\ol \beta - a \beta)
 + y^2 (\beta - a \ol \beta)
 \over y^2 - a }
 =  { ( \ol \beta/\beta - a )
 + y^2 ( 1- a \ol \beta/\beta)
 \over 2(y^2 - a) }\cr
 & = { ({y^2 -a \over 1 - a y^2} - a )
 + y^2 ( 1- a {y^2 -a \over 1 - a y^2})
 \over 2(y^2 - a) } =1,
 }$$
 $$\log k(1) = 0,$$
 $$k'(v) =
 g'(y) y {2v (\ol \beta - a \beta) (y^2 - av^2)
 +2av[v^2 (\ol \beta - a \beta)
 + y^2 (\beta - a \ol \beta)]
 \over (y^2 - av^2)^2 },
 $$
 $$
 k'(v)/k(v)
 = {2v (\ol \beta - a \beta) (y^2 - av^2)
 +2av[v^2 (\ol \beta - a \beta)
 + y^2 (\beta - a \ol \beta)]
 \over (y^2 - av^2)[v^2 (\ol \beta - a \beta)
 + y^2 (\beta - a \ol \beta)] },
 $$
 $$\eqalign{
 (\log k)'(1) &= k'(1) / k(1)
 ={2 (\ol \beta - a \beta) (y^2 - a)
 +2a[ (\ol \beta - a \beta)
 + y^2 (\beta - a \ol \beta)]
 \over (y^2 - a)[ (\ol \beta - a \beta)
 + y^2 (\beta - a \ol \beta)] }\cr
 &= {2 (\ol \beta - a \beta)
 \over  (\ol \beta - a \beta)
 + y^2 (\beta - a \ol \beta) }
 +
 {2a \over (y^2 - a) }\cr
 &= {2 \over  1
 + y^2 \left({\beta - a \ol \beta\over
 \ol \beta - a \beta}\right) }
 +
 {2a \over (y^2 - a) }\cr
  &= {2 \over  1
 + y^2 \left({ 1+a^2 -2a y^2 \over y^2 - 2a + a^2 y^2}\right) }
 +
 {2a \over (y^2 - a) }\cr
 &= {2(y^2 - 2a + a^2 y^2) \over
 2y^2 -2a + 2a^2y^2 -2 a y^4 }
 +
 {2a \over (y^2 - a) }\cr
 &= { (1-a^2) y^2 \over (y^2 - a)( 1-ay^2)}\cr
 &= (1-a^2) {1 \over (1-a/y^2)(1-ay^2)}\cr
 & = (1-a^2) {1 \over |1- ay^2|^2} \in \R.
 }$$
Since $\Re \log k(1) = 0$ and $(\log k)'(1)$ is real, we can apply
(2.1) to obtain
 $$ {1\over \pi} \int_{\prt U}
 { \Re \log k(v) \over |1-v|^2 } |dv| = -(\log k)'(1),
 $$
and
$$\eqalignno{
{1\over \pi}& \int_{\prt U} \int_ {\prt U}
 { \log \Re \left( {x g'(x) \over |g'(x)|}\cdot {|g'(y)| \over y g'(y)}
 \right) \over |x-y|^2 } |dx| |dy|
 =  \int_ {\prt U}{1\over \pi} \int_{\prt U}
 { \Re \log k(v) \over |1-v|^2 } |dv| |dy|\cr
 &= -\int_ {\prt U}(\log k)'(1) |dy|
 = -\int_ {\prt U}{ (1-a^2) y^2 \over (y^2 - a)( 1-ay^2)}|dy|\cr
 & = i\int_ {\prt U}{ (1-a^2) y^2 \over (y^2 - a)( 1-ay^2)}
 {1\over y} dy.&(2.2)
 }$$
The function ${ (1-a^2) y \over (y^2 - a)( 1-ay^2)}$ has two poles
inside $U$, at $y= \pm \sqrt{a}$, and the residue is equal to 1/2
at each of these points. Hence, by the residue theorem, the right
hand side of (2.2) is equal to $2\pi$. This completes the proof of
the proposition.
 \qed

\bigskip

The last result raises some questions, but before we state them as
a formal conjecture, we rush to add that it is very easy to see
that $\Lambda(D)>0$ for exteriors of some convex smooth domains,
for example, those that have ``approximate'' corners.

\bigskip

\noindent{\bf Conjecture 2.5} {\sl (i) If $D$ is the exterior of a
simply connected domain then $\Lambda(D) \geq 0$.

(ii) If $D$ is the exterior of a simply connected domain and
$\Lambda(D) =0$ then $D^c$ is a disc or an ellipse.

}
\bigskip

Another problem, hard to state as a formal conjecture, is to find
an (easy) way to derive $\Lambda(D)$ for the exterior of an
ellipse from the value of this constant for the exterior of a
disc. We point out an obvious fact that $\Lambda(D)$ is not
invariant under conformal mappings. It is not hard to see that
$\Lambda(D)$ is not invariant under the transformation $(x_1,x_2)
\mapsto (cx_1, x_2)$.

\bigskip
\noindent{\bf Problem 2.6}. {\sl Let $D$ be the exterior of a
disc. Is it true that $\dist (X_t, Y_t) \to 0$ as $t\to\infty$,
a.s.?}
\bigskip

The last problem might be hard because it deals with the
``critical'' case, i.e., the case when $\Lambda(D)=0$. On the
other hand, the symmetries of the disc might be the basis of a
reasonably easy proof, specific to this domain.

\bigskip
\noindent{\bf 3. Analysis of Skorokhod transforms}.

{\it Notation}. The following notation will be used throughout the
paper.

All constants $c_1,c_2, \dots$ will take values in $(0,\infty)$
unless stated otherwise. We will write $a\lor b = \max(a,b)$ and
$a\land b = \min(a,b)$. Recall that the distance between
$x,y\in\R^2$ is denoted as $\dist(x,y)$; the same symbol will be
used to denote the distance between a point and a set, etc. Our
arguments will involve elements of $\R$ or $\R^2$, and one- or
two-dimensional vectors. We will use $|\,\cdot\,|$ to denote the
usual Euclidean norm in all such cases. For $x,y\in \R^2$, the
meaning of $|x-y|$ is the same as that of $\dist(x,y)$ but we will
nevertheless find it convenient to use both pieces of notation.
The disc with center $x$ and radius $r$ will be denoted $\B(x,r)$.
Recall the definition of curvature $\nu(x)$ at a point $x\in \prt
D$, from the Introduction and let $\nu^* = \sup_{x\in \prt D}
|\nu(x)|$. The unit inward normal vector at $x\in \prt D$ will be
denoted as $\n(x)$. We will indicate coordinates of points and
components of vectors by writing $X_t=(X^1_t, X^2_t)$,
$Y_t=(Y^1_t, Y^2_t)$, $B_t=(B^1_t, B^2_t)$, and $\n(x) =(\n_1(x),
\n_2(x))$, but this notation may refer to a coordinate system
specific to a proof and different from the usual one. The angle
between vectors ${\bf p}$ and $\bf r$ will be denoted $\angle({\bf
p}, {\bf r})$, with the convention that it takes values in
$[0,\pi]$. Recall that $\alpha(x,y)=\angle(\n(x),
\n(y))\land(\pi-\angle(\n(x), \n(y)))$, for $x,y\in\prt D$.

The area of $D$ and the length of its boundary will be denoted
$|D|$ and $|\prt D|$, resp.

The distribution of the solution $\{(X_t,Y_t), t\geq 0\}$ to
(1.1)-(1.2) will be denoted $\bP^{x,y}$ and the distribution of
$\{X_t, t\geq 0\}$ will be denoted $\bP^x$. We will suppress the
superscripts when no confusion may arise. We will denote the usual
Markov shift operator by $\theta_t$.

\bigskip

In the first of our lemmas, we will prove that for an arbitrarily
small $\eps_0>0$, for any two points $x_0, y_0\in\ol D$,
two synchronously coupled reflecting Brownian motions $X$ and $Y$
starting from $x_0$ and $y_0$ respectively will come within
$\eps$-distance from each other in finite time a.s.
 This claim is very similar to Lemma 3.3 of [BC] but
sufficiently different to make it impossible for us to use that
lemma in the present paper. Regrettably, we could not find a
shorter proof of this seemingly quite intuitive result.

We will write $\tau^+_\eps =
\tau^+(\eps) = \inf\{t>0: \dist(X_t,Y_t) \geq \eps\}$ and
$\tau^-_\eps = \tau^-(\eps) = \inf\{t>0: \dist(X_t,Y_t) \leq
\eps\}$.

We remark that the following lemma  holds for smooth domains in
any dimension.

\bigskip
\noindent{\bf Lemma 3.1}. {\sl Consider any $\eps_0 >0$, any
$x_0,y_0 \in \ol D$, and assume that $(X_0,Y_0) = (x_0,y_0)$. Then
$\tau^-(\eps_0) < \infty$ a.s. }

\bigskip
\noindent{\bf Proof}. The proof will consist of several steps. In
the first two steps, we will prove some properties of the
deterministic Skorokhod mapping.

\bigskip
\noindent{\it Step 1}.
Let $\gamma=(\gamma^1, \gamma^2) :\
[0,\infty)\to \R^2$ be a continuous function
with $\gamma (0)\in \ol D$
and finite variation on each bounded interval of $[0, \infty)$.
Let $\lfloor\gamma\rfloor_{s,t}$ denote the total variation of
$\gamma$ on $[s,t]$.
We will use analogous notation for other functions.
 By the results of [LS], there exists
a unique pair of continuous functions $\beta: [0,\infty)\to \ol D$
and $\eta: [0,\infty)\to \R^2$ with the following properties:

\item{(i)}
$\lfloor \eta\rfloor_{s,t} \leq \lfloor\gamma\rfloor_{s,t}$ for
every $0\leq s \leq t$,

\item{(ii)} $ \int_0^\infty
\bone_{\{\beta_s\in D\}} d\ell_s =0$,
where $\ell_t \df \lfloor \eta\rfloor_{0,t}$,

\item{(iii)} $\eta_t = \int_0^t
\n(\beta_s) d\ell_s$ for every $t\geq 0$, and

\item{(iv)} $\beta_t = \gamma_t +\eta_t$, for
all $t\geq 0$.

We will call $(\beta, \eta)$ the Skorokhod transform of $\gamma$;
sometimes we will call $\beta$ the Skorokhod transform and denote
$\beta$ by $\S(\gamma)$.

We will show that

\item{(1.a)} $\lfloor\beta\rfloor_{s, t} \leq
\lfloor\gamma\rfloor_{s,t}$ for all $t>s\geq 0$, and

\item{(1.b)} for $c_1\in(0,1)$,
$c_2,c_3\in(0,\infty)$, there exists $c_4>0$ such that if
$\beta_{t_1}^1 - \beta_0^1 \leq (1-c_1)(\gamma_{t_1}^1 -
\gamma_0^1 )$, $\gamma_{t_1}^1 - \gamma_0^1 >c_2$, and
$\lfloor\gamma\rfloor_{0,t_1}\leq c_3$ then
$\lfloor\gamma\rfloor_{0,t_1} - \lfloor\beta\rfloor_{0,t_1} \geq
c_4$.

\medskip
According to (8') of [LS],
$$ dl_t = -1_{\partial D}(\beta_t) \< \n
(\beta_t), \, d \gamma_t \>
\eqno (3.1)
$$
and so
$$ d \beta_t= d\gamma_t + \n (\beta_t) dl_t
 = d\gamma_t - 1_{\partial D} (\beta_t) \n (\beta_t)
\< \n (\beta_t), d\gamma_t \>. \eqno (3.2)
$$
This proves that $\[\beta \]_{s, t} \leq \[\gamma \]_{s, t}$ for
any $0\leq s<t$,
i.e., this proves (1.a).

Let $\theta_t$ denote the angle between
 $\n (\beta_t)$ and $d\gamma_t$ whenever $\beta_t\in \partial D$
and $d\gamma_t$ is defined. Otherwise, define $\theta_t =\pi /2$.
Note that $\ell_t $ is non-decreasing, by its definition in (ii),
so (3.1) implies that $\theta_t \in [{\pi \over 2}, \, \pi]$.
Recall that $\n(x) =(\n_1(x), \n_2(x))$. By (3.2),
$$\eqalignno{ \int_0^{t_1} |\cos \theta_s|  |d\gamma_s|
&\geq  \int_0^{t_1} \n_1 (\beta_s) \cos \theta_s d \gamma_s \cr
&=(\gamma^1_{t_1}-\gamma^1_0)-(\beta^1_{t_1}-\beta^1_0) \cr &\geq
c_1 ( \gamma^1_{t_1}-\gamma^1_0) \geq c_1 c_2. \cr}
$$
Since $\int_0^{t_1} |d \gamma_s| = [\gamma ]_{0, t_1} \leq  c_3$,
for  $\delta :=\min\{1,\,  c_1c_2/(2c_3) \}$ we have from the above
$$ \int_0^{t_1} |\cos \theta_s| 1_{\{\theta_s\in [
{\pi \over 2}+\delta, \, \pi] \}}  |d\gamma_s|
\geq c_1c_2-c_3 \sin \delta \geq c_1c_2/2. \eqno(3.3)
$$
On the other hand, $|d\beta_t|= | d\gamma_t| \sin \theta_t $ and
so
$$\eqalignno{
 \[\gamma \]_{0, t_1}- \[\beta \]_{0, t_1}
&= \int_0^{t_1} (1- \sin \theta_s  ) |d \gamma_s| \cr
&\geq  {1\over 2} \int_0^{t_1} \cos^2 \theta_s |d\gamma_s| \cr
&\geq  {1\over 2} \int_0^{t_1} 1_{\{\theta_s\in [
{\pi \over 2}+\delta, \, \pi] \}} \cos^2 \theta_s |d\gamma_s| \cr
&\geq {\sin \delta \over 2}
 \int_0^{t_1} 1_{\{\theta_s\in [
{\pi \over 2}+\delta, \, \pi] \}} |\cos \theta_s|\,  |d\gamma_s| \cr
&\geq {c_1c_2 \sin \delta \over 4}:=c_4. \cr}
$$
This proves (1.b).

\bigskip

\noindent{\it Step 2}. Since $D$ is bounded and has a smooth
boundary, there is a constant $c_1< \infty$ such that any two
points $x,y\in \ol D$ can be connected by a $C^\infty$ curve
inside $\ol D$ of length $t_1=t_1(x,y) < c_1$. Consider any
$x,y\in \ol D$, and fix some $C^\infty$ curve $\gamma: [0,t_1] \to
\ol D$ with the natural (length) parametrization, and such that
$t_1<c_1$, $\gamma_0 = x$ and $\gamma_{t_1} = y$.

In this step, we will extend the definition of $\gamma$ from $[0,
t_1]$ to $[0, \infty)$. We will show that for any $D$ and
$\eps>0$, there exists a constant $c_2\in[c_1, \infty)$ such that
any curve $\gamma$ defined initially on $[0,t_1]$ may be extended
to $[0,\infty)$ in such a way that for some $t\leq c_2$,
$$ |\gamma_t - \S (\gamma + y-x)_t | \leq \eps. \eqno(3.4)
$$
Recall that $\S( \gamma)$ is the Skorokhod transform of $ \gamma$
(see Step 1).

Let $\{\beta_t, 0 \leq t \leq t_1\}$ be the Skorokhod transform of
$\{\gamma_t + \gamma_{t_1} - \gamma_0, 0 \leq t \leq t_1\}$,
defined as in Step 1. We will inductively define $\gamma_t$ for
all $t\geq 0$. Let $\gamma_t = \beta_{t-t_1} $ for $ t\in [t_1, 2
t_1]$, and let $\{\beta_t, t\in [t_1, 2 t_1]\}$ be the Skorokhod
transform of $\{\gamma_t + \gamma_{2t_1} - \gamma_{t_1}, t\in
[t_1, 2 t_1]\}$. We continue by induction, i.e., we let $\gamma_t
= \beta_{t-kt_1}$ for $ t\in [kt_1, (k+1) t_1],, \ k\geq 2$, and
we let $\{\beta_t, t\in [kt_1, (k+1) t_1]\}$ be the Skorokhod
transform of $\{\gamma_t+ \gamma_{(k+1)t_1} - \gamma_{kt_1}, t\in
[kt_1, (k+1) t_1]\}$. Note that both $\gamma_t$ and $\beta_t$ stay
in $\ol D$ for all $t\geq 0$. Clearly
$$\beta = \S (\gamma +y-x) \qquad \hbox{ and } \qquad
\gamma_t=\beta_{t-t_1} \quad
\hbox{for every } t\geq t_1. \eqno(3.5)
$$
By (1.a) in Step 1 and (3.5),
 we have $\lfloor\gamma\rfloor_{s,t} \leq t-s$ for all
$0\leq s \leq t < \infty$.

If $|\gamma_0 - \beta_0 | \leq \eps$ then we are done. Otherwise,
at least one of the following inequalities holds, $|\gamma^1_0 -
\beta^1_0 | \geq \eps/2$ or  $|\gamma^2_0 - \beta^2_0 | \geq
\eps/2$. We will assume without loss of generality that it is the
first of the two inequalities that holds and we will make another
harmless assumption that in fact $\gamma^1_0 - \beta^1_0  \leq
-\eps/2$, or, equivalently, $\gamma^1_{t_1} - \gamma^1_0 \geq
\eps/2$. Let $c_3$ be the diameter of $D$. Fix some $c_4 \in
(0,1)$ and integer $j>1$ such that $\sum_{k=1}^j c_4^{k-1} \eps /2
> 2 c_3$.
If we had
$$ \gamma_{kt_1}^1 - \gamma_{(k-1)t_1}^1
 \geq c_4 (\gamma_{(k-1)t_1}^1 - \gamma_{(k-2)t_1}^1)
\quad \hbox{ for every } 2\leq k \leq j,
\eqno(3.6)
$$
then we would obtain
$$
 \gamma_{jt_1}^1 - \gamma_0^1
 = \sum_{k=1}^j \gamma_{kt_1}^1 -\gamma_{(k-1)t_1}^1
 \geq \sum_{k=1}^j c_4^{k-1} \eps /2 >2 c_3,
$$
and that would contradict the definition of $c_3$ as the diameter
of $D$.
 So there must be some $k_0\leq j$ such that
 $$ \gamma_{k_0 t_1}^1 - \gamma_{(k_0-1)t_1}^1
 \leq c_4 (\gamma_{(k_0-1)t_1}^1 - \gamma_{(k_0-2)t_1}^1).  \eqno(3.7)
$$
Let $k_0$ be the smallest integer
with this property. Then
$$ \gamma_{kt_1}^1 - \gamma_{(k-1)t_1}^1
 \geq c_4 (\gamma_{(k-1)t_1}^1 - \gamma_{(k-2)t_1}^1)$$
for all $k< k_0$ and so $\gamma_{(k_0-1)t_1}^1 -
\gamma_{(k_0-2)t_1}^1 \geq c_4^j \eps/2$. The following is
equivalent to (3.7),
 $$ \beta_{(k_0-1)t_1}^1 - \beta_{(k_0-2)t_1}^1
 \leq c_4 (\gamma_{(k_0-1)t_1}^1 - \gamma_{(k_0-2)t_1}^1).$$
By (1.b) of Step 1, for some $c_5>0$,
 $$\lfloor\gamma\rfloor_{(k_0-2)t_1,(k_0-1)t_1}
 - \lfloor\beta\rfloor_{(k_0-2)t_1,(k_0-1)t_1}
 \geq  c_5 .\eqno(3.8)$$

If $\lfloor\gamma\rfloor_{(k-1)t_1,kt_1} \leq \eps$  for some
$k\leq j+1$ then
$$|\gamma_{(k-1)t_1} - \beta_{(k-1)t_1}|
=  |\gamma_{(k-1)t_1} - \gamma_{kt_1}| \leq
\lfloor\gamma\rfloor_{(k-1)t_1,kt_1} \leq \eps,
$$
and  we can take  $c_2 = jt_1$, i.e., there exists $t\leq jt_1$
with $|\gamma_t -\beta_t|\leq \eps$.

If $\lfloor\gamma\rfloor_{(k-1)t_1,kt_1} > \eps$ for all $k\leq
j+1$, then by Step 1 and (3.8),
 $$\eqalign{
 \lfloor\gamma\rfloor_{0,jt_1} -  \lfloor\beta\rfloor_{0,jt_1}
 &=  \sum_{k=1}^j \lfloor\gamma\rfloor_{(k-1)t_1,kt_1} -
 \lfloor\beta\rfloor_{(k-1)t_1,kt_1}\cr
 &\geq\lfloor\gamma\rfloor_{(k_0-2)t_1,(k_0-1)t_1}
 - \lfloor\beta\rfloor_{(k_0-2)t_1,(k_0-1)t_1}
 \geq c_5 .
 }$$

Thus we
have shown that either there exists $0\leq t\leq jt_1$ with
$|\gamma_t -\beta_t|\leq \eps$ or
$\lfloor\gamma\rfloor_{t_1,(j+1)t_1} =
\lfloor\beta\rfloor_{0,jt_1} \leq \lfloor\gamma\rfloor_{0,jt_1} -
c_5 $. The same argument shows that either there exists
$t\in[kt_1, (k+j)t_1]$ with $|\gamma_t -\beta_t|\leq \eps$ or
 $$\lfloor\gamma\rfloor_{(k+1) t_1,(k+1+j)t_1} =
 \lfloor\beta\rfloor_{k t_1,(k+j)t_1} \leq \lfloor\gamma\rfloor_{k
 t_1,(k+j)t_1} - c_5 .\eqno(3.9)$$
Recall that $\lfloor\gamma\rfloor_{0,t_1} = t_1 \leq c_1$, and
$\lfloor\gamma\rfloor_{k t_1,(k+1)t_1} =
\lfloor\beta\rfloor_{(k-1) t_1,kt_1} \leq
\lfloor\gamma\rfloor_{(k-1) t_1,kt_1}$ for all $k$. Hence,
$\lfloor\gamma\rfloor_{0,jt_1} \leq j c_1$, and if (3.9) holds
for all $k\leq m$, then
 $$0\leq \lfloor\gamma\rfloor_{m t_1,(m+j)t_1}
 \leq jc_1 - m c_5 .$$
This can be true only if $m\leq jc_1 /c_5 $. Hence, for some $k
\leq jc_1 /c_5  + 1$ and some $t\in[0,(k+j)t_1]$ we have
$|\gamma_t -\beta_t|\leq \eps$.

\bigskip
\noindent{\it Step 3}. First, we will present a version of the
``support theorem'' stronger than that given in Theorem I (6.6) in
[Ba]. Recall that one calls a continuous non-decreasing function
$\psi:[0,\infty) \to [0,\infty)$ with $\psi(0)=0$ a modulus of
continuity for a function $\gamma: [0,t_1]\to\R^2$ if for all
$s,t\in[0,t_1]$ we have $|\gamma_t - \gamma_s| \leq \psi(|t-s|)$.
Let $\K_{\psi,t_1}$ denote the family of all functions $\gamma:
[0,t_1]\to\R^2$ with modulus of continuity $\psi$. Let $\bP$ denote
the Wiener measure on $C[0,t_1]^2$, i.e., the distribution of the
planar Brownian motion. It follows easily from the existence of
``L\'evy's modulus of continuity'' (see Theorem 2.9.25 in [KS]),
that for every $t_1\in(0,\infty)$ and $p_0<1$ there exists $\psi$
such that $\bP(\K_{\psi,t_1})>p_0$. This fact can be used to modify
the proof of Proposition I (6.5) of [Ba] to show that there exists
$\psi$ such that for any $\eps>0$ one can find $p_1>0$ with
 $$\bP(\{\gamma\in\K_{\psi,t_1}: \sup_{0\leq t\leq t_1} |\gamma_t|\leq
 \eps\}) \geq p_1.\eqno(3.10)$$
Let $\psi_\lambda(t) = \psi(t) + \lambda t$ and for
$\phi:[0,t_1]\to \R^2$, let
 $$\K_{\psi,t_1,\phi,\eps}=
 \{\gamma\in\K_{\psi,t_1}: \sup_{0\leq t\leq t_1}
 |\gamma_t-\phi_t|\leq \eps\}.$$
If $\gamma \in \K_{\psi,t_1}$ and $\phi$ is Lipschitz with
constant $\lambda $ then $\gamma + \phi \in
\K_{\psi_\lambda,t_1}$. The proof of Theorem I (6.6) in [Ba] can
be easily modified to yield the following version of the support
theorem. Suppose that $\eps,p_1>0$ and $\psi$ satisfy (3.10). Then
for every $\lambda,t_1 <\infty$ and $\eps' >0$ one can find
$p_2>0$ such that for any function $\phi:[0,t_1]\to \R^2$ which is
Lipschitz with constant $\lambda$ and satisfies $\phi(0)=0$, we
have $\bP(\K_{\psi_\lambda,t_1,\phi,\eps'}) \geq p_2$. The important
aspect of the last assertion is that $p_2$ does not depend on
$\phi$. Fix a function $\psi$ satisfying this statement for the
rest of the proof.

Recall that ${\cal S} (\gamma)$ denotes the Skorokhod transform of
$\gamma$ (see Step 1) and let $\eps_0$ be as in the statement of
the lemma. Let $c_1$ be the constant defined in Step 2 and
$c_2>c_1$ be the constant in Step 2 relative to $\eps_0/5$ in
place of $\eps$. By Theorem 1.1 in [LS], the Skorokhod mapping
${\cal S}: C([0,c_2],\R^2) \to C([0,c_2],\R^2)$ is H\"older
continuous on compact sets. Let $\K = \{\gamma\in \K_{\psi_2,
c_2}: \gamma_0 \in \ol D\}$. The set $\K$ is compact so one can
find $\eps_1\in(0,\eps_0/5)$ such that if $\gamma, \wh
\gamma\in\K$ and $|\wh \gamma_t - \gamma_t| \leq \eps_1$ for $0
\leq t \leq c_2$, then $|{\cal S}(\wh \gamma)_t - {\cal
S}(\gamma)_t| \leq \eps_0/5$ for $0 \leq t \leq c_2$.

Recall from Step 2 that for every pair of points $x,y\in \ol D$
there is a curve $\gamma = \gamma^{x,y}: [0,\infty) \to \ol D$
such that $\gamma_0=x$, $\gamma_{t_1}=y$ for some $0< t_1 < c_1$,
and $\lfloor\gamma\rfloor_{s,t}\leq t -s$, for all $s$ and $t$ in
$[0, t_1]$. By Step 2, we can extend $\gamma$ to be a curve in
$\ol D$ satisfying (3.4) and (3.5). Note that $\gamma$ is a
Lipschitz curve on $[0, \infty)$ with Lipschitz constant 1.

Recall that reflected Brownian motions $ X_t$ and $Y_t$ are
defined in (1.1)-(1.2) relative to a Brownian motion $B_t$ and
assume that $X_0=x$ and $Y_0=y$. Find $p_2>0$ such that
$\bP(\K_{\psi_2,c_2,\phi,\eps_1}) \geq p_2$ for every Lipschitz
function $\phi$ with Lipschitz constant $1$ satisfying
$\phi(0)=0$. It follows that
$$\bP\left(\{B_t+ x ,
 0\leq t \leq c_2\} \in \K_{\psi_2,c_2,  \gamma^{x,y},\eps_1} \right)
 \geq p_2.\eqno(3.11)$$
Consider $\omega$ such that $B_\cdot(\omega)+x
\in \K_{\psi_2,c_2, \gamma^{x,y},\eps_1}\subset
\K_{\psi_2,c_2}$.
Then
$$ |{\cal S}(B_\cdot+x)_t - {\cal S}(\gamma^{x, y})_t|
\leq \eps_0/5 \qquad \hbox{for every } 0\leq t \leq c_2.
$$
Clearly $B_\cdot +y\in \K_{\psi_2, c_2}$. Since
$|B_t + y - (\gamma^{x,y}_t+y-x)|\leq \eps_1$ for
$t\in [0, c_2]$,
we have
$$ |{\cal S}(B_\cdot+y)_t - {\cal S}(\gamma^{x,y}+y-x)_t|\leq
\eps_0/5 \qquad \hbox{for } 0\leq t \leq c_2.
$$

Note that by Step 2 and our choice of $c_2$, there is some $t_0\in
[0, c_2]$ such that $|\gamma^{x,y}_{t_0} - {\cal
S}(\gamma^{x,y}+y-x)_{t_0}| \leq \eps_0/5$, for some $ t_0\in[0,
c_2]$. Note also that since $\gamma [0, \infty) \subset \ol D$, by
the uniqueness of the Skorokhod problem, $S(\gamma)=\gamma$.
Combining these observations, we conclude that
$$\eqalign{
 |&X_{t_0} - Y_{t_0}|
 = |{\cal S}(B_\cdot + x)_{t_0} - {\cal S}(B_\cdot + y)_{t_0}|\cr
 &= |{\cal S}(B_\cdot+x)_{t_0} - \gamma_{t_0}|
     +| \gamma_{t_0}-{\cal S} (\gamma +y-x)_{t_0}|
     +| {\cal S} (\gamma +y-x)_{t_0} -{\cal S}( B_\cdot +y)_{t_0}| \cr
 &\leq  \eps_0/5 + \eps_0/5+ \eps_0/5 <
 \eps_0.
 }$$
It follows from (3.11) that there exists $p_2>0$ such that for any
$x,y\in\ol D$, $X_0=x$, $Y_0=y$, the probability that there exists
$t_0 \leq c_2$ with $\dist(X_{t_0} , Y_{t_0})\leq \eps_0$ is
greater than $p_2$. By the Markov property applied at times
$jc_2$, $j=1,2, ...$, the probability that there is no $t_0\leq
kc_2$ with $\dist(X_{t_0} , Y_{t_0})\leq \eps_0$ is bounded above
by $(1-p_2)^k$. This implies easily that with probability one,
there exists $t< \infty$ with $\dist(X_t, Y_t) \leq  \eps _0$.
\qed

\bigskip
\noindent{\bf Lemma 3.2}. {\sl Let $\delta (x)$ denote the
Euclidean distance between $x$ and $\partial D$. Define
$\tau_D=\inf\{t\geq 0: X_t \notin D\}$  and $\tau_{\B(x,
r)}=\inf\{t\geq 0: X_t \notin \B(x, r)\}$. Then there exists
$c_1<\infty$ such that for $X_0=x_0\in D$,
$$\bP( \tau_{\B(x_0, r)} \leq \tau_D ) \leq c_1 \delta (x_0)/r
\qquad \hbox{ for } r\geq \delta (x_0).
\eqno (3.12)
$$
}

\bigskip
\noindent{\bf Proof}.
We are going to prove that (3.12) holds for any bounded
$C^{1,1}$-smooth domain in $\R^n$, for any $n\geq 2$.

Since $D$ is a bounded $C^{1,1}$-smooth domain, the ``uniform''
boundary Harnack principle holds for $D$ (see [A]), that is,
there exist $r_0>0$ and $c>0$ such that for $z\in \partial D$,
$r\in (0, r_0]$ and any non-negative harmonic functions $u$ and
$v$ in $D\cap \B(z, 2r)$ that vanish continuously on $\partial
D\cap \B(z, 2r)$, we have
$$ {u(x)\over v(x)} \leq c \, {u(y)\over v(y)}
\qquad \hbox{for any } x, y \in D\cap \B(z, r).
$$

Let $\{K_D(x, z); \, x\in D, \, z\in \partial D\}$
denote the Poisson kernel of the Brownian motion
$W$ killed upon leaving $D$; that is,
$$ {\Bbb E}^x \left[ \phi (W_{\tau_D} )\right]
=\int_{\partial D} K_D(x, z) \phi (z) \sigma (dz)
$$
for every continuous function $\phi$ on $\partial D$,
where $\sigma$ denotes the surface area measure.
Since $D$ is bounded $C^{1,1}$-smooth, it is known (see [Z]) that
there are constants $c_3>c_2>0$ such that
$$ {c_2 \delta (x) \over |x-z|^n} \leq K_D (x, z)
 \leq {c_3 \delta (x) \over |x-z|^n}
\qquad \hbox{for every } x\in D \hbox{ and } z\in \partial D.
\eqno (3.13)
$$
Note that $x\mapsto K_D(x, z)$ is a harmonic function in $D$
and vanishes continuously on $\partial D \setminus \{z\}$.

The lemma clearly holds when $\delta (x_0) \geq r_0/8$. This is
because, since $D$ is bounded, there is $R>0$ such that for every
$x_0 \in D$, $D\subset \B(x_0, R)$ and so  (3.12) holds trivially
for $r>R$. Thus in the case of $\delta (x_0) \geq r_0/8$, (3.12)
holds for every $r>0$ by choosing $c_1$ sufficiently large.

We now assume $\delta (x_0)<r_0/8$. Without loss of generality, we
may and do assume that $r> 8 \delta (x_0)$ and $\B(x_0, 2r)^c \cap
D \not= \emptyset$. We can further assume that $r\leq r_0$ since
$D$ is bounded.

Define $h(x)=\bP_x ( \tau_{\B(x_0, r)} \leq \tau_D )$. Clearly,
$h$ is a harmonic function in $D\cap \B(x_0, r)$ and vanishes
continuously on $\partial D \cap \B(x_0, r)$. Let $y_0\in \partial
D$ be such that $\delta (x_0)={\rm dist}(x_0, y_0)$. By the
triangle inequality, $\B(y_0, r/2)\subset \B(x_0, r)$. Now take
$z\in D\setminus \B(x_0, 2 r)$. Since $x_0\in \B(y_0, r/4)$, we
have by the boundary Harnack inequality,
$$ {h(x_0)\over K_D (x_0, z)} \leq c {h(x)\over K_D(x,z)}
  \qquad \hbox{for every } x\in \B(y_0, r/4).
$$
Let $x=x_0+{r\over 8}(x_0-y_0)$. Note that $h(x)\leq 1$,
$$  {1\over 2}|z-x_0|\leq |z-x|\leq 2 |z-x_0|,
$$
and $\delta (x) \geq r/8$.
These facts and (3.13) imply that $h(x_0) \leq c_1 \delta
(x_0)/r$. This proves the lemma. \qed

\bigskip
We fix parameters $a_1, a_2>0$ for the rest of the paper. We will
impose bounds on their values later on. Let $S_0=U_0=0$ and for
$k\geq 1$ define
 $$\eqalign{
 S_k& = \inf\{t> U_{k-1}: \dist(X_t, \prt D)
 \lor \dist(Y_t, \prt D) \leq a_2 \dist(X_t, Y_t)^2\},  \cr
 U_k &= \inf\{t> S_k:
 \dist(X_t, X_{S_k})  \lor\dist(Y_t, Y_{S_k}) \geq a_1
 \dist(X_{S_k}, Y_{S_k}) \}. \cr}
$$

We will assume that $a_1<1/4$. Then it is easy to see that
$\bP(U_k < \infty\mid S_k < \infty)=1$, for every $k$. Finiteness
of $S_k$'s is less obvious. The next lemma contains a result that
is significantly stronger than the finiteness of $S_k$'s. This
stronger result is needed in later arguments.

\bigskip
\noindent{\bf Lemma 3.3}. {\sl There exist
$c_1,c_2,c_3,c_4\in(0,\infty)$ and $\eps_0,r_0,p_0>0$ with the
following properties. Assume that $X_0\in \prt D$, $\dist(X_0,Y_0)
= \eps$, $\dist(Y_0,\prt D) = r$ and let
 $$T_1 = \inf\{t\geq 0: \dist(X _t, X_0) \lor \dist(Y_t, Y_0) \geq c_1
 r\}.$$

(i) If $\eps\leq\eps_0$ and $r\leq r_0$ then $\bP( S_1 \leq T_1,
L^X_{S_1} - L^X_{0} \leq c_2 r ) \geq p_0$.

(ii) If $\eps\leq\eps_0$ and $r\leq c_3 \eps$ then $\bE( L^X_{S_1
\land \tau^+(\eps_0)} - L^X_{0}) \leq c_4 r $.

 }

\bigskip
\noindent{\bf Proof}. (i) Recall the notation from the beginning
of this section. Let $CS_1$ be the orthonormal coordinate system
such that $X_0=0$ and $\n(X_0) $ lies on the second axis. Assume
that $ r_0<\eps_0< 1/(200\nu^*)$. Let $c_5\in(0,1/6)$ be a small
constant whose value will be chosen later. The following
definitions refer to the coordinates in $CS_1$,
 $$\eqalign{
 T_2 &= \inf\{t\geq 0: Y^2_t \geq 2r\},\cr
 T_3&=\inf\{t\geq 0: |Y^1_t-Y^1_{0}| \geq c_5 r\},\cr
 T_4&=\inf\{t\geq 0: Y_t \in \prt D\},\cr
 A_1& = \{T_4 \leq T_2 \land T_3\},\cr
 T_5 &=\inf\{t\geq 0: |X^1_t-X^1_{0}| \geq 2 c_5 r\}.
 }$$

First we will assume that $r\leq \eps/2$. We will show that $T_5
\geq T_2 \land T_3 \land T_4$ if $A_1$ holds. We will argue by
contradiction. Assume that $A_1$ holds and $T_5 < T_2 \land T_3
\land T_4$. Then $B^1_t-B^1_{0} = Y^1_t-Y^1_{0}$ for $t\in[0,T_5]$
so $|B^1_t-B^1_{0}| \leq c_5 r$ for the same range of $t$'s. We
have
  $$X^1_{T_5} - X^1_{0} = B^1_{T_5} - B^1_{0} +
  \int_{0}^{T_5} \n_1(X_t) dL^X_t,
  $$
so $\left|\int_{0}^{T_5} \n_1(X_t) dL^X_t\right| \geq c_5 r$. We
assume that $\eps_0>0$ is so small that for $r< \eps_0$ and
$x\in\B(0,2c_5 r)$, we have $\n_2(x) \geq |\n_1(x)|/(2\nu^* \cdot
2 c_5 r)$. It follows that
  $$\int_{0}^{T_5} \n_2(X_t) dL^X_t
  \geq c_5 r/(2\nu^* \cdot 2 c_5 r)= 1/(4\nu^*).$$
Note that $B^2_t-B^2_{0} = Y^2_t-Y^2_{0}$ for
$t\in[0,T_4]$. Since $\dist(X_0,Y_0) = \eps$, we have
$$\dist(X_0, Y_t)\leq \eps+\sqrt{(2 r)^2 +(c_5r)^2} < \eps +3r < 3\eps
\qquad \hbox{ for } t\leq T_2\land T_3\land T_4.
$$
Therefore for $t\leq T_2\land T_3\land T_4$,
 $|Y^2_t|\leq 3\eps$.
Since $T_5 < T_2\land T_3\land T_4$, it follows that
$$ |B^2_t-B^2_s|=|Y^2_t-Y^2_s| \leq |Y^2_t|+|Y^2_s|
\leq 6\eps
\qquad \hbox{ for } s,t\in[0,T_5].
$$
Thus
  $$\eqalignno{
  X^2_{T_5} - X^2_{0} &\geq -| B^2_{T_5} - B^2_{0}| +
  \int_{0}^{T_5} \n_2(X_t) dL^X_t \cr
&\geq -6\eps + 1/(4\nu^*) \geq - 6\eps_0 + 1/(4\nu^*)  \geq 44\eps_0 ,
 }$$
and $X^2_{T_5} \geq 44\eps_0 + X^2_{0}= 44\eps_0$. Let $T_6
=\sup\{t\leq T_5: X_t\in \prt D\}$. Then $B^2_{T_5}-B^2_{T_6} =
X^2_{T_5}-X^2_{T_6}\geq 44\eps_0 - r\geq 43\eps_0$, a
contradiction with the fact that $|B^2_t-B^2_s| \leq 6\eps\leq
6\eps_0$ for $s,t\in[0,T_5]$. This proves that $T_5 \geq T_2 \land
T_3 \land T_4$ if $A_1$ holds.

We will show that if $A_1$ holds then $S_{1} \leq T_4$. Assume
that $A_1$ holds and let $T_7 = \sup\{t\leq T_4: X_t \in \prt
D\}$. Note that neither $X_t$ nor $Y_t$ visit $\prt D$ on the
interval $(T_7,T_4)$. Hence, $X_{T_7} - Y_{T_7} = X_{T_4} -
Y_{T_4}$. If $\eps_0$ and $r_0$ are sufficiently small then
$|X^1_0 - Y^1_0| \geq \eps/2 $ because $r \leq \eps/2$ and
$\dist(Y_0,\prt D) =r$. We have assumed that $A_1$ holds so $
|Y^1_{T_4}- Y^1_{0}| \leq c_5 r $. We have proved that $T_5 \geq
T_4$ on $A_1$, so $ |X^1_{T_4}- X^1_{0}| \leq 2c_5 r $. Recall
that $c_5\leq 1/6$ and $r\leq \eps/2$. It follows that
 $$\eqalign{
 \dist(X_{T_7},Y_{T_7}) &=\dist(X_{T_4},Y_{T_4})
 \geq |X^1_{T_4} - Y^1_{T_4}| \cr
 &\geq |X^1_0 - Y^1_0| - |Y^1_{T_4}- Y^1_{0}| - |X^1_{T_4}- X^1_{0}|
 \geq \eps/2 - 3c_5 r \geq \eps/4.
 }$$
On the other hand, assuming $\eps_0>0$ is small,
$$\eqalign{
 \dist(X_{T_7},Y_{T_7}) &\leq \dist(X_{T_7},X_0) +
\dist(X_0,Y_{0}) + \dist(Y_0,Y_{T_7})\cr & \leq  2 |X^1_{T_7} -
X^1_0| +\eps + \dist(Y_0,Y_{T_7}) \leq 2 \cdot 2 c_5 r + \eps + 3r
\leq 3\eps.}$$
 We have
$$|Y^1_{T_4}- Y^1_{T_7}| =|Y^1_{T_4}- Y^1_{0}| +
|Y^1_{0}- Y^1_{T_7}| \leq c_5 r + c_5 r =2 c_5 r.
$$
Since $Y_{T_4} \in \prt D$, $X_{T_7} \in \prt D$, $X_{T_7} -
Y_{T_7} = X_{T_4} - Y_{T_4}$, $\dist(X_{T_7},Y_{T_7}) \leq 3\eps$,
and $|Y^1_{T_4}- Y^1_{T_7}| \leq 2 c_5 r$, we have
$\angle(\n(Y_{T_4}), \n(X_{T_7})) \leq 2 \nu^* \cdot 2(3\eps + 2
c_5 r)\leq 16 \nu^* \eps$. This and easy geometry show that
$\dist(Y_{T_7}, \prt D) \leq 2\cdot 2 c_5 r \cdot 16 \nu^*\eps= 64
c_5 r \nu^* \eps$. Hence,
 $${\dist(Y_{T_7}, \prt D) \over\dist(X_{T_7}, Y_{T_7})}
 \leq {64 c_5 r \nu^* \eps\over \eps/4}
 = 256c_5 r\nu^*
 \leq 128 c_5 \nu^* \eps
 \leq 32 c_5 \nu^* \dist(X_{T_7}, Y_{T_7}).$$
We choose $c_5>0$ so small that $32 c_5 \nu^* \leq a_2$. Then
$\dist(Y_{T_7}, \prt D) \leq a_2 \dist(X_{T_7}, Y_{T_7})^2$. We
obviously have $\dist(X_{T_7}, \prt D) \leq a_2 \dist(X_{T_7},
Y_{T_7})^2$ because $X_{T_7} \in \prt D$. This shows that $S_1
\leq T_7$ and completes the proof that if $A_1$ holds then $S_{1}
\leq T_4$.

Assume that $A_1$ holds and suppose that $\int_{0}^{T_4} \n_2(X_t)
dL^X_t \geq 20r$. We will show that this leads to a contradiction.
Recall that $|B^2_t-B^2_s| \leq 8 r$ for $s,t\in[0,T_4]$. We
obtain
 $$
 X^2_{T_4}=X^2_{T_4} - X^2_{0} \geq -| B^2_{T_4} - B^2_{0}| +
 \int_{0}^{T_4} \n_2(X_t) dL^X_t
 \geq -8r + 20 r = 12 r.
 $$
Recall that $T_7 =\sup\{t\leq T_4: X_t\in \prt D\}$. We have
$B^2_{T_4}-B^2_{T_7} = X^2_{T_4}-X^2_{T_7}\geq 12r - r = 11r$, a
contradiction with the fact that $|B^2_t-B^2_s| \leq 8 r$ for
$s,t\in[0,T_4]$. Hence, if $A_1$ holds then $\int_{0}^{T_4}
\n_2(X_t) dL^X_t \leq 20r$. Note that $\n_2(x) \geq 1/2$ for all
$x\in \prt D \cap \B(0, 6r)$, assuming that $\eps_0>0$ is small
and $r\leq r_0< \eps_0$. We have shown that if $A_1$ holds then
$T_5 \geq T_4$, so $\n_2(X_t) \geq 1/2$ for $t\in[0,T_4]$ such
that $X_t\in\prt D$. This implies that,
 $$ (1/2)(L^X_{S_{1}} - L^X_{0})
 \leq (1/2)(L^X_{T_4} - L^X_{0}) \leq \int_{0}^{T_4} \n_2(X_t) dL^X_t
 \leq 20 r .$$
We have shown that $\{S_1 \leq T_1, L^X_{S_1} - L^X_{0} \leq 40
r\} \subset A_1$. It is easy to see that $\bP(A_1)> p_1$ for some
$p_1>0$ which depends only on $c_5$. This completes the proof of
part (i) in the case $r\leq \eps/2$, with $c_1 = 2$ and $c_2 =40$.

Next consider the case when $r\geq \eps/2$. Let
 $$\eqalign{T_8&=\inf\{t> 0: \dist(Y_t, X_{0}) \geq
 2\eps\},\cr
 T_{9}&= \inf\{t>0: X_t\in\prt D, \dist(Y_t,\prt D) \leq
 \dist(X_t,Y_t)/2\},\cr
 T_{10}&=\inf\{t>0: L^X_t - L^X_{0} \geq 20\eps\},\cr
 A_2 &= \{T_4 \leq T_8\},\cr
 A_3 &= \{T_{9} \leq T_8\land T_{10}\}.
 }$$
We will show that $A_2 \subset A_3$. Assume that $A_2$ holds.
First, we will prove that $L^X_{T_4} - L^X_{0} \leq 20\eps$.
Suppose otherwise, i.e., $L^X_{T_4} - L^X_{0} \geq 20\eps$. Recall
that we are using the coordinate system $CS_1$ with the origin at
$X_{0}\in\prt D$. Let $T_{11} =\inf\{t\geq 0: |X^1_t-X^1_{0}| \geq
5\eps\}$. We will show that $T_{11} \geq T_4$. We will argue by
contradiction. Assume that $T_{11} < T_4$. We have assumed that
$A_2$ holds, so $T_{11} < T_8$. Then $B^1_t-B^1_{0} =
Y^1_t-Y^1_{0}$ for $t\in[0,T_{11}]$ and $|B^1_t-B^1_{0}| \leq
4\eps$ for the same range of $t$'s. We have
 $$\left|\int_{0}^{T_{11}} \n_1(X_t) dL^X_t\right|
 = |X^1_{T_{11}} - X^1_{0} -( B^1_{T_{11}} - B^1_{0})|
 \geq |X^1_{T_{11}} - X^1_{0}| -| B^1_{T_{11}} - B^1_{0}|
 \geq 5\eps - 4\eps = \eps.
 $$
If $\eps_0>0$ is sufficiently small and $\eps<\eps_0$ then
$\n_2(x) \geq |\n_1(x)|/(2\nu^* \cdot 5\eps)$ for $x\in \prt D
\cap \B(0, 5\eps)$, so $\int_{0}^{T_{11}} \n_2(X_t) dL^X_t \geq
\eps/(2\nu^* \cdot 5\eps)= 1/(10\nu^*)$. We have $B^2_t-B^2_{0} =
Y^2_t-Y^2_{0}$ for $t\in[0,T_{11}]$, because $T_{11} \leq T_4$, so
$|B^2_t-B^2_s| \leq 4\eps$ for $s,t\in[0,T_{11}]$. Recall that
$\eps<\eps_0< 1/(100\nu^*)$. We obtain,
 $$\eqalignno{
 X^2_{T_{11}} &=X^2_{T_{11}}- X^2_{0}
 \geq -| B^2_{T_{11}} - B^2_{0}| +
 \int_{0}^{T_{11}} \n_2(X_t) dL^X_t \cr
 &\geq -4\eps + 1/(10\nu^*)  \geq 6\eps.
 }$$
Let $T_{12} =\sup\{t\leq T_{11}: X_t\in \prt D\}$. Then
$B^2_{T_{11}}-B^2_{T_{12}} = X^2_{T_{11}}-X^2_{T_{12}}\geq 6\eps -
\eps = 5\eps$, a contradiction, because $|B^2_t-B^2_s| \leq 4\eps$
for $s,t\in[0,T_{11}]$. This proves that $T_{11} \geq T_4$.

Recall that we have assumed that $L^X_{T_4} - L^X_{0} \geq
20\eps$. We have $\n_2(x) \geq 1/2$ for $x\in\prt D\cap \B(0,
10\eps)$, assuming $\eps_0>0$ is small and $\eps\leq \eps_0$.
Since $T_{11} \geq T_4$, $\n_2(X_t) \geq 1/2$ for $t \leq T_4$
such that $X_t\in \prt D$, so
 $$\eqalignno{
 X^2_{T_4} &=X^2_{T_4}- X^2_{0} \geq -| B^2_{T_4} - B^2_{0}| +
 \int_{0}^{T_4} \n_2(X_t) dL^X_t
 \geq -4\eps + (1/2) (L^X_{T_4} - L^X_{0})\cr
 &\geq  - 4\eps + 10\eps = 6\eps.
 }$$
Recall that $T_7 =\sup\{t\leq T_4: X_t\in \prt D\}$. Then
$B^2_{T_4}-B^2_{T_7} = X^2_{T_4}-X^2_{T_7}\geq 6\eps - \eps =
5\eps$, a contradiction, because $|B^2_t-B^2_s| \leq 4\eps$ for
$s,t\in[0,T_{11}]$. This proves that if $A_2$ holds then
$L^X_{T_4} - L^X_{0} \leq 20\eps\leq 40r$.

Note that $X_{T_4} - Y_{T_4} = X_{T_7} - Y_{T_7}$, $Y_{T_4},
X_{T_7} \in \prt D$, and $T_7 \leq T_4 \leq T_{11}$. Assuming that
$\eps_0>0$ is small, these facts easily imply that the angle
between $X_{T_7} - Y_{T_7}$ and the tangent line to $\prt D$ at
$X_{T_7}$ is smaller than $\pi/8$, so $\dist(Y_{T_7}, \prt D) \leq
\dist (X_{T_7}, Y_{T_7})/2$. Hence, $T_{9} \leq T_4$ and,
therefore, if $A_2$ occurs then $T_{9} \leq T_4 \leq T_8 \land
T_{10}$. This completes the proof that $A_2 \subset A_3$.

It is easy to see that $\bP(A_2)>p_2>0$, where $p_2$ depends only on
$\eps_0$ and $D$. It follows that $\bP(A_3)>p_2$.

We may now apply the strong Markov property at the stopping time
$T_9$ and repeat the argument given in the first part of the
proof, discussing the case $r\leq \eps/2$. It is straightforward
to complete the proof of part (i), adjusting the values of
$c_1,c_2,\eps_0,r_0$ and $p_0$, if necessary.

(ii) Let $c_1$ and $c_2$ be as in part (i) of the lemma, let
$T_5^0 = 0$, and for $k\geq 1$ let
 $$\eqalign{
 T^k_1 &=\inf\{t\geq T_5^{k-1}: \dist(X _{T_5^{k-1}}, X_t)
 \lor \dist(Y_{T_5^{k-1}}, Y_t) \geq c_1
  \dist(Y_{T_5^{k-1}},\prt D)\},\cr
 T_2^k &= \inf\{t\geq T_5^{k-1}: L^X_{t} - L^X_{T_5^{k-1}} \geq c_2
 \dist(Y_{T_5^{k-1}},\prt D) \},\cr
 T_3^k&=\inf\{t\geq T_5^{k-1}: Y_t \in \prt D\},\cr
 T_4^k &= T_1^k \land T_2^k \land T_3^k , \cr
 T_5^k &= \inf\{t\geq T_4^{k}: X_t \in \prt D\}.\cr
 }$$
Let $\eps_0>0$ be the constant which works for part (i) of the
lemma. An examination of the proof of part (i) shows that we have
in fact proved a statement stronger than that in part (i) of the
lemma, namely, using the notation of the first part of the proof,
 $$\bP( S_1 \leq T_1\land T_4,L^X_{S_1} - L^X_{0}
 \leq c_2 r ) \geq p_0.\eqno(3.14)$$

Next we will estimate $\bE\dist(Y_{T_5^k\land \tau^+(\eps_0)},\prt
D)$. By Lemma 3.2,
 $$\eqalignno{
 \bP(\sup_{t\in[T_4^k, T_5^k]}
 \dist(X_t, X_{T_4^k}) \in [2^{-j-1},2^{-j}] \mid {\cal F}_{T^k_4})
 &\leq c_6 \dist(X_{T_4^k},\prt D)/2^{-j}\cr
 &\leq c_7 \dist(Y_{T_5^{k-1}},\prt D)/2^{-j} .&(3.15)
 }$$
Write $\gamma=\dist(Y_{T_5^{k-1}},\prt D)$, and let $j_0$ be the
largest integer such that $2^{-j_0} \geq \diam(D)$. Consider $j$
such that $\sup_{t\in[T_4^k, T_5^k\land \tau^+(\eps_0)]}
\dist(X_t, X_{T_4^k}) \leq 2^{-j} $. It is not hard to show that
if $j_0\leq j\leq |\log \eps_0|$ then $\dist(Y_{T_5^k\land
\tau^+(\eps_0)},\prt D) \leq c_7 \eps_0 2^{-j}$ for some
$c_7<\infty$. If $j\geq |\log \eps_0| $ then $\dist(Y_{T_5^k\land
\tau^+(\eps_0)},\prt D) \leq \gamma + c_8 \eps_0 2^{-j}$. This and
(3.15) imply that
 $$\eqalign{\bE&(\dist(Y_{T_5^k\land \tau^+(\eps_0)},\prt D)
 \mid {\cal F}_{T^k_4}) \cr
 &\leq
 \sum_{j_0\leq j\leq |\log \eps_0|} c_7 \eps_0 2^{-j}
 \bP(\sup_{t\in[T_4^k, T_5^k]}
 \dist(X_t, X_{T_4^k}) \in [2^{-j-1},2^{-j}] \mid {\cal F}_{T^k_4})\cr
 & \quad + \sum_{ | \log \eps_0|\leq j \leq |\log \gamma|}
 (\gamma + c_8 \eps_0 2^{-j}) \bP(\sup_{t\in[T_4^k, T_5^k]}
 \dist(X_t, X_{T_4^k}) \in [2^{-j-1},2^{-j}] \mid {\cal F}_{T^k_4})\cr
 & \quad + \sum_{ j \geq | \log \gamma|}
 (\gamma + c_8 \eps_0 2^{-j}) \bP(\sup_{t\in[T_4^k, T_5^k]}
 \dist(X_t, X_{T_4^k}) \in [2^{-j-1},2^{-j}] \mid {\cal F}_{T^k_4})\cr
 &\leq \sum_{j_0\leq j\leq |\log \eps_0|}
 c_9 \eps_0 2^{-j} ( \gamma/ 2^{-j}) \cr
 &\quad + \gamma + \sum_{|\log \eps_0| \leq j \leq
 | \log \gamma|}  c_{10} \eps_0 2^{-j}(\gamma/2^{-j})
 + \sum_{ j \geq
 | \log \gamma|}  c_{10} \eps_0 2^{-j}\cr
 &\leq c_{11} \eps_0 \gamma |\log \eps_0| +
 \gamma + c_{12} \gamma\eps_0 |\log \gamma| + c_{13}\gamma\eps_0
 \leq \gamma(1+ c_{13}\eps_0 |\log \eps_0|).
 }$$
 Thus
 $$\eqalign{
 \bE &(\dist(Y_{T_5^k\land \tau^+(\eps_0)},
 \prt D)\mid {\cal F}_{T^k_4})
 \bone_{\{ T_5^{k-1} < \tau^+(\eps_0) \}}\cr
 &\leq  \bone_{\{ T_5^{k-1} < \tau^+(\eps_0) \}}
 (1+ c_{13} \eps_0 |\log \eps_0|)
 \dist(Y_{T_5^{k-1}\land \tau^+(\eps_0)},\prt D).}$$
 This, (3.14) and the strong Markov property yield,
 $$\eqalign{
 \bE( &\dist(Y_{T_5^k\land \tau^+(\eps_0)},\prt D)
 \bone _{\{S_1 \geq T_4^{k}\}}
 \bone_{\{ T_5^{k-1} < \tau^+(\eps_0) \}})\cr
 &=  \bE( \bone _{\{S_1 \geq T_4^{k}\}}
 \bone_{\{ T_5^{k-1} < \tau^+(\eps_0) \}}
 \bE (\dist(Y_{T_5^k\land \tau^+(\eps_0)},\prt D)
 \mid {\cal F}_{T^k_4}))\cr
 &\leq (1+ c_{13} \eps_0 |\log \eps_0|)
 \bE( \bone _{\{S_1 \geq
 T_4^{k}\}} \bone_{\{ T_5^{k-1} < \tau^+(\eps_0) \}}
 \dist(Y_{T_5^{k-1}\land \tau^+(\eps_0)},\prt D))\cr
 &\leq (1+ c_{13} \eps_0 |\log \eps_0|)
 \bE(  \dist(Y_{T_5^{k-1}\land \tau^+(\eps_0)},\prt D)
  \bone_{\{ T_5^{k-1} < \tau^+(\eps_0) \}}
 \bE(\bone _{\{S_1 \geq T_4^{k}\}} \mid {\cal F}_{T^{k-1}_5}))\cr
 &\leq (1+ c_{13} \eps_0 |\log \eps_0|)
 \bE (\dist(Y_{T_5^{k-1}\land \tau^+(\eps_0)},\prt D)
 \bone _{\{S_1 \geq T_4^{k-1}\}}
 \bone_{\{ T_5^{k-1} < \tau^+(\eps_0) \}} \cr
 &\quad \times (1- \bP(T_5^{k-1}\leq S_1 \leq T_4^{k}
 \mid S_1 \geq T_5^{k-1}))\cr
 &\leq (1+ c_{13} \eps_0 |\log \eps_0|)(1-p_0)
 \bE (\dist(Y_{T_5^{k-1}\land \tau^+(\eps_0)},\prt D)
 \bone _{\{S_1 \geq T_4^{k-1}\}}
 \bone_{\{ T_5^{k-1} < \tau^+(\eps_0) \}}) \cr
 &\leq (1+ c_{13} \eps_0 |\log \eps_0|)(1-p_0)
 \bE (\dist(Y_{T_5^{k-1}\land \tau^+(\eps_0)},\prt D)
 \bone _{\{S_1 \geq T_4^{k-1}\}}
 \bone_{\{ T_5^{k-2} < \tau^+(\eps_0) \}})
 .}$$
We obtain by induction,
 $$\eqalign{
 \bE( &\dist(Y_{T_5^k\land \tau^+(\eps_0)},\prt D)
 \bone _{\{S_1 \geq T_4^{k}\}}
 \bone_{\{ T_5^{k-1} < \tau^+(\eps_0) \}})\cr
 &\leq (1+ c_{13} \eps_0 |\log \eps_0|)^k(1-p_0)^k
 \bE \dist(Y_{T_5^{0}},\prt D)
 .}$$
It follows that
 $$\eqalign{
 \bE&( L^X_{S_1 \land \tau^+(\eps_0)} - L^X_{0})
 = \sum_{k=0}^\infty \bE\left(( L^X_{S_1 \land \tau^+(\eps_0)} - L^X_{0})
 \bone_{\{ S_1 \in [T_5^k, T_5^{k+1})\}} \right)\cr
 & = \sum_{k=0}^\infty \bE\left(
 \bone_{\{ S_1 \in [T_5^k, T_5^{k+1})\}}
 \sum_{j=0}^ {k}
 \bone_{\{ T_5^{j} < \tau^+(\eps_0) \}}
 (L^X_{T_5^{j+1}\land \tau^+(\eps_0)}
 - L^X_{T_5^j\land \tau^+(\eps_0)})\right )\cr
 & \leq \sum_{k=0}^\infty \bE\left(
 \bone_{\{ S_1 \in [T_5^k, T_5^{k+1})\}}
 \sum_{j=0}^ {k}  \bone_{\{ T_5^{j-1} < \tau^+(\eps_0) \}}
 c_2 \dist(Y_{T_5^j\land \tau^+(\eps_0)} , \prt D) \right )\cr
 & \leq \sum_{k=0}^\infty c_2
 \bE (\dist(Y_{T_5^k\land \tau^+(\eps_0)},\prt D)
 \bone _{\{S_1 \geq T_4^{k}\}}
  \bone_{\{ T_5^{k-1} < \tau^+(\eps_0) \}}) \cr
 & \leq \sum_{k=0}^\infty c_2
 (1+ c_{13} \eps_0 |\log \eps_0|)^k(1-p_0)^k
 \bE \dist(Y_{T_5^{0}},\prt D)
 . }$$
If we assume that $\eps_0>0$ is sufficiently small, this is
bounded by $c_{14}\bE \dist(Y_{T_5^{0}},\prt D)=c_{14}
\dist(Y_{0},\prt D)$.
 \qed
\bigskip

Recall that $a_1$ and $a_2$ are parameters in the definitions of
$S_k$'s and $U_k$'s stated at the paragraph
preceding Lemma 3.3.

 \bigskip
\noindent{\bf Corollary 3.4}. {\sl For any $a_1,a_2 >0$, and any
starting points $(X_0,Y_0) = (x_0,y_0) \in \ol D\times \ol D$, all
stopping times $S_k$ are finite a.s. }

\bigskip
\noindent{\bf Proof}. Let $\eps_0$ and $r_0$ be as in Lemma 3.3.
By Lemma 3.1, there is a finite stopping time $T_1$ such that
$\dist(X_{T_1}, Y_{T_1}) \leq \eps_0 \land r_0$. So there exist
$p_1>0$ and $c_1<\infty$ such that $\bP(T_1\leq c_1)\geq p_1>0$.
Let $T_2$ be the first time after $T_1$ when either $X$ or $Y$
hits $\prt D$, and note that $\dist(X_{T_2}, Y_{T_2}) =
\dist(X_{T_1}, Y_{T_1}) \leq \eps_0 \land r_0$. Let
$$T_3 = \inf\{t\geq T_2: \dist(X _t, X_{T_2})
 \lor \dist(Y_t, Y_{T_2}) \geq c_2 r_0\} .
$$
It is easy to see that $\bP(T_2\leq T_3<\infty\mid T_1<\infty)=1$
when $c_2>0$ is small.
Select such $c_2>0$ and apply the strong Markov property at $T_2$
and Lemma 3.3 (i) to see that there exists $p_2>0$ such
 that $\bP(S_1 \leq T_3\mid T_1<\infty) \geq p_2$.
On the other hand, for some $c_3< \infty$, we have $\bP(T_3 \leq
T_1 + c_3 \mid T_1 <\infty) \geq 1-p_2/2$. It follows that
$\bP(S_1 \leq T_3\leq c_1 +c_3)> p_1p_2/2$ and so $\bP(S_1 > c_1
+c_3)< 1-{p_1p_2\over 2}$. By the Markov property, $\bP(S_1
> k(c_1+c_3)) \leq (1-{p_1p_2\over 2})^k$ for $k\geq 1$, so $S_1 <
\infty$, a.s.

Recall that $\bP(U_k< \infty\mid S_k<\infty)=1$ for every $k$,
according to the remark made before the statement of Lemma 3.3. By
induction and the strong Markov property applied at $S_k$'s and
$U_k$'s, all stopping times $S_k$ and $U_k$ are finite a.s.
 \qed

\bigskip
\noindent{\bf Lemma 3.5}. {\sl For any $c_1>0$, one can choose
$a_1,a_2>0$ and $\eps_0>0$ so that for every $k\geq 1$ and all
$s,t\in(S_k\land \tau^+(\eps_0), U_k\land \tau^+(\eps_0))$, a.s.,
 $$\angle(X_t-Y_t,X_s-Y_s)\leq c_1 \dist(X_{S_k\land
 \tau^+(\eps_0)}, Y_{S_k\land \tau^+(\eps_0)}).$$
 }

 \bigskip
\noindent{\bf Proof}. Recall that $D$ is assumed to be
$C^4$-smooth. Elementary geometry shows that for any $c_1>0$ there
exist $\eps_0,a_1,a_2>0$ with the following properties. Suppose
that $x,y \in \prt D$ and $r=\dist(x,y) \leq \eps_0/2$. Let
 $$\eqalignno {
 A_1 &= \B(x,2 a_1 r)\cap \prt D,&(3.16)\cr
 A_2 &= \B(y,2 a_1 r)\cap \prt D,\cr
 A_3 &= \{z\in \ol D: \dist(x,z) \leq 2 a_1r,
 \dist(z, \prt D)\leq 2a_2 r^2\},\cr
 A_4 &= \{z\in \ol D: \dist(y,z) \leq 2 a_1r,
 \dist(z, \prt D)\leq 2a_2 r^2\},\cr
 A_5 & = \{z\in \ol D: \exists v\in A_1, u,w\in A_2: z= v-u +
 w\} , \cr
 A_6 & = \{z\in \ol D: \exists v\in A_2, u,w\in A_1: z= v-u +
 w\} , \cr
 A_7 &= A_3 \cup A_5 ,\cr
 A_8 &= A_4 \cup A_6 ,\cr
 \beta &= \sup \{\angle(x_0 - y_0, x_1-y_1): x_0,x_1\in A_7,
 y_0,y_1 \in A_8\}.
 }$$
With a suitable choice of small $\eps_0,a_1,a_2>0$, we have $\beta
\leq (c_1/4) r$.

We will assume that $S_k < \tau^+(\eps_0)$ because otherwise
$(S_k\land \tau^+(\eps_0), U_k\land \tau^+(\eps_0))=\emptyset$ and
there is nothing to prove. Let $x\in \prt D$ be the closest point
to $X_{S_k}$ and let $y\in \prt D$ be the closest point to
$Y_{S_k}$. Note that if $\eps_0$ and $a_2$ are small then
$\dist(x,y) = r \leq 2 \dist (X_{S_k}, Y_{S_k})$. We use points
$x$ and $y$ to define sets $A_k$, as in (3.16). We will argue that
$\angle(X_t-Y_t,X_s-Y_s)\leq 2\beta$, for all $s,t\in(S_k,
U_k\land \tau^+(\eps_0))$. Note that $X_{S_k} \in A_3$ and
$Y_{S_k} \in A_4$. This and the definition of $\beta$ imply that
$\angle(X_{S_k}- Y_{S_k}, x-y) \leq \beta$.

Suppose that there exists $t\in(S_k, U_k\land \tau^+(\eps_0))$
with $\angle(X_{t}- Y_{t}, x-y) > \beta$ and let $T=\inf\{t\geq
S_k: \angle(X_{t}- Y_{t}, x-y) > \beta\}$. By continuity,
$\angle(X_{T}- Y_{T}, x-y) = \beta$.
 It is impossible that both
$X_T$ and $Y_T$ are in $D$, because then we would have $X_t - Y_t
= X_T - Y_T$ for some $t_0 >0$ and all $t\in (T-t_0, T+ t_0)$.
This would imply that $\angle(X_{t}- Y_{t}, x-y) = \beta$ for
$t\in (T-t_0, T+ t_0)$, and, therefore, $\inf\{t\geq S_k:
\angle(X_{t}- Y_{t}, x-y)> \beta\} \geq T+t_0$, a contradiction.
We will show that it cannot happen that $X_T,Y_T\in \prt D$.
Suppose that it is true that $X_T,Y_T\in \prt D$ and recall that
we are working under assumption that $T < U_k$. The definition of
$U_k$ implies that $X_T \in A_3 \subset A_7$ and $Y_T \in A_4
\subset A_8$. Since $\angle(X_{T}- Y_{T}, x-y) = \beta$, it
follows that the supremum in the definition of $\beta$ is attained
for points $x_0, x_1, y_0, y_1 \in \prt D$ (take $x_0=x, y_0=y,
x_1 = X_T $ and $y_1 = Y_T$). Easy geometry shows that this cannot
be the case because we can slightly move either $y_0$ or $y_1$
into the interior of $D$ to increase the value of $\angle(x_0 -
y_0, x_1-y_1)$.

Suppose without loss of generality that $X_T\in \prt D$ and
$Y_T\in D$. For some random $t_1>0$, the process $Y$ will not
touch the boundary within $[T,T+t_1]$, while the local time $L^X$
will have a non-zero increment, a.s. It is easy to see that the
local-time-term push that $X$ will get over $[T, T+t_1]$ will make
$\angle(X_{t}- Y_{t}, x-y)$ smaller, and hence $\angle(X_{t}-
Y_{t}, x-y)\leq \beta$ for $t \in [T, T+t_1]$, contradicting the
definition of $T$. We conclude that $\angle(X_{t}- Y_{t}, x-y)
\leq \beta$ for $t\in(S_k, U_k\land \tau^+(\eps_0))$. This and the
fact that $\angle(X_{S_k}- Y_{S_k}, x-y) \leq \beta$ imply that
 $$\angle(X_t-Y_t,X_s-Y_s)\leq 2\beta\leq (c_1/2)r
 \leq c_1 \dist(X_{S_k}, Y_{S_k}),$$
for all $s,t\in(S_k, U_k\land \tau^+(\eps_0))$.
 \qed

\bigskip

Recall the definition of the stopping times $S_k$, $U_k$ from the
paragraph preceding Lemma 3.3. For $k\geq 1$,  define
 $$\eqalign{
  \rho_t &= {\dist(X_t,Y_t)\over \dist(X_0,Y_0)},\cr
 \wt \rho_t &=
 \prod_{j=1}^\infty {\dist(X_{U_j\land t}, Y_{U_j\land t})
 \over  \dist(X_{S_j\land t},Y_{S_j\land t}) }, \cr
 \ol \rho_t &= \prod_{j=0}^\infty {\dist(X_{S_{j+1}\land t},
 Y_{S_{j+1}\land t}) \over \dist(X_{U_j\land t}, Y_{U_j\land t})
 },
 }$$
with the convention $0/0=1$. Note that $\rho_t = \wt \rho_t \ol
\rho_t$. Let ${\cal T} = \bigcup_{k\geq 1} (S_k,U_k]$ and ${\cal
T}^c = (0,\infty) \setminus {\cal T}$.

\bigskip
\noindent{\bf Lemma 3.6}. {\sl For any $c_1 >0$ there exist
$a_0,\eps_0>0$ such that if $a_1, a_2\in(0,a_0)$ and $\dist(X_0,
Y_0)  \leq \eps_0$ then for all $t \geq 0$, a.s.,
$$
\left| \log \wt \rho_{t \land \tau(\eps_0)} + {1\over 2}
\int_{[0,t\land \tau(\eps_0)] \cap {\cal T}} \left( \nu(X_s) dL^X_s  +
\nu(Y_s) dL^Y_s \right) \right| \leq c_1 \left( L^X_{t\land \tau(\eps_0)} +
L^Y_{t\land \tau(\eps_0)} \right).
$$
}

 \bigskip
\noindent{\bf Proof}. Since $D$ is assumed to be $C^4$-smooth, for
any $c_2\in (0, 1)$ and $c_3>0$, we can find $\eps_0>0$ so small
that for any $x,y\in \prt D$ with $\dist(x,y) \leq 2 \eps_0$,
 $$\big(-(1+c_2)(1/2) \nu(x)-c_3/2 \big) \dist(x, y)
 \leq {x-y \over \dist(x,y)} \cdot \n(x)
 \leq \big( -(1-c_2)(1/2) \nu(x)+c_3/2 \big) \dist(x, y). $$
This, Lemma 3.5, differentiability of $\nu$ and simple geometry
show that one can choose small $a_1,a_2>0$ and $\eps_0>0$ so that
for every $k\geq 1$ and all $t\in[S_k, U_k\land \tau^+(\eps_0)]$
such that $X_t \in \prt D$, assuming $S_k < \tau^+(\eps_0)$,
 $$\eqalign{
 (-(1+c_2)(1/2) \nu(X_t)-c_3) \dist(X_{S_k}, Y_{S_k})
 &\leq {X_t - Y_t\over \dist(X_t , Y_t)}
 \cdot\n(X_t)\cr
 &\leq (-(1-c_2)(1/2) \nu(X_t)+c_3)\dist(X_{S_k}, Y_{S_k}).
 }$$
Analogous estimates hold for $ {Y_t - X_t\over \dist(Y_t, X_t)}
\cdot\n(Y_t)$. We obtain for $t \in [S_k, U_k\land
\tau^+(\eps_0)]$,
 $$\eqalign{
 &\dist (X_t, Y_t) -  \dist(X_{S_k}, Y_{S_k})
 = \int_{S_k}^t {(X_s - Y_s)\over \dist(X_s , Y_s)}
 \cdot\n(X_s) dL^X_s
 + \int_{S_k}^t {(Y_s - X_s)\over \dist(Y_s, X_s)}
 \cdot\n(Y_s) dL^Y_s\cr
 &\leq \int_{S_k}^t (-(1-c_2)(1/2) \nu(X_s)+c_3)
 \dist(X_{S_k}, Y_{S_k}) dL^X_s\cr
 &\quad+ \int_{S_k}^t (-(1-c_2)(1/2) \nu(Y_s)+c_3)
 \dist(X_{S_k}, Y_{S_k}) dL^Y_s.
 }$$
Thus
 $$
 { \dist (X_t, Y_t) \over  \dist(X_{S_k}, Y_{S_k}) } \leq  1 -
(1-c_2)(1/2) \int_{S_k}^t  (\nu(X_s)
 dL^X_s +\nu(Y_s) dL^Y_s) + c_3
 \int_{S_k}^t  (dL^X_s + dL^Y_s).$$
We obtain in a similar way,
 $$
 { \dist (X_t, Y_t) \over  \dist(X_{S_k}, Y_{S_k}) } \geq  1 -
(1+c_2)(1/2) \int_{S_k}^t  (\nu(X_s)
 dL^X_s +\nu(Y_s) dL^Y_s) - c_3
 \int_{S_k}^t  (
 dL^X_s + dL^Y_s).$$
Note that $ \dist (X_t, Y_t) / \dist(X_{S_k}, Y_{S_k}) \in
[1-2a_1, 1+2a_1]$ for $t \in [S_k, U_k\land \tau^+(\eps_0)]$, so
 $$\eqalign{
 &{ \dist (X_t, Y_t) \over  \dist(X_{S_k}, Y_{S_k}) } \cr
 &\leq
 (1+2a_1) \land \left(1 -
 (1-c_2)(1/2) \int_{S_k}^t  (\nu(X_s)
 dL^X_s +\nu(Y_s) dL^Y_s) + c_3
 \int_{S_k}^t  (dL^X_s + dL^Y_s)\right),}$$
and
 $$\eqalign{
 &{ \dist (X_t, Y_t) \over  \dist(X_{S_k}, Y_{S_k}) }\cr
&\geq (1 -2a_1) \lor \left(  1 -
(1+c_2)(1/2) \int_{S_k}^t  (\nu(X_s)
 dL^X_s +\nu(Y_s) dL^Y_s) - c_3
 \int_{S_k}^t  (
 dL^X_s + dL^Y_s)\right).}$$
We have $ 1+ a \leq e^a$ for all $a$. For any $c_4>0$ we can
choose $a_1>0$ so small that $1+a \geq e^a e^{-c_4|a|}$ for
$a\in[-2a_1, 2a_1]$. Hence, for sufficiently small $a_1$, and $t
\in [S_k, U_k\land \tau^+(\eps_0)]$, assuming $S_k <
\tau^+(\eps_0)$,
 $$\eqalign{
 \wt \rho_t  &= {\dist(X_t, Y_t) \over  \dist(X_{S_k},
 Y_{S_k}) }
 \prod_{j=1}^{k-1}
 {\dist(X_{U_j}, Y_{U_j}) \over  \dist(X_{S_j}, Y_{S_j}) } \cr
 &\leq \left(1 -  (1/2-c_2/2) \int_{S_k}^t  (\nu(X_s)
 dL^X_s +\nu(Y_s) dL^Y_s) + c_3
 \int_{S_k}^t  (dL^X_s + dL^Y_s)\right) \cr
 &\qquad \times \prod_{j=1}^{k-1}\left( 1 -  (1/2-c_2/2)
 \int_{S_j}^{U_j}  (\nu(X_s)
 dL^X_s +\nu(Y_s) dL^Y_s)  + c_3
 \int_{S_j}^{U_j}  (dL^X_s + dL^Y_s)\right) \cr
 & \leq \exp\Bigg( -  (1/2-c_2/2)
 \int_{[0,t\land \tau(\eps_0)] \cap {\cal T}}(\nu(X_s)
 dL^X_s +\nu(Y_s) dL^Y_s) \cr
 &\qquad + c_3
 \int_{[0,t\land \tau(\eps_0)] \cap {\cal T}}
 (dL^X_s + dL^Y_s) \Bigg),
 }$$
and
 $$\eqalign{
 \wt \rho_t
 \geq \exp\Bigg( &-  (1/2+c_2)
 \int_{[0,t\land \tau(\eps_0)] \cap {\cal T}}(\nu(X_s)
 dL^X_s +\nu(Y_s) dL^Y_s) \cr
 &- 2c_3
 \int_{[0,t\land \tau(\eps_0)] \cap {\cal T}}
 (dL^X_s + dL^Y_s)\cr
 &-c_4  \Big((1/2+c_2)
 \int_{[0,t\land \tau(\eps_0)] \cap {\cal T}}(|\nu(X_s)|
 dL^X_s +|\nu(Y_s)| dL^Y_s) \cr
 &\qquad + 2c_3
 \int_{[0,t\land \tau(\eps_0)] \cap {\cal T}}
 (dL^X_s + dL^Y_s) \Big )
 \Bigg).
 }$$
Since $|\nu|$ is bounded and $c_2,c_3$ and $c_4$ are arbitrarily
small, the last two estimates yield the lemma.
 \qed
\bigskip

Recall that we have assumed that for every $x\in \prt D$, there
are only finitely many points $y\in \prt D$ with $\alpha(x,y) =0$.

\bigskip

\noindent{\bf Lemma 3.7}. {\sl Suppose that $z\in \prt D$ and let
$K=\{y\in\prt D: \alpha(z,y) =0\}$ and $M_k=\{ y\in \prt D:
\alpha(z,y ) \in [2^{-k}, 2^{-k+1}]\}$. There exist
$k_0,c_1<\infty$ and $c_2>0$ not depending on $z$ such that for
$k\geq k_0$, the arc length measure of $M_k$ is less than $c_1
2^{-k/2}$ and the distance from $M_k$ to $K$ is bounded below by
$c_2 2^{-k}$.
 }

\bigskip
\noindent{\bf Proof}. We have assumed that the boundary of $D$ is
$C^4$-smooth and that there exist at most a finite number of
points $x_1, x_2, \dots, x_n$ such that $\nu(x_k)=0$,
$k=1,\dots,n$. Moreover, we have assumed that the third derivative
of the function representing the boundary does not vanish at any
$x_k$. This implies that there exist $\delta_0,c_3,c_4>0$ such
that if $x\in \prt D$ and $\dist(x, x_k) \leq \delta_0$ for some
$k$ then $|\nu(x)| \geq c_3 \dist (x,x_k)$; moreover, if $x\in
\prt D$ and $\dist(x, x_k) \geq \delta_0$ for every $k=1,\dots,n$,
then $|\nu(x)| \geq c_4$. We make $\delta_0$ smaller, if
necessary, so that $\dist(x_j,x_k) \geq 4 \delta_0$ for all $j\ne
k$. It is elementary to see that there exists $c_5>0$ with the
following properties (i)-(iii).

(i) For every point $x\in \prt D$ such that $\dist(x, x_k) \geq
2\delta_0$ for every $k=1,\dots,n$, and every $y\in \prt D$ with
$\dist(x,y) \leq \delta_0$, we have $|\alpha(x,y)| \geq c_5 \dist
(x,y)$.

(ii) If $x\in \prt D$ and $\dist(x, x_k) < 2 \delta_0$ for some
$k$, $y\in \prt D$, $\dist(x,y) \leq \delta_0$, and $y$ lies on
the same side of $x_k$ as $x$ then $|\alpha(x,y)| \geq c_5
\dist(x,y) \dist(x,x_k)$.

(iii) If $x=x_k$ for some $k$, $y\in \prt D$ and $\dist(x,y) \leq
\delta_0$ then $|\alpha(x,y)| \geq c_5 \dist(x,y)^2$.

Make $\delta_0$ smaller, if necessary, so that for any $x,y\in
\prt D$ with $|\alpha(x,y)| \geq \pi/4$, we have $\dist(x,y) \geq
4\delta_0$.

Consider any $z\in \prt D$ and let $z_1, z_2, \dots, z_m$ be all
points in $\prt D$ such that $\alpha(z,z_k) =0$ or $\alpha(z,z_k)
=\pi/2$. The number $m$ of such points is bounded by a constant
$m_0$ depending on $D$ but not on $z$. The family of points $\{
z_1,\dots, z_m, x_1, \dots, x_n\}$ divides $\prt D$ into $n+m$
Jordan arcs $\Gamma_k$, $k=1,\dots, n+m$. Let $\lambda$ denote the
arc length measure on $\prt D$, i.e., $\lambda(dx)$ is an
alternative notation for $dx$.

Fix some $\Gamma_k$ and note that the curvature $\nu(x)$ has a
constant sign on this arc because there are no $x_j$'s between the
endpoints of $\Gamma_k$. Since there are no points $z_j$ between
the endpoints of $\Gamma_k$, the function $x\to \alpha(z,x)$ is
monotone on this arc. For an arc $\Gamma_k$, let $y_k^-$ and
$y_k^+$ denote its endpoints and assume that $\alpha(z,x)$ takes
the maximum on $\Gamma_k$ at $x=y_k^-$. It is elementary to deduce
from (i)-(iii) that for some $c_6>0 $ depending only on $D$, and
all $j$,
$$\lambda(\{x\in\Gamma_k: \alpha(x,y_k^-)
 \in [2^{-j}, 2^{-j-1}]\}) \leq c_{6} 2^{-j/2}. $$
Since the number of $\Gamma_k$'s is bounded by a constant
independent of $z$,
 $$\lambda(\{x\in\prt D: \alpha(x,z)
 \in [2^{-j}, 2^{-j-1}]\}) \leq c_{7} 2^{-j/2}. $$

Conditions (i)-(iii) easily imply that $\dist(M_j, K)\geq c_8
2^{-j}$ for some $c_8$ depending only on $D$.
 \qed

\bigskip
\noindent{\bf Lemma 3.8}. {\sl There exists $c_1<\infty$ such that
for any $s>0$,
 $$\rho_s \leq \exp(c_1 (L^X_s + L^Y_s)).$$
 }

\bigskip
\noindent{\bf Proof}. Since $D$ is assumed to be $C^4$-smooth,
there exists $c_2<\infty$ such that for any $x\in \prt D$ and
$y\in D$,
 $$ {x-y \over \dist(x,y)} \cdot \n(x)
 \leq c_2 \dist(x, y). \eqno(3.17)$$
Let $T_0 = 0$, and for $k\geq 1$,
 $$\eqalign{
 T_k &= \inf \{t\geq T_{k-1}: \dist (X_t, Y_t) \notin
 (\thalf \dist(X_{T_{k-1}}, Y_{T_{k-1}}) , \,
 2\dist(X_{T_{k-1}}, Y_{T_{k-1}})) \}\cr
 & \quad\land \inf \{t\geq T_{k-1}: L^X_t - L^X_{T_{k-1}} \geq 1 \}
 \land \inf \{t\geq T_{k-1}: L^Y_t - L^Y_{T_{k-1}} \geq 1 \}.
 }$$
Then, by (3.17), for any $k\geq 1$ and $t\in (T_{k-1}, T_k]$,
 $$\eqalign{
 &\dist (X_t, Y_t) -  \dist(X_{T_{k-1}}, Y_{T_{k-1}})\cr
 &= \int_{T_{k-1}}^t {(X_s - Y_s)\over \dist(X_s , Y_s)}
 \cdot\n(X_s) dL^X_s
 + \int_{T_{k-1}}^t {(Y_s - X_s)\over \dist(Y_s, X_s)}
 \cdot\n(Y_s) dL^Y_s\cr
 &\leq \int_{T_{k-1}}^{T_k} c_2 \dist(X_s , Y_s)dL^X_s
 + \int_{T_{k-1}}^{T_k} c_2 \dist(X_s , Y_s) dL^Y_s\cr
 &\leq 2c_2 \dist(X_{T_{k-1}} , Y_{T_{k-1}})
 \int_{T_{k-1}}^{T_k} (dL^X_s+ dL^Y_s).\cr
 }$$
This implies that for any $t\in (T_{k-1}, T_k]$,
 $$\eqalign{
 {\dist (X_t, Y_t ) \over \dist(X_{T_{k-1}}, Y_{T_{k-1}})}
 &\leq 1+ 2 c_2 (L^X_{T_k} - L^X_{T_{k-1}} + L^Y_{T_k} -
 L^Y_{T_{k-1}})\cr
 &\leq \exp (2 c_2 (L^X_{T_k} - L^X_{T_{k-1}} + L^Y_{T_k} -
 L^Y_{T_{k-1}})),
 }$$
and
 $$\eqalign{
 \dist (X_t, Y_t)
 &= {\dist (X_t, Y_t) \over \dist(X_{T_{k-1}}, Y_{T_{k-1}}) } \, \prod_{j=1}^{k-1}
{\dist (X_{T_j}, Y_{T_j}) \over \dist(X_{T_{j-1}}, Y_{T_{j-1}})} \cr
 & \leq \prod_{j=1}^k
 \exp (2 c_2 (L^X_{T_j} - L^X_{T_{j-1}} + L^Y_{T_j} -
 L^Y_{T_{j-1}}))\cr
 &\leq \exp (2c_2  (L^X_{T_k} + L^Y_{T_k})).
 }$$
This proves the lemma.
 \qed

\bigskip

We define a partial order for two distinct points $x=(x_1,x_2)$
and $y=(y_1,y_2)$ by saying that $x\prec y$ if $x_1 < y_1$, or
$x_1 = y_1$ and $x_2 < y_2$. Let $Z_t$ be the closest point in
$\prt D$ to the pair $\{X_t, Y_t\}$, if there is only one such
point. In the case when there are multiple points in $\prt D$ with
the minimum distance to $\{X_t, Y_t\}$, we let $Z_t$ be the point
which is the smallest one according to $\prec$; an easy argument
based on compactness of $\prt D$ shows that there exists such a
point. Our choice of the tie-breaking convention is arbitrary---it
plays no role in the proofs. Note that if $T_0=\inf\{t\geq 0: X_t
\in \prt D
 \hbox{ or  } Y_t \in \prt D\}$ then
$Z_{T_0} = X_{T_0}$ if $X_{T_0} \in \prt D$ and $Y_{T_0} \notin
\prt D$; $Z_{T_0} = Y_{T_0}$ if $Y_{T_0} \in \prt D$ and $X_{T_0}
\notin \prt D$; $Z_{T_0}$ can be either $X_{T_0}$ or $Y_{T_0}$ if
both  $X_{T_0} \in \prt D$ and $Y_{T_0} \in \prt D$.

The following piece of notation will be used in many lemmas,
 $$F(s,u, x, a) = \left\{\sup_{s\leq t \leq u} \dist(X_t,x) \leq a \right\}.
$$
The proof of the next lemma is the most complicated and delicate
argument in this paper.

\bigskip

\noindent{\bf Lemma 3.9}. {\sl Let $T_0=\inf\{t\geq 0: X_t \in
\prt D \hbox {  or  } Y_t\in \prt D\}$ and $\eps =
\dist(X_0,Y_0)$. There exist $\beta_0 \in (1/2,1)$ and
$c_1,c_2<\infty$ such that the following hold. Assume that,
 $$|\pi/2- \angle(X_0-Y_0, \n(Z_0))| \leq c_1 \dist(X_0,Y_0)^{\beta_0}
 \ \hbox {  and } \ \dist(X_0,\prt D)\leq c_2 \dist(X_0,Y_0).
 \eqno(3.18) $$

(i) There exist $c_3< \infty $ and $\eps_0>0$ such that whenever
$\eps \leq \eps_0$,
$$\bE |\log \dist(X_{S_1}, Y_{S_1}) - \log \dist (X_0, Y_0)|\leq
c_3 \eps.$$

(ii) For some $\beta_1 >0$, $\beta_2 >1$, $c_4< \infty $ and
$\eps_0>0$, we have for all $\eps \leq \eps_0$,
$$\bE \left( \bone_{F^c(T_0, S_1, Z(T_0),\eps^{\beta_1})}\,
|\log \dist(X_{S_1}, Y_{S_1}) - \log \dist (X_0, Y_0) |  \right)
\leq c_4 \eps^{\beta_2}.$$

(iii) Let $K = \{x\in \prt D: \tan \alpha(Z_0, x) \geq
\eps^{-\beta_1}\}$. For some $\beta_1 >0$, $\beta_2 >1$, $c_5<
\infty $ and $\eps_0>0$, we have for all $\eps \leq \eps_0$,
$$\bE \left( \bone_{\{Z_{T_0}\in K\}}\,
|\log \dist(X_{S_1}, Y_{S_1}) - \log \dist (X_0, Y_0) |
\right)\leq c_5 \eps^{\beta_2}.$$
 }

\bigskip
\noindent{\bf Proof}. (i) {\it Step 1}. For some $c_6<\infty$, let
 $$\eqalign{T_1 &= \inf\{t>T_0:
 \dist(X_t ,\prt D)\land \dist(Y_t ,\prt D)
 \leq \dist(X_t,Y_t), \cr
 & \qquad\qquad|\pi/2 -\angle(\n(Z_{t}), X_{t}-Y_{t})| \leq
 c_1\dist(X_t,Y_t)^{\beta_0}
  \}, \cr
 T_2 & = \inf\{t>T_0: \dist (X_t, X_{T_0})
 \geq c_6 \dist(X_{T_0}, Y_{T_0})^{\beta_0} \},\cr
 T_3& = \inf\{t> T_0: X_t \in \prt D\} \bone_{\{Y_{T_0} \in \prt
 D\}} + \inf\{t> T_0: Y_t \in \prt D\} \bone_{\{Y_{T_0} \notin \prt
 D\}},\cr
 T_4 & = \inf\{t > T_2: X_t \in \prt D \hbox{ or } Y_t \in \prt
 D\}, \cr
 T_5 & = \inf\{t>T_4: \dist (X_t, X_{T_4})
 \geq c_6 \dist(X_{T_4}, Y_{T_4})^{\beta_0} \},\cr
 T_6& = \inf\{t> T_4: X_t \in \prt D\} \bone_{\{Y_{T_4} \in \prt
 D\}} + \inf\{t> T_4: Y_t \in \prt D\} \bone_{\{Y_{T_4} \notin \prt
 D\}}
  . }$$
It is elementary to see that for a suitable choice of $c_6$,
$\{T_3 < T_2\} \subset \{T_1 < T_2\}$, and similarly $\{T_6 <T_5\}
\subset \{T_1 < T_5\}$.

We will now estimate changes in the distance between $X_t$ and
$Y_t$ over the interval $[T_0, T_1\land T_2]$ under various
scenarios, and probabilities of these scenarios.

Let $M_*=\{x\in\prt D: \alpha(x, Z_0) = \pi/2\}$. For integer $k$
and any $x\in M_*$, let $M_k= \{y \in \prt D: \tan \alpha(y,x )
\in [2^{-k}, 2^{-k+1})\}$. Let $ N $ be such that $ Z_{T_0} \in
M_N$. Let $k_1$ be the largest integer with $\dist(M_{k_1}, M_*)
\geq 4 c_6 \dist(X_0,Y_0)^{\beta_0}$. Since we are concerned with
the case when $\dist(X_0,Y_0)$ is small, we can assume that
$k_1>0$.

Suppose that $-k_1 \leq k \leq 0$. Then, by Lemmas 3.2 and 3.7,
$\bP(N=k) \leq c_7 \dist (X_0,Y_0) 2^{-k}$.
 If $N=k$ then $\dist(X_{T_0}, \prt D) \lor \dist(Y_{T_0}, \prt D)
\leq c_8 \dist (X_0,Y_0)2^k$. This and Lemma 3.2 imply that
 $$\bP(T_1 \geq T_2 \mid N=k)\leq \bP(T_3 \geq T_2\mid N=k)
 \leq c_9 \dist(X_{T_0}, Y_{T_0})^{1-\beta_0}2^k,
 $$
and, therefore,
 $$\bP(N=k, T_1 \geq T_2) \leq c_{10}
 \dist(X_{T_0}, Y_{T_0})^{2-\beta_0}.$$
Elementary geometry shows that the distance between $X_t$ and
$Y_t$ is reduced by at most a factor of $1- c_{11} 2^{2k}$ over
the interval $[T_0, T_1\land T_2]$, so we have $\dist(X_{T_1\land
T_2}, Y_{T_1\land T_2}) \geq (1- c_{11} 2^{2k})\dist (X_0,Y_0)$.

Next assume that $0 < k \leq k_1$. It follows from Lemmas 3.2 and
3.7 that $\bP(N=k) \leq c_{12} \dist (X_0,Y_0) 2^{-k/2}$. We
obviously have $\dist(X_{T_0}, \prt D) \lor \dist(Y_{T_0}, \prt D)
\leq \dist (X_0,Y_0)$. This and Lemma 3.2 imply that
 $$\bP(T_1 \geq T_2 \mid N=k)\leq \bP(T_3 \geq T_2\mid N=k)
 \leq c_{13} \dist(X_{T_0}, Y_{T_0})^{1-\beta_0},
 $$
and, therefore,
 $$\bP(N=k, T_1 \geq T_2) \leq c_{14} \dist(X_{T_0}, Y_{T_0})^{2-\beta_0}
 2^{-k/2}.$$
We have $\angle(\n(Z_t), X_t-Y_t) \in (c_{15} 2^{-k}, \pi - c_{15}
2^{-k})$ for $t \in [T_0, T_2]$. It follows that the distance
between $X_t$ and $Y_t$ is reduced by at most a factor of $c_{16}
2^{k}$ over the interval $[T_0, T_1\land T_2]$, so
$\dist(X_{T_1\land T_2}, Y_{T_1\land T_2}) \geq c_{16} 2^{-k}\dist
(X_0,Y_0)$.

The next case is $N \leq -k_1$. We trivially have $\bP(N\leq -k_1)
\leq 1$. If $N\leq -k_1$ then $\dist(X_{T_0}, \prt D) \lor
\dist(Y_{T_0}, \prt D) \leq  c_{17} \dist (X_0,Y_0)^{1+\beta_0}$.
Lemma 3.2 implies that
 $$\eqalignno{ \bP(N\leq -k_1, T_1 \geq T_2) &\leq
 \bP(T_1 \geq T_2 \mid N\leq -k_1)&(3.19)\cr
 &\leq \bP(T_3 \geq T_2\mid N\leq -k_1)
 \leq c_{18} \dist(X_{T_0}, Y_{T_0}).
 }
 $$
The distance between $X_t$ and $Y_t$ is reduced by at most a
factor of $1- c_{19} \dist(X_{T_0}, Y_{T_0})^{2\beta_0}$ over the
interval $[T_0, T_1\land T_2]$, so $\dist(X_{T_1\land T_2},
Y_{T_1\land T_2}) \geq (1- c_{19} \dist(X_{T_0},
Y_{T_0})^{2\beta_0})\dist (X_0,Y_0)$.

Let $N'$ be such that $X_{T_4} \in M_{N'}$ if $X_{T_4} \in \prt D$
and $Y_{T_4} \in M_{N'}$ if $Y_{T_4} \in \prt D$. We will analyze
the change to $\dist(X_t,Y_t)$ over the interval $[T_0, T_1\land
T_5]$ for different values of $N'$, assuming that $T_1 \geq T_2$
and $N\leq -k_1$.

Suppose that $-k_1 \leq k \leq 0$. Since $\dist(X_{T_2}, \prt D)
\leq c_{20} \dist(X_0, Y_0)^{\beta_0}$, Lemmas 3.2 and 3.7 imply
that $\bP(N'=k\mid N\leq -k_1, T_1 \geq T_2) \leq c_{21} \dist
(X_0,Y_0)^{\beta_0} 2^{-k}$. If $N'=k$ then $\dist(X_{T_4}, \prt
D) \lor \dist(Y_{T_4}, \prt D) \leq c_{22} \dist (X_0,Y_0)2^k$.
This and Lemma 3.2 imply that
 $$\eqalign{
 \bP(T_1 \geq T_5 \mid N\leq -k_1, T_1 \geq T_2, N'=k)
 &\leq \bP(T_6 \geq T_5\mid N\leq -k_1, T_1 \geq T_2, N'=k)\cr
 &\leq c_{23} \dist(X_{T_0}, Y_{T_0})^{1-\beta_0}2^k.}
 $$
We combine estimates of probabilities in this paragraph with
(3.19) to obtain,
 $$\bP(N\leq -k_1, N'=k, T_1 \geq T_5)
 \leq c_{24} \dist(X_{T_0}, Y_{T_0})^{2}.$$
The distance between $X_t$ and $Y_t$ is reduced by at most a
factor of $1- c_{25} 2^{2k}$ over the interval $[T_0, T_1\land
T_5]$, so $\dist(X_{T_1\land T_5}, Y_{T_1\land T_5}) \geq (1-
c_{25} 2^{2k})\dist (X_0,Y_0)$.

Next consider the case $0 < k \leq k_1$. By Lemmas 3.2 and 3.7,
 $$\bP(N'=k\mid N\leq -k_1, T_1 \geq T_2) \leq c_{26} \dist
 (X_0,Y_0)^{\beta_0} 2^{-k/2}.
 $$
If the event $\{ N\leq -k_1, T_1 \geq T_2\}$ holds then
$\dist(X_{T_4}, \prt D) \lor \dist(Y_{T_4}, \prt D) \leq \dist
(X_0,Y_0)$. This and Lemma 3.2 imply that
 $$\eqalign{
 \bP(T_1 \geq T_5 \mid N\leq -k_1, T_1 \geq T_2, N'=k)
 &\leq \bP(T_6 \geq T_5\mid N\leq -k_1, T_1 \geq T_2, N'=k)\cr
 &\leq c_{27} \dist(X_{T_0}, Y_{T_0})^{1-\beta_0},}
 $$
and, using (3.19),
 $$\bP(N\leq -k_1, N'=k, T_1 \geq T_5)
 \leq c_{28} \dist(X_{T_0}, Y_{T_0})^{2}
 2^{-k/2}.$$
We have $\angle(\n(Z_t), X_t-Y_t) \in (c_{29} 2^{-k}, \pi - c_{29}
2^{-k})$ for $t \in [T_0, T_5]$. It follows that the distance
between $X_t$ and $Y_t$ is reduced by at most a factor of $c_{30}
2^{k}$ over the interval $[T_0, T_1\land T_5]$, so
$\dist(X_{T_1\land T_5}, Y_{T_1\land T_5}) \geq c_{30} 2^{-k}\dist
(X_0,Y_0)$.

Consider the case $N'\leq -k_1 $. We trivially have $\bP(N'\leq
-k_1\mid N\leq -k_1, T_1 \geq T_2) \leq 1$. If $N'\leq -k_1$ then
$\dist(X_{T_4}, \prt D) \lor \dist(Y_{T_4}, \prt D) \leq c_{31}
\dist (X_0,Y_0)^{1+\beta_0}$. This and Lemma 3.2 imply that
 $$\eqalign{
 \bP(T_1 \geq T_5 \mid N\leq -k_1, T_1 \geq T_2, N'\leq -k_1)
 &\leq \bP(T_6 \geq T_5\mid N\leq -k_1, T_1 \geq T_2, N'\leq -k_1)\cr
 &\leq c_{32} \dist(X_{T_0}, Y_{T_0}),}
 $$
and, using (3.19),
 $$\bP(N\leq -k_1, N'\leq -k_1, T_1 \geq T_5)
 \leq c_{33} \dist(X_{T_0}, Y_{T_0})^{2}.$$
If $N'\leq -k_1 $, the distance between $X_t$ and $Y_t$ is reduced
by at most a factor of $1- c_{34} \dist(X_{T_0},
Y_{T_0})^{2\beta_0} $ over the interval $[T_0, T_1\land T_5]$, so
we have $\dist(X_{T_1\land T_5}, Y_{T_1\land T_5}) \geq (1- c_{34}
\dist(X_{T_0}, Y_{T_0})^{2\beta_0})\dist (X_0,Y_0)$.

An argument similar to those given above yields
$$\bP(N\leq -k_1, T_1 \geq T_2, N' \geq k_1) \leq c_{35}
\dist(X_0,Y_0)^{1+ 3\beta_0/2}.$$
 If $\{N\leq -k_1, T_1 \geq T_2, N' \geq k_1\}$ holds then
$\dist(X_{T_4}, Y_{T_4}) \geq (1- c_{36} \dist(X_{T_0},
Y_{T_0})^{2\beta_0})\dist (X_0,Y_0)$.

Finally, by Lemmas 3.2 and 3.7, we have $\bP(N \geq k_1) \leq c_{37}
\dist(X_0,Y_0)^{1+ \beta_0/2}$.

\bigskip
\noindent {\it Step 2}. Recall that $\theta$ denotes the usual
Markov shift operator and let
 $$\eqalign{
 A &= \big(\{|N|\leq k_1\}\cap\{ T_1 \leq T_2\} \big)
 \cup \big( \{N< -k_1\}\cap \{N'\leq k_1\}
 \cap \{T_1 \leq T_5\}\big),\cr
 T_ 7 &= T_1 \bone_{A} + T_0 \circ \theta_{T_2}
 \bone_{\{|N|\leq k_1\}\cap\{ T_1 > T_2\}}
 + T_0 \bone_{\{N > k_1\}}\cr
 &\qquad + T_0 \circ \theta_{T_5}
 \bone_{\{N< -k_1\}\cap \{N'\leq k_1\} \cap \{T_1 > T_5\}}
 + T_4 \bone_{\{N< -k_1\}\cap \{N'> k_1\}}
 .
 }$$
Note that $\dist(X_{T_7},Y_{T_7}) \leq \dist(X_0,Y_0)$.

Let $D(a) = \{x\in D: \dist(x, \prt D) \leq a\}$. We will define a
number of stopping times and events involving a parameter
$\beta_3>1$ whose value will be chosen later. Recall that $Z_t$ is
the closest point on $\prt D$ to the pair $\{X_t, Y_t\}$. Let
$\delta = \dist(X_{T_7}, Y_{T_7})$ and for $k\geq 0$,
 $$\eqalignno{
 V_k &= \inf\{t\geq T_7 : X_t, Y_t \in D(\delta^{\beta_3^k})\}, \cr
 G_k &= \{|\pi/2 -\angle(\n(Z_{V_k}), X_{V_k}-Y_{V_k})| > c_1
 \dist(X_{V_k},Y_{V_k})^{\beta_0} \}.
 }$$
We will define some stopping times $V_k^j$ and related events
$A^j_k$ assuming that $V_k >T_7$ (otherwise $V_k^j$'s and
$A^j_k$'s can be defined in an arbitrary way). If $G_k^c$ holds,
we let $V_k^j = V_k$ for all $j$. We will state the definitions in
the case when $G_k $ holds and  $X_{V_k} \in \prt
D(\delta^{\beta_3^k}) \setminus \prt D$. In the case when $G_k$
holds and $Y_{V_k} \in \prt D(\delta^{\beta_3^k}) \setminus \prt
D$, the roles of $X_t$ and $Y_t$ should be interchanged in the
definitions of $V_k^j$'s and $A^j_k$'s. Let $CS_k$ be the
orthonormal coordinate system with the origin at the point in
$\prt D$ that is closest to $X_{V_k}$, whose first axis is tangent
to $\prt D$. We will write $X_t=(X^1_t, X^2_t)$ in this coordinate
system. Note that $X^1_{V_k} =0$ in $CS_k$. Let
 $$\eqalign{
 V_k^1 & = \inf\{t\geq V_k: X_t \in D(\delta^{\beta_3^k}/2)\}, \cr
 A_k^1 & = \left\{V^1_k \leq
  \inf\{t\geq V_k: X_t \in D(2\delta^{\beta_3^k})^c \hbox{  or  }
  |X^1_t | = \delta^{\beta_3^k} \}\right\},\cr
  V_k^2 & = \inf\{t\geq V^1_k: X_t \in D(\delta^{\beta_3^k})^c\}, \cr
 A_k^2 & = \left\{V^2_k \leq
  \inf\{t\geq V^1_k: X_t \in D(\delta^{\beta_3^k}/4) \hbox{  or  }
  |X^1_t | = 2\delta^{\beta_3^k} \} \right\}.\cr
 }
 $$
If $\delta$ is small then $\prt D\cap \B(X_{V_k}, 2
\dist(X_{V_k},Y_{V_k}) \lor 2 \delta^{\beta_3^k})$ is almost flat.
If events $A_k^1$ and $A_k^2$ occur, the process $X_t$ moves
towards the boundary of $D$ and then away from the boundary,
without moving too much in the horizontal direction in $CS_k$. The
result is that the distance from $Y_{V_k^2}$ to $\prt D$ is
greater than $\delta^{\beta_3^k}/8$.

It follows from Lemma 3.7 and its proof that there exists
$\delta_0, c_{38}>0$ depending only on $D$, such that if
$\delta\leq \delta_0$ then either (i) $\angle(\n(x), \n(y)) \geq
c_{38} \delta^{2\beta_3^k}$ for all $x=(x^1,x^2)\in \prt D$ and
$y=(y^1,y^2)\in \prt D$ with $3 \delta^{\beta_3^k} \leq -x^1 \leq
9\delta^{\beta_3^k}$ and $3 \delta^{\beta_3^k} \leq y^1 \leq
9\delta^{\beta_3^k}$, or (ii) $\angle(\n(x), \n(y)) \geq c_{38}
\delta^{2\beta_3^k}$ for all $x,y\in \prt D$ with $10
\delta^{\beta_3^k}\leq -x^1 \leq 16\delta^{\beta_3^k}$ and $10
\delta^{\beta_3^k}\leq y^1 \leq 16\delta^{\beta_3^k}$, in $CS_k$.
We have to consider cases (i) and (ii) because there might be a
(single) $z\in \prt D$ with $-16 \delta^{\beta_3^k}\leq z^1 \leq
16\delta^{\beta_3^k}$ and $\nu(z)=0$. Depending on the sign of
$X^1 _{V_k^2} - Y^1 _{V_k^2}$, one of the following events holds,
 $$\eqalignno{
 \angle(X _{V_k^2} - Y _{V_k^2}, \n(x)) &\geq
 c_{38} \delta^{2\beta_3^k}, \qquad \hbox{  for  }
 x\in \prt D, \ 3
 \delta^{\beta_3^k} \leq x^1 \leq 9\delta^{\beta_3^k}, &(3.20)\cr
 \angle(X _{V_k^2} - Y _{V_k^2}, \n(x)) &\geq
 c_{38} \delta^{2\beta_3^k}, \qquad \hbox{  for  }
 x\in \prt D, \ 3
 \delta^{\beta_3^k} \leq -x^1 \leq 9\delta^{\beta_3^k}, &(3.21)\cr
 \angle(X _{V_k^2} - Y _{V_k^2}, \n(x)) &\geq
 c_{38} \delta^{2\beta_3^k}, \qquad \hbox{  for  }
 x\in \prt D, \ 10 \delta^{\beta_3^k}\leq x^1 \leq
 16\delta^{\beta_3^k}. &(3.22)\cr
 \angle(X _{V_k^2} - Y _{V_k^2}, \n(x)) &\geq
 c_{38} \delta^{2\beta_3^k}, \qquad \hbox{  for  }
 x\in \prt D, \ 10 \delta^{\beta_3^k}\leq -x^1 \leq
 16\delta^{\beta_3^k}. &(3.23)}$$
We will discuss only cases (3.20) and (3.22). The other cases are
symmetric---we leave them to the reader. In case (3.20) we let
 $$\eqalign{
  V_k^{3} & = \inf\{t\geq V^2_k:
  X^1_t  = 6\delta^{\beta_3^k}\}, \cr
 A_k^{3} & = \left\{V^{3}_k \leq
 \inf\{t\geq V^2_k: X^1_t  = -3\delta^{\beta_3^k} \hbox{  or  }
  X_t \in  D(2\delta^{\beta_3^k})^c \cup D(\delta^{\beta_3^k}/2) \}
  \right\},\cr
  V_k^{4} & = \inf\{t\geq V^3_k:
  X_t \in  D(\delta^{\beta_3^{k+1}}) \}, \cr
 A_k^{4} & = \left\{ V^4_{k} \leq
 \inf\{t\geq V^{3}_k:  X_t \in  D(3\delta^{\beta_3^k})^c \hbox{  or  }
 |X^1_t - 6 \delta^{\beta_3^k}| = \delta^{\beta_3^k} \}\right\}.\cr
 }
 $$
In case (3.22), we let
 $$\eqalign{
  V_k^{3} & = \inf\{t\geq V^2_k:
  X^1_t  = 13\delta^{\beta_3^k}\}, \cr
 A_k^{3} & = \left\{V^{3}_k \leq
 \inf\{t\geq V^2_k: X^1_t  = -3\delta^{\beta_3^k} \hbox{  or  }
  X_t \in  D(2\delta^{\beta_3^k})^c \cup D(\delta^{\beta_3^k}/2) \}
  \right\},\cr
    V_k^{4} & = \inf\{t\geq V^3_k:
  X_t \in  D(\delta^{\beta_3^{k+1}}) \}, \cr
 A_k^{4} & = \left\{ V^4_{k} \leq
 \inf\{t\geq V^{3}_k:  X_t \in  D(3\delta^{\beta_3^k})^c \hbox{  or  }
 |X^1_t - 13 \delta^{\beta_3^k}| = \delta^{\beta_3^k} \}\right\}.\cr
 }$$
In either case, let $ A_k = A_k^1 \cap A_k^2 \cap A_k^3 \cap
A_k^{4}$, and similarly in cases (3.21) and (3.23).

We will assume that $\delta_0>0$ is so small that
$\delta^{\beta_3^{k+1}} < \delta^{\beta_3^k}/4$ for $\delta\leq
\delta_0$. We will later impose an upper bound on $\beta_3$ which,
in turn, will impose an upper bound on $\delta_0$. Note that given
this assumption about $\delta_0$, if $\dist(X_0, Y_0) \leq
\delta_0$ and $A_k$ holds then $V^4_k \leq V_{k+1}$, so we can
estimate the probability of the intersection of consecutive
$A_k$'s using the strong Markov property at times $V_k$. Let
 $$\eqalign{
 V^5_k &= \inf\{t > V_{k}: X_t \in \prt D \},\cr
 C_k &= \left\{\sup_{t\in[V_k, V_k^5]}
 \dist(X_{t} , X_{V_k}) \leq \delta^{\beta_3^{k-1}}/2\right\}.
 }
 $$
By Lemma 3.2,
 $$\bP(C_k^c\mid X_{V_k} \in \prt D(\delta^{\beta_3^k})
\setminus \prt D) \leq c_{39}
 \delta^{\beta_3^k}/\delta^{\beta_3^{k-1}} = c_{39}
 \delta^{\beta_3^{k-1}(\beta_3 -1)}.\eqno(3.24)$$

We will find a lower bound for $\sup_{t\in[V_k,
V_{k}^5]}\dist(X_t,Y_t)$ assuming that $\dist(X_{T_7},Y_{T_7})
\geq \delta^{\beta_3^k}$ and $G_k \cap A_k \cap C_{k+1}$ occurred.
First consider the case when $|X^1_{V_k} - Y^1_{V_k}| >
\delta^{\beta_3^k}/8$. Then it is easy to see that the distance
between $X_t$ and $Y_t$ is reduced between times $V_k$ and
$V_{k}^5$ by at most a constant factor $c_{40}$, so
$\dist(X_t,Y_t) \geq c_{41} \delta^{\beta_3^k}$ for $t\in[V_k,
V_{k}^5]$. Next suppose that $|X^1_{V_k} - Y^1_{V_k}| \leq
\delta^{\beta_3^k}/8$. Then $|X^1_{t} - Y^1_{t}| \leq
\delta^{\beta_3^k}/4$ for $t\in[V_k, V_{k}^5]$, because the
boundary of $D$ is ``flat'' in the neighborhood under
consideration. After time $V_k^2$, processes $X_t$ and $Y_t$ move
along $\prt D$ without touching it, to the place where the angle
between the line passing through both particles and the normal to
the boundary of the domain is bounded below by $ c_{38}
\delta^{2\beta_3^k}$. It follows that for $t\in [V_{k}, V_{k}^5]$,
the process $Y_t$ is reflecting on the part of the boundary where
the angle between the line passing through both particles and the
normal to the boundary of the domain is bounded below by $ c_{38}
\delta^{2\beta_3^k}$. Hence, the distance between $X_t$ and $Y_t$
is reduced between times $V_k$ and $V_{k}^5$ by at most a factor
of $c_{42} \delta^{2\beta_3^k}$. This implies that if $G_k \cap
A_k \cap C_{k+1}$ holds then $\dist(X_t,Y_t) \geq c_{43}
\delta^{3\beta_3^k}$ for $t\in[V_k, V_{k}^5]$.

Let
 $$F_k = \bigcup_{k\leq m < 2k} \left (G_m^c
 \cup (G_m \cap A_m \cap C_{m+1} )\right),$$
and note that
 $$\bigcup_{k\leq m < 2k} A_m \cap \bigcap_{k\leq m < 2k} C_{m+1}
 \subset \bigcup_{k\leq m < 2k} \left (G_m^c
 \cup (G_m \cap A_m )\right)
 \cap \bigcap_{k\leq m < 2k} C_{m+1} \subset F_k.$$
It is elementary to see that the probability of $A_k$ is bounded
below by $p_1>0$, not depending on $k$ or $\beta_3$, assuming
$\delta$ is small. By the strong Markov property applied at
stopping times $V_m$, we have $\bP(\bigcap _{k\leq m < 2k} A_m^c)
\leq (1-p_1)^k$, so
 $$\bP\left(\bigcup_{k\leq m< 2k} A_m\right) \geq 1-
 (1-p_1)^{k}.\eqno(3.25)
 $$
By (3.24),
 $$\bP\left(\bigcap_{k\leq m< 2k}C_m\right) \geq 1 - \sum_{k\leq m< 2k}
 c_{39} \delta^{\beta_3^{m-1}(\beta_3 -1)}.
 $$
The quantity on the right hand side is bounded below by $1-
(1-p_1)^{k}$, for small $\delta$. This and (3.25) imply that
$\bP(F_k) \geq 1- 2(1-p_1)^{k}$, for small $\delta$. Let
 $$\eqalign{T_8 = \inf\{t>T_7: \dist&(X_t ,\prt D)
 \land \dist(Y_t ,\prt D)
 \leq \dist(X_t,Y_t), \cr
 & |\pi/2 -\angle(\n(Z_{t}), X_{t}-Y_{t})| \leq c_1
 \dist(X_t,Y_t)^{\beta_0}
  \}. }$$
Note that if $G^c_k\cup (G_k \cap A_k \cap C_{k+1})$ occurs then
$T_8 \leq V_{k}^5$ and $\dist(X_{T_8}, Y_{T_8}) \geq c_{43}
\delta^{3\beta_3^k}$. In view of the estimate for the probability
of $F_k$, we see that for $k\geq 1$,
 $$\bP(\dist(X_{T_8}, Y_{T_8}) \leq c_{43} \delta^{3\beta_3^{2k}} )\leq
  2(1-p_1)^{k}.$$
We choose $\beta_3 >1$ so that $\beta^2_3 (1-p_1) < 1$.

\bigskip
\noindent{\it Step 3}. Let $c_{44}$ be the same as $c_2$ in the
statement of Lemma 3.3, and let $c_{45}$ be the same as $c_1$ in
the statement of that lemma. Let $T_9^1 = T_8$ and for $k \geq 1$,
$$\eqalign{
 T^k_{10} &= \inf\{t\geq T_9^k:
 \dist(X_t, X_{T^k_9}) \lor \dist(Y_t, Y_{T^k_9})
 \geq 2(\dist(X_{T^k_9},\prt D)\lor \dist(Y_{T^k_9},\prt D))\},\cr
 T^k_{11} &= \inf\{t \geq T_9^k: X_t \in \prt D \}, \cr
 T^k_{12} &= \inf\{t \geq T_9^k: Y_t \in \prt D \}, \cr
 T^k_{13} &= \inf\{t \geq T_{11}^k: L^X_t -L^X_{T^k_{11}}
 \geq c_{44} \dist(Y_{T^k_{11}}, \prt D) \}
 \bone_{\{T^k_{11} \leq T^k_{12}\}}, \cr
 T^k_{14} &= \inf\{t \geq T_{12}^k: L^Y_t -L^Y_{T^k_{11}}
 \geq c_{44} \dist(X_{T^k_{11}}, \prt D) \}
 \bone_{\{T^k_{11}> T^k_{12}\}}, \cr
 T^k_{15} &= \inf\{t \geq T_{11}^k: Y_t \in \prt D \}
 \bone_{\{T^k_{11} \leq T^k_{12}\}}, \cr
 T^k_{16} &= \inf\{t \geq T_{11}^k: X_t \in \prt D \}
 \bone_{\{T^k_{11}> T^k_{12}\}}, \cr
 T^{k+1}_9 &= \inf\{t\geq T_9^k:
 \dist(X_t, X_{T^k_9}) \lor \dist(Y_t, Y_{T^k_9})
 \geq 2 c_{45} (\dist(X_{T^k_9},\prt D)\lor \dist(Y_{T^k_9},\prt D))\}\cr
 &\quad \land (T^k_{13} + T^k_{14}) \land ( T^k_{15} +T^k_{16}).
 }$$
Note that $\dist(X_{T_9^1},Y_{T_9^1})=\dist(X_{T_8},Y_{T_8}) \leq
\dist(X_0,Y_0)$. It is easy to see that if $\eps_0>0$ is small and
$\dist(X_0,Y_0)\leq\eps_0$, we have $\bP(T^k_{11} < T^k_{10}) \geq
p_2>0$, where $p_2$ depends only on $D$. By the strong Markov
property applied at $T^k_{11}\land T^k_{12}$ and Lemma 3.3 (i),
for some $p_3>0$,
 $$\bP(S_1 \leq T^{k+1}_{9}
 \mid \F_{T_9^k}, \{S_1 > T_9^k\}) \geq p_3.$$
Let $k_2$ be such that $(1-p_3)^{k_2-1} \leq 1/2$. Then
 $$\bP(S_1 \geq T_{9}^{k_2}) \leq (1-p_3) ^{k_2-1}\leq 1/2.$$

\bigskip
\noindent{\it Step 4}. Let $T_{17} = T_9^{k_2}$ and note that
$$
 \dist(X_{T^1_9}, X_{T_{17}}) \lor \dist(Y_{T^1_9}, Y_{T_{17}})
 \leq c_{46}(\dist(X_{T^1_9},\prt D)\lor \dist(Y_{T^1_9},\prt D)),
$$
where $c_{46} = (2c_{45})^{k_2}$. For $t$ between $T^1_9$ and
$T_{17}$, the angle between the vector of reflection for any of
the processes and $X_t-Y_t$ is bounded below by a quantity
depending on $N$; we will next discuss this dependence and its
consequences. We will examine various cases in the same order as
in Step 1 and we will also recall some estimates from Step 1.
There will be many cases to consider---we will label them for
future reference.

We start with a general remark that applies to many of the cases
discussed below. If $T_1 \leq T_2$ then the lower bounds for
$\dist(X_{T_1\land T_2}, Y_{T_1\land T_2})$ obtained in Step 1
apply also to $\dist(X_{T_{17}}, Y_{T_{17}}) $, for the same
reasons, but with constants that may be different. The same is
true when $T_1\leq T_5$.

(a) Consider the case when $N=k$ and $-k_1 \leq k \leq 0$. Then we
have $\bP(N=k) \leq c_{47} \dist (X_0,Y_0) 2^{-k}$,
 $$\bP(N=k, T_1 \geq T_2) \leq c_{48}
 \dist(X_{T_0}, Y_{T_0})^{2-\beta_0},$$
and $(1- c_{49} 2^{2k})\dist (X_0,Y_0) \leq \dist(X_{T_1\land
T_2}, Y_{T_1\land T_2}) \leq \dist (X_0,Y_0)$. If $T_1 \leq T_2$
then $\dist(X_{T_{17}}, Y_{T_{17}})\geq (1- c_{50} 2^{2k})\dist
(X_0,Y_0)$. Note that the distance between $X_t$ and $Y_t$ does
not increase before time $T_8=T_9^1$. The increase of the local
time $L^X_t$ between times $T_9^1$ and $T_{17}$ is bounded by
$c_{51} \dist (X_0,Y_0)$, and a similar bound holds for the
increment of $L^Y_t$, so, according to Lemma 3.8,
$\dist(X_{T_{17}},Y_{T_{17}}) \leq \dist (X_0,Y_0)(1+c_{52}
\dist(X_0,Y_0))$, assuming that $\dist(X_0,Y_0)$ is small.
Combining the two estimates, we obtain
 $$(1- c_{50} 2^{2k})\dist(X_0,Y_0) \leq
 \dist(X_{T_{17}},Y_{T_{17}}) \leq
 \dist (X_0,Y_0)(1+c_{52}\dist(X_0,Y_0))
 .$$

(b) Next consider the case when $N=k$, $-k_1 \leq k \leq 0$, and
$T_2 < T_1$. Then, $\dist(X_7,Y_7) \leq \dist(X_0,Y_0)$,
$\dist(X_{17},Y_{17}) \leq c_{53}\dist(X_8,Y_8)$, and using Step
2, for $n\geq 1$,
 $$\bP(N=k, T_1 \geq T_2, \dist(X_{T_{17}}, Y_{T_{17}}) \leq
 c_{54}\dist(X_0,Y_0)^{3\beta_3^{2n}} )
 \leq c_{55} \dist(X_{T_0}, Y_{T_0})^{2-\beta_0}(1-p_1)^{n+1}.$$

(c) Assume that $0 < k \leq k_1$. Then $\bP(N=k) \leq c_{56} \dist
(X_0,Y_0) 2^{-k/2}$,
 $$\bP(N=k, T_1 \geq T_2) \leq c_{57} \dist(X_{T_0}, Y_{T_0})^{2-\beta_0}
 2^{-k/2},$$
and $c_{58} 2^{-k}\dist (X_0,Y_0) \leq \dist(X_{T_1\land T_2},
Y_{T_1\land T_2}) \leq \dist (X_0,Y_0)$. If $T_1 \leq T_2$ then we
have $\dist(X_{T_{17}}, Y_{T_{17}})\geq c_{59} 2^{-k}\dist
(X_0,Y_0)$. The bound $\dist(X_{T_{17}},Y_{T_{17}}) \leq \dist
(X_0,Y_0)(1+c_{60} \dist(X_0,Y_0))$ holds for the same reason as
in case (a). Combining the two estimates, we see that
 $$c_{58} 2^{-k}\dist(X_0,Y_0) \leq
 \dist(X_{T_{17}},Y_{T_{17}}) \leq
 \dist (X_0,Y_0)(1+c_{60}\dist(X_0,Y_0))
 .$$

(d) Next consider the case when $N=k$, $0 < k \leq k_1$, and $T_2
< T_1$. Then, $\dist(X_7,Y_7) \leq \dist(X_0,Y_0)$,
$\dist(X_{17},Y_{17}) \leq c_{61}\dist(X_8,Y_8)$, and using Step
2, for $n\geq 1$,
 $$\eqalign{
 \bP&(N=k, T_1 \geq T_2, \dist(X_{T_{17}}, Y_{T_{17}}) \leq
 c_{62}(c_{58} 2^{-k}\dist(X_0,Y_0))^{3\beta_3^{2n}} )\cr
 &\leq c_{63} \dist(X_{T_0}, Y_{T_0})^{2-\beta_0}2^{-k/2}(1-p_1)^{n+1}.
 }$$

(e) The next case is when $N \leq -k_1$. We will use the trivial
estimate $\bP(N\leq -k_1) \leq 1$. We have $\dist(X_{T_1\land T_2},
Y_{T_1\land T_2}) \geq (1- c_{64} \dist(X_{T_0},
Y_{T_0})^{2\beta_0})\dist (X_0,Y_0)$. If $T_1 \leq T_2$ then
 $$(1- c_{65}
 \dist(X_{T_0}, Y_{T_0})^{2\beta_0})\dist (X_0,Y_0) \leq
 \dist(X_{T_{17}},Y_{T_{17}}) \leq
 \dist (X_0,Y_0)(1+c_{66}\dist(X_0,Y_0))
 .$$

(f) Recall that $N'$ is defined by the following conditions,
$X_{T_4} \in M_{N'}$ if $X_{T_4} \in \prt D$, and $Y_{T_4} \in
M_{N'}$ if $Y_{T_4} \in \prt D$. Suppose that $-k_1 \leq k \leq
0$. Then
 $$\bP(N\leq -k_1, N'=k, T_1 \geq T_5)
 \leq c_{67} \dist(X_{T_0}, Y_{T_0})^{2},$$
and $(1- c_{68} 2^{2k})\dist (X_0,Y_0) \leq \dist(X_{T_1\land
T_5}, Y_{T_1\land T_5}) \leq \dist (X_0,Y_0)$. If $N\leq -k_1$,
$N'=k$, $-k_1 \leq k \leq 0$, and $T_1 \leq T_5$ then we have
$\dist(X_{T_{17}}, Y_{T_{17}})\geq(1- c_{69} 2^{2k})\dist
(X_0,Y_0)$ and
 $$(1- c_{69} 2^{2k})\dist (X_0,Y_0) \leq
 \dist(X_{T_{17}},Y_{T_{17}}) \leq
 \dist (X_0,Y_0)(1+c_{70}\dist(X_0,Y_0))
 .$$

(g) If $N\leq -k_1$, $N'=k$, $-k_1 \leq k \leq 0$, and $T_5 < T_1$
then, $\dist(X_7,Y_7) \leq \dist(X_0,Y_0)$, $\dist(X_{17},Y_{17})
\leq c_{71}\dist(X_8,Y_8)$, and using Step 2, for $n\geq 1$,
 $$\eqalign{
 \bP&(N\leq -k_1, N'=k, T_1 \geq T_5, \dist(X_{T_{17}}, Y_{T_{17}}) \leq
 c_{72} \dist(X_0,Y_0)^{3\beta_3^{2n}} )\cr
 &\leq c_{73} \dist(X_{T_0}, Y_{T_0})^{2}(1-p_1)^{n+1}.
 }$$

(h) If $0 < k \leq k_1$ then
 $$\bP(N\leq -k_1, N'=k, T_1 \geq T_5)
 \leq c_{74} \dist(X_{T_0}, Y_{T_0})^{2}
 2^{-k/2},$$
and $c_{75} 2^{-k}\dist (X_0,Y_0) \leq \dist(X_{T_1\land T_5},
Y_{T_1\land T_5}) \leq \dist (X_0,Y_0)$. If $N\leq -k_1$, $N'=k$,
$0 < k \leq k_1$, and $T_1 \leq T_5$ then $\dist(X_{T_{17}},
Y_{T_{17}})\geq c_{76} 2^{-k}\dist (X_0,Y_0)$ and
 $$c_{76} 2^{-k}\dist (X_0,Y_0) \leq
 \dist(X_{T_{17}},Y_{T_{17}}) \leq
 \dist (X_0,Y_0)(1+c_{77}\dist(X_0,Y_0))
 .$$

(i) If $N\leq -k_1$, $N'=k$, $0 < k \leq k_1$, and $T_5 < T_1$
then, $\dist(X_7,Y_7) \leq \dist(X_0,Y_0)$, $\dist(X_{17},Y_{17})
\leq c_{78}\dist(X_8,Y_8)$, and using Step 2, for $n\geq 1$,
 $$\eqalign{
 \bP&(N\leq -k_1, N'=k, T_1 \geq T_5, \dist(X_{T_{17}}, Y_{T_{17}}) \leq
 c_{79} (c_{75} 2^{-k} \dist(X_0,Y_0))^{3\beta_3^{2n}} )\cr
 &\leq c_{80} \dist(X_{T_0}, Y_{T_0})^{2}2^{-k/2}(1-p_1)^{n+1}.
 }$$

(j) The next case to be considered is when $N\leq -k_1$ and
$N'\leq -k_1 $. We have $\bP(N'\leq -k_1\mid N\leq -k_1, T_1 \geq
T_2) \leq 1$,
 $$\bP(N\leq -k_1, N'\leq -k_1, T_1 \geq T_5)
 \leq c_{81} \dist(X_{T_0}, Y_{T_0})^{2},$$
and $\dist(X_{T_1\land T_5}, Y_{T_1\land T_5}) \geq (1- c_{82}
\dist(X_{T_0}, Y_{T_0})^{2\beta_0})\dist (X_0,Y_0)$. If $N\leq
-k_1$, $N'\leq -k_1 $, and $T_1 \leq T_5$ then $\dist(X_{T_{17}},
Y_{T_{17}})\geq (1- c_{83} \dist(X_{T_0},
Y_{T_0})^{2\beta_0})\dist (X_0,Y_0)$ and
 $$(1- c_{83}\dist(X_{T_0}, Y_{T_0})^{2\beta_0})\dist (X_0,Y_0) \leq
 \dist(X_{T_{17}},Y_{T_{17}}) \leq
 \dist (X_0,Y_0)(1+c_{84}\dist(X_0,Y_0))
 .$$

(k) If $N\leq -k_1$, $N'\leq -k_1 $, and $T_5 < T_1$ then,
$\dist(X_7,Y_7) \leq \dist(X_0,Y_0)$, $\dist(X_{17},Y_{17}) \leq
c_{85}\dist(X_8,Y_8)$, and using Step 2, for $n\geq 1$,
 $$\eqalign{
 \bP&(N\leq -k_1, N'\leq -k_1, T_1 \geq T_5, \dist(X_{T_{17}}, Y_{T_{17}}) \leq
 c_{86} \dist(X_0,Y_0)^{3\beta_3^{2n}} )\cr
 &\leq c_{87} \dist(X_{T_0}, Y_{T_0})^{2}(1-p_1)^{n+1}.
 }$$

(l) Consider the case when $N'\geq k_1$ and note that $\bP(N\leq
-k_1, T_1 \geq T_2, N' \geq k_1) \leq c_{88} \dist(X_0,Y_0)^{1+
3\beta_0/2}$. If $N\leq -k_1$, $N' \geq k_1$, and $T_1 \geq T_2$
then we have $\dist(X_{T_4}, Y_{T_4}) \geq (1- c_{89}
\dist(X_{T_0}, Y_{T_0})^{2\beta_0})\dist (X_0,Y_0)$. If, in
addition, $T_1 \leq T_5$ then
 $$\dist(X_{T_{17}}, Y_{T_{17}})\geq (1- c_{90} \dist(X_{T_0},
 Y_{T_0})^{2\beta_0})\dist (X_0,Y_0)$$
and
 $$(1- c_{90}\dist(X_{T_0}, Y_{T_0})^{2\beta_0})\dist (X_0,Y_0) \leq
 \dist(X_{T_{17}},Y_{T_{17}}) \leq
 \dist (X_0,Y_0)(1+c_{91}\dist(X_0,Y_0))
 .$$

(m) If $N\leq -k_1$, $N' \geq k_1$, and $T_5 < T_1$ then,
$\dist(X_7,Y_7) \leq \dist(X_0,Y_0)$, $\dist(X_{17},Y_{17}) \leq
c_{92}\dist(X_8,Y_8)$, and using Step 2, for $n\geq 1$,
 $$\eqalign{
 \bP&(N\leq -k_1, N'\leq -k_1, T_1 \geq T_5, \dist(X_{T_{17}}, Y_{T_{17}}) \leq
 c_{93} \dist(X_0,Y_0)^{3\beta_3^{2n}} )\cr
 &\leq c_{94} \dist(X_{T_0}, Y_{T_0})^{2}(1-p_1)^{n+1}.
 }$$

(n) Finally, we consider the case $N\geq k_1$. We have $\bP(N \geq
k_1) \leq c_{95} \dist(X_0,Y_0)^{1+ \beta_0/2}$. Then,
$\dist(X_7,Y_7) \leq \dist(X_0,Y_0)$, $\dist(X_{17},Y_{17}) \leq
c_{96}\dist(X_8,Y_8)$, and using Step 2, for $n\geq 1$,
 $$\eqalign{
 \bP&(N\leq -k_1, N'\leq -k_1, T_1 \geq T_5,
 \dist(X_{T_{17}}, Y_{T_{17}}) \leq
 c_{97} \dist(X_0,Y_0)^{3\beta_3^{2n}} )\cr
 &\leq c_{98} \dist(X_{T_0}, Y_{T_0})^{1+\beta_0/2}(1-p_1)^{n+1}.
 }$$

The estimates for values of $\dist(X_{T_{17}}, Y_{T_{17}})$ and
the corresponding probabilities listed above as (a)-(n) yield the
following inequality. Its lines are labelled according to the case
they represent.
 $$\eqalignno{
 \bE& | \log \dist(X_{T_{17}},Y_{T_{17}})
 -\log \dist(X_0,Y_0)|&(3.26)\cr
 &\leq \sum_{-k_1 \leq k \leq 0} c_{99} \dist(X_0,Y_0) 2^{-k}
 (c_{100} 2^{2k} + c_{101} \dist(X_0,Y_0))
 &(a) \cr
 &+ \sum_{-k_1 \leq k \leq 0}
 \sum_{n\geq 1} c_{102} \dist(X_0,Y_0)^{2-\beta_0}
 (1-p_1)^{n+1} (c_{103} + (3 \beta_3^{2n} - 1) |\log \dist(X_0,Y_0)|)
 &(b) \cr
 &+ \sum_{0 \leq k \leq k_1} c_{104} \dist(X_0,Y_0) 2^{-k/2}
 (c_{105} k + c_{106} \dist(X_0,Y_0))
 &(c)\cr
 &+ \sum_{0 \leq k \leq k_1}\sum_{n\geq 1}
 c_{107} \dist(X_0,Y_0)^{2-\beta_0} 2^{-k/2}
 (1-p_1)^{n+1} \cr
 &\qquad \times (c_{108} + c_{109} \beta_3^{2n}
 + c_{110} k \beta_3^{2n} + (3\beta_3^{2n}
  - 1) |\log \dist(X_0,Y_0)|)
 &(d) \cr
 &+ c_{111} \dist(X_0,Y_0)^{2\beta_0} + c_{112} \dist(X_0,Y_0)
 &(e) \cr
 &+ \sum_{-k_1 \leq k \leq 0} c_{113} \dist(X_0,Y_0)^2
 (c_{114} 2^{2k} + c_{115} \dist(X_0,Y_0))
 &(f) \cr
 &+ \sum_{-k_1  \leq k \leq 0}\sum_{n\geq 1}
 c_{116} \dist(X_0,Y_0)^{2} (1-p_1)^{n+1}
 (c_{117} + (3 \beta_3^{2n} - 1) |\log\dist(X_0,Y_0)|)
 &(g)\cr
 &+ \sum_{0 \leq k \leq k_1} c_{118} \dist(X_0,Y_0)^2 2^{-k/2}
 (c_{119} k + c_{120} \dist(X_0,Y_0))
 &(h)\cr
 &+ \sum_{0 \leq k \leq k_1}\sum_{n\geq 1}
 c_{121} \dist(X_0,Y_0)^{2} 2^{-k/2}
 (1-p_1)^{n+1} \cr
 &\qquad \times (c_{122} + c_{123} \beta_3^{2n}
 + c_{124} k \beta_3^{2n} + (3\beta_3^{2n}
  - 1) |\log \dist(X_0,Y_0)|)
 &(i) \cr
 &+ c_{125} \dist(X_0,Y_0)^{2\beta_0} + c_{126} \dist(X_0,Y_0)
 &(j) \cr
 &+ \sum_{n\geq 1}
 c_{127} \dist(X_0,Y_0)^{2} (1-p_1)^{n+1}
 (c_{128} + (3 \beta_3^{2n} - 1) |\log\dist(X_0,Y_0)|)
 &(k)\cr
 &+ c_{129} \dist(X_0,Y_0)^{1+3\beta_0/2}
 (c_{130} \dist(X_0,Y_0)^{2\beta_0} + c_{131} \dist(X_0,Y_0))
 &(l) \cr
 &+ \sum_{n\geq 1}
 c_{132} \dist(X_0,Y_0)^{2} (1-p_1)^{n+1}
 (c_{133} + (3 \beta_3^{2n} - 1) |\log\dist(X_0,Y_0)|)
 &(m) \cr
 &+ \sum_{n\geq 1}
 c_{134} \dist(X_0,Y_0)^{1+\beta_0/2} (1-p_1)^{n+1}
 (c_{135} + (3 \beta_3^{2n} - 1) |\log\dist(X_0,Y_0)|).
 &(n) \cr
 }$$
Recall that $\beta_0 \in (1/2,1)$ and $\beta^2_3 (1-p_1) < 1$.
Given these constraints on the values of the parameters, it is
straightforward to check that (3.26) implies that
 $$\bE | \log \dist(X_{T_{17}},Y_{T_{17}})
 -\log \dist(X_0,Y_0)| \leq c_{136} \dist(X_0,Y_0).$$

\bigskip \noindent{\it Step 5}. It is easy to check that our
estimates on the size of $\dist(X_t,Y_t)$ apply not only at
$T_{17}$ but on the whole interval $[0,T_{17}]$. Hence,
 $$\bE  \sup_{t \in [0,T_{17}]}
 |\log \dist(X_{t},Y_{t})
 -\log \dist(X_0,Y_0)| \leq c_{136} \dist(X_0,Y_0).$$

At several places in our argument we have assumed that
$\dist(X_0,Y_0)$ is small. Let $\eps_1>0$ be such that the last
inequality holds if $\dist(X_0,Y_0) \leq \eps_1$. Let $Q_0=0$,
$Q_1 = T_{17} \land \tau^+( \eps_1)$ and $Q_k = Q_1 \circ
\theta_{Q_{k-1}}$ for $k\geq 2$. Note that if
$\dist(X_{Q_{k}},Y_{Q_{k}}) = \eps_1$ then
$\dist(X_{Q_{n}},Y_{Q_{n}}) = \eps_1$ for all $n\geq k$. If $Q_k =
T_{17} \circ \theta_{Q_{k-1}}< \tau^+( \eps_1) \circ
\theta_{Q_{k-1}}$ then we can apply the argument given in Steps
1-4 to the post-$Q_k$ process, by the strong Markov property,
because condition (3.18) is satisfied for $t=Q_k$ in place of
$t=0$. It follows that if $\dist(X_0,Y_0) \leq \eps_1$ then
$$\bE  \sup_{t \in [0, Q_{1}]}
 |\log \dist(X_{t},Y_{t})
 -\log \dist(X_{0},Y_{0})| \leq c_{136} \dist(X_{0},Y_{0}),$$
 and
 $$\eqalign{
 \bE&(\sup_{t \in [Q_{k-1}, Q_{1}]}
 | \log \dist(X_{t},Y_{t})
 -\log \dist(X_{Q_{k-1}},Y_{Q_{k-1}})| \mid \F_{Q_{k-1}}) \cr
 &\leq c_{137} \dist(X_{Q_{k-1}},Y_{Q_{k-1}}).}$$
The argument given in part (a) of Step 4 shows that, a.s.,
 $$ \dist(X_{Q_{k}},Y_{Q_{k}}) \leq \dist(X_{Q_{k-1}},Y_{Q_{k-1}})
  + c_{138} \dist(X_{Q_{k-1}},Y_{Q_{k-1}})^2.
 $$
If $\dist(X_{Q_{k-1}},Y_{Q_{k-1}}) \leq 2 \dist(X_{0},Y_{0})$ then
 $$c_{138} \dist(X_{Q_{k-1}},Y_{Q_{k-1}})^2
 \leq 4 c_{138} \dist(X_{Q_{0}},Y_{Q_{0}})^2$$
and, therefore, $\dist(X_{Q_{k}},Y_{Q_{k}}) \leq 2
\dist(X_{0},Y_{0})$ for
 $$k \leq \dist(X_{Q_{0}},Y_{Q_{0}}) / (4
 c_{138} \dist(X_{Q_{0}},Y_{Q_{0}})^2) = 1/(4 c_{138}
 \dist(X_{0},Y_{0})).$$
This implies that for $k \leq 1/(4 c_{138} \dist(X_{0},Y_{0}))$,
 $$\bE ( \sup_{t \in [Q_{k-1}, Q_{k}]}
 |\log \dist(X_{t},Y_{t})
 -\log \dist(X_{Q_{k-1}},Y_{Q_{k-1}})|  \mid \F_{Q_{k-1}})
 \leq c_{139} \dist(X_{0},Y_{0}),$$
and
 $$\bE ( \sup_{t \in [Q_{k-1}, Q_{k}]}
 |\log \dist(X_{t},Y_{t})
 -\log \dist(X_{0},Y_{0})|  \mid \F_{Q_{k-1}})
 \leq c_{140} k\dist(X_{0},Y_{0}).$$
For $k \geq 1/(4 c_{138} \dist(X_{0},Y_{0}))$, we use the bound
$\dist(X_{Q_{k-1}},Y_{Q_{k-1}}) \leq \eps_1$ to conclude that
 $$\bE ( \sup_{t \in [Q_{k-1}, Q_{k}]}
 |\log \dist(X_{t},Y_{t})
 -\log \dist(X_{Q_{k-1}},Y_{Q_{k-1}})|
 \mid \F_{Q_{k-1}}) \leq  c_{141},$$
and
 $$\bE ( \sup_{t \in [Q_{k-1}, Q_{k}]}
 |\log \dist(X_{t},Y_{t})
 -\log \dist(X_{0},Y_{0})|  \mid \F_{Q_{k-1}}) \leq  c_{142}k.$$

By Step 3, if $\dist(X_0,Y_0) \leq \eps_1$ then $\bP(S_1 \geq Q_k)
\leq 2^{-k}$. Hence
 $$\eqalignno{
 \bE & \left(|\log \dist(X_{S_1},Y_{S_1})
 -\log \dist(X_{0},Y_{0})|\cdot \bone\left( \bigcup_{k\geq 0}
 \{ S_1 \in [Q_k, Q_{k+1}]\}\right)\right)
 &(3.27)\cr
 & \leq
 \sum_{k\geq 0} \bE (\bone _{\{ S_1 \in [Q_k, Q_{k+1}]\}}
   \sup_{t \in [Q_k, Q_{k+1}]}
 |\log \dist(X_{t},Y_{t})
 -\log \dist(X_{0},Y_{0})|) \cr
 & \leq
 \sum_{k\geq 0} \bE (\bone _{\{ S_1 \geq Q_k\}}
  \sup_{t \in [Q_k, Q_{k+1}]}
 |\log \dist(X_{t},Y_{t})
 -\log \dist(X_{0},Y_{0})|) \cr
 & \leq \sum_{k \leq  1/(4 c_{138} \dist(X_{0},Y_{0})) }
 2^{-k}  c_{140} k\dist(X_{0},Y_{0})
 + \sum_{k >  1/(4 c_{138} \dist(X_{0},Y_{0})) }
 2^{-k}  c_{142}k \cr
 & \leq c_{143} \dist(X_{0},Y_{0})
 + c_{144} \dist(X_{0},Y_{0})^{-1}
 \exp( - c_{145} \dist(X_{0},Y_{0})^{-1}) \cr
 & \leq c_{146} \dist(X_{0},Y_{0}).
 }$$

Suppose that $ \dist(X_{0},Y_{0}) \leq \eps_1 ^2$ and let
 $$\eqalign{
 W' &= \min\{Q_k: \dist(X_{Q_k},Y_{Q_k}) = \eps_1\},\cr
 W'' &= \inf\{t \geq 0: X_t \in \prt D,
 \n(X_t)\cdot(Y_t - X_t) \leq 0\}\cr
  & \quad\quad \land
 \inf\{t \geq 0: Y_t \in \prt D,
 \n(Y_t)\cdot(X_t-Y_t) \leq 0\},
 }$$
with the convention that $\inf \emptyset =\infty$. Note that
$\sup_{t\in[0, W'')}\dist(X_{t},Y_{t}) \leq \dist(X_{0},Y_{0})$.
This implies that $W''\leq W'$, and if $W''<\infty$ then,
 $$\dist(X_{W''},\prt D) \lor \dist(Y_{W''},\prt D)
 \leq c_{147} \dist(X_{W''},Y_{W''})^2
 \leq c_{147} \dist(X_{0},Y_{0})^2.$$
By Lemma 3.8,
 $$\eqalignno{
 (L^X_{W'} -L^X_{W''})+ (L^Y_{W'} -L^Y_{W''})
 &\geq c_{148} |\log \dist(X_{W''},Y_{W''}) -\log
 \dist(X_{W'},Y_{W'})|&(3.28)\cr
 &\geq c_{148} |\log \dist(X_{0},Y_{0}) -\log
 \dist(X_{W'},Y_{W'})| \geq c_{149}.
 }$$
Suppose $\eps_1$ is less than $\eps_0$ in Lemma 3.3. Then, by
Lemma 3.3 (ii) and the strong Markov property applied at $W''$,
 $$\bE (L^X_{S_1\land W'} -L^X_{ W''})
 + \bE (L^Y_{S_1\land W'} -L^Y_{ W''})\leq
 c_{150} \dist(X_{0},Y_{0})^{2}.$$
This and (3.28) imply that,
 $$\bP(W' \leq S_1) \leq c_{151} \dist(X_{0},Y_{0})^{2}. \eqno(3.29)$$

Let $W'''=\inf\{t\geq W':  \dist(X_{t},Y_{t}) =
\dist(X_{0},Y_{0})\}$ and $W_1 = T_1 \circ \theta_{W'''}$. By the
last formula in Step 2 and the strong Markov property,
 $$
 \bE  |\log \dist(X_{W_1},Y_{W_1})
 -\log \dist(X_{W'''},Y_{W'''})|  \leq c_{152}.\eqno(3.30)$$
Let $Q_k^{W_1} = Q_k \circ \theta_{W_1}$. Then, by the strong
Markov property at $W_1$, (3.27), (3.29) and (3.30),
  $$\eqalign{
 \bE  &\left(|\log \dist(X_{S_1},Y_{S_1})
 -\log \dist(X_{0},Y_{0})|\cdot \bone\left( \bigcup_{k\geq 0}
 \{ S_1 \in [Q^{W_1}_k, Q^{W_1}_{k+1}]\}\right)
 \bone_{\{W'\leq S_1\}}\right)\cr
 & \leq c_{153} \dist(X_{0},Y_{0})^2.
 }$$

Let $W_k = W_{1} \circ \theta_{W_{k-1}}$ and $Q_k^{W_j} = Q_k
\circ \theta_{W_j}$. By induction we have
  $$\bP(W_k \leq S_1) \leq c_{154} \dist(X_{0},Y_{0})^{2k},$$
and for $j\geq 1$,
  $$\eqalignno{
 \bE  &\left(|\log \dist(X_{S_1},Y_{S_1})
 -\log \dist(X_{0},Y_{0})|\cdot \bone\left( \bigcup_{k\geq 0}
 \{ S_1 \in [Q^{W_j}_k, Q^{W_j}_{k+1}]\}\right)
 \bone_{\{W_j\leq S_1\}}\right)\cr
 & \leq c_{155} \dist(X_{0},Y_{0})^{2j}.&(3.31)
 }$$

Note that $\dist(X_{0},Y_{0})$ is bounded by the diameter of the
domain so
 $$\log \dist(X_{S_1},Y_{S_1}) - \log \dist(X_{0},Y_{0})
 \leq c_{156} + c_{157} |\log \dist(X_{0},Y_{0})|.$$
This, (3.29) and (3.30) imply that
 $$\bE ( |\log \dist(X_{S_1},Y_{S_1})
 -\log \dist(X_{0},Y_{0})| \cdot \bone_{\{W'\leq S_1\leq W_1\}})
 \leq c_{158} \dist(X_{0},Y_{0})^2 |\log \dist(X_{0},Y_{0})|.
 \eqno(3.32)$$
 Let $W'_k = W' \circ \theta_{W_k}$. Then, for $k\geq 1$,
 $$\bE ( |\log \dist(X_{S_1},Y_{S_1})
 -\log \dist(X_{0},Y_{0})| \cdot \bone_{\{W'_k\leq S_1\leq W_{k+1}\}})
 \leq c_{159} \dist(X_{0},Y_{0})^{2k+2} |\log \dist(X_{0},Y_{0})|.
 \eqno(3.33)$$

It is straightforward to check that
 $$\bigcup_k [Q_k, Q_{k+1}] \cup \bigcup_j \bigcup_k
 [Q_k^{W_j}, Q_{k+1}^{W_j}] \cup [W', W_1]
 \cup \bigcup_k [W'_k, W_{k+1}] = [0,\infty).$$
Hence, part (i) of the lemma follows from (3.27), (3.31), (3.32)
and (3.33).

(ii) Let $A^* = F^c(T_0, S_1, Z(T_0),\dist (X_0, Y_0)^{\beta_1})$
and assume that $\beta_1\in(0,1)$. Recall the estimates for the
probabilities that $N$ takes values in $[-k_1,0]$ or $[0,k_1]$
from Step 1 of part (i) of the proof. Analogous estimates, the
strong Markov property applied at time $T_0$, and an argument
similar to that given in Step 3 but with $k_2$ replaced by $k_3 =
\min\{k: 2^{-k} \leq \dist (X_0, Y_0)^{\beta_1}\}$, yield for some
$\beta_4 >0$,
 $$\eqalignno{
  \bP( A^*) &\leq c_{160}
 \dist(X_0,Y_0)^{\beta_4},&(3.34)\cr
 \bP(\{-k_3 \leq N \leq 0\} \cap A^*) &\leq c_{161}
 \dist(X_0,Y_0)^{1+\beta_4} 2^{k_3},
 \cr
  \bP(\{0 < N \leq k_3\} \cap A^*) &\leq c_{162}
 \dist(X_0,Y_0)^{1+\beta_4} .
 }$$
We will estimate
 $$\bE(| \log \dist(X_{T_{17}},Y_{T_{17}})
 -\log \dist(X_0,Y_0)|
 \bone_{ A^*})
 $$
by splitting the integral into the sum of integrals over various
events, as in (3.26). An upper bound for the above expectation can
be obtained by using the same estimates as in (3.26), lines (b),
(d), and (f)-(n), and replacing estimates in lines (a), (c) and
(e) by the following estimates. By estimates similar to those in
Step 2 (a), for some $\beta_5>0$,
 $$\eqalignno{
 \bE& (| \log \dist(X_{T_{17}},Y_{T_{17}})
 -\log \dist(X_0,Y_0)|
 \bone_{\{ -k_1 \leq N \leq 0, T_1 \leq T_2\} \cap A^*})&(3.35)
 \cr
 &\leq \sum_{-k_1 \leq k \leq -k_3}
 \bP(N= k,T_1 \leq T_2)
 (c_{163} 2^{2k} + c_{164} \dist(X_0,Y_0))\cr
 &\qquad + \bP(\{-k_3 < N \leq 0\} \cap A^*)
 (c_{165} 2^{2k_3} + c_{166} \dist(X_0,Y_0))
 \cr
 &\leq \sum_{-k_1 \leq k \leq -k_3} c_{167} \dist(X_0,Y_0) 2^{-k}
 (c_{168} 2^{2k} + c_{169} \dist(X_0,Y_0))\cr
 &\qquad + c_{170} \dist(X_0,Y_0)^{1+\beta_4} 2^{k_3}
 (c_{171} 2^{-2k_3} + c_{172} \dist(X_0,Y_0))
 \cr
 &\leq c_{173} \dist(X_0,Y_0)^{1+\beta_1} + c_{174}
 \dist(X_0,Y_0)^{2-\beta_0}\cr
 &\qquad + c_{175} \dist(X_0,Y_0)^{1+\beta_4+\beta_1}
 + c_{176}
 \dist(X_0,Y_0)^{2+\beta_4-\beta_1} \cr
 & \leq c_{177} \dist(X_0,Y_0)^{1+\beta_5}.
 }$$
Similarly, by estimates similar to Step 2 (c), for some
$\beta_6>0$,
 $$\eqalignno{
 \bE& (| \log \dist(X_{T_{17}},Y_{T_{17}})
 -\log \dist(X_0,Y_0)|
 \bone_{\{ 0 < N \leq k_1, T_1 \leq T_2\}\cap A^*})&(3.36)
 \cr
 &\leq \sum_{k_3 \leq k \leq k_1}
 \bP(N=k, T_1 \leq T_2)
 (c_{178} k + c_{179} \dist(X_0,Y_0))
 \cr
 &\qquad + \bP(\{0 < N < k_3\}\cap A^*)
 (c_{180} k_3 + c_{181} \dist(X_0,Y_0))
 \cr
 &\leq \sum_{k_3 \leq k \leq k_1} c_{182} \dist(X_0,Y_0) 2^{-k/2}
 (c_{183} k + c_{184} \dist(X_0,Y_0))
 \cr
 &\qquad + c_{185} \dist(X_0,Y_0)^{1+\beta_4}
 (c_{186} k_3 + c_{187} \dist(X_0,Y_0))
 \cr
 &\leq c_{188} \dist(X_0,Y_0)^{1+\beta_4/2} |\log \dist(X_0,Y_0)|
 + c_{189} \dist(X_0,Y_0)^{2+\beta_4/2}\cr
 &\qquad
 + c_{190} \dist(X_0,Y_0)^{1+\beta_4}  |\log \dist(X_0,Y_0)|
 + c_{191} \dist(X_0,Y_0)^{2+\beta_4} \cr
 & \leq c_{192} \dist(X_0,Y_0)^{1+\beta_6}.
 }$$
Recall that $\beta_0\in(1/2,1)$. We apply estimates similar to
those in Step 2 (e) to see that for some $\beta_7>0$,
 $$\eqalignno{
 \bE& (| \log \dist(X_{T_{17}},Y_{T_{17}})
 -\log \dist(X_0,Y_0)|
 \bone_{\{  N \leq -k_1, T_1 \leq T_2\}\cap A^*})&(3.37)
 \cr
 & \leq \bP( A^*) (
 c_{193} \dist(X_0,Y_0)^{2\beta_0} + c_{194} \dist(X_0,Y_0))
 \cr
 & \leq c_{195} \dist(X_0,Y_0)^{\beta_4}(
 c_{196} \dist(X_0,Y_0)^{2\beta_0} + c_{197} \dist(X_0,Y_0))
 \cr
 & \leq c_{198} \dist(X_0,Y_0)^{1+\beta_7}.
 }$$
Note that the sum of lines (b), (d), and (f)-(n) in (3.26) is
bounded by $ c_{199} \dist(X_0,Y_0)^{ \beta_8}$ for some $c_{199}$
and $\beta_8>1$. This and (3.35)-(3.37) imply that for some
$c_{200}$ and $\beta_9>1$,
 $$\bE(| \log \dist(X_{T_{17}},Y_{T_{17}})
 -\log \dist(X_0,Y_0)|
 \bone_{ A^*}) \leq  c_{200} \dist(X_0,Y_0)^{ \beta_9}.
 $$
As in part (i), we note that in fact we have proved that
 $$\bE(\sup_{t\in[0, T_{17}]}| \log \dist(X_{t},Y_{t})
 -\log \dist(X_0,Y_0)|
 \bone_{ A^*}) \leq  c_{200} \dist(X_0,Y_0)^{ \beta_9}.
 $$
Recall from part (i) that for $k \leq 1/(4 c_{138}
\dist(X_{0},Y_{0}))$ we have $\dist(X_{Q_{k}},Y_{Q_{k}}) \leq 2
\dist(X_{0},Y_{0})$. Let
 $$A^*_k =F^c(T_0 \circ \theta_{Q_{k-1}}, S_1, Z(T_0 \circ
 \theta_{Q_{k-1}}), \dist(X_{Q_{k-1}}, Y_{Q_{k-1}})^{\beta_1})
 .$$
Then for $k \leq 1/(4 c_{138} \dist(X_{0},Y_{0}))$,
 $$\eqalign{
  \bE&(\sup_{t\in[Q_{k-1}, Q_k]}| \log \dist(X_{t},Y_{t})
 -\log \dist(X_{Q_{k-1}},Y_{Q_{k-1}})|
 \bone_{ A^*})\cr
 &\leq \bE(\sup_{t\in[Q_{k-1}, Q_k]}| \log \dist(X_{t},Y_{t})
 -\log \dist(X_{Q_{k-1}},Y_{Q_{k-1}})|
 \bone_{ A^*_k})\cr
 &\leq  c_{201} \dist(X_0,Y_0)^{ \beta_9}.
 }$$
We have the following estimate analogous to (3.27),
 $$\eqalignno{
 \bE & \left(|\log \dist(X_{S_1},Y_{S_1})
 -\log \dist(X_{0},Y_{0})|\cdot \bone\left( \bigcup_{k\geq 0}
 \{ S_1 \in [Q_k, Q_{k+1}]\}\right)\bone_{A^*}\right)
 &(3.38)\cr
 & \leq
 \sum_{k\geq 0} \bE (\bone _{\{ S_1 \in [Q_k, Q_{k+1}]\}}
   \sup_{t \in [Q_k, Q_{k+1}]}
 |\log \dist(X_{t},Y_{t})
 -\log \dist(X_{0},Y_{0})|\bone_{A^*} ) \cr
 & \leq
 \sum_{k\geq 0} \bE (\bone _{\{ S_1 \in [Q_k, Q_{k+1}]\}}
   \sup_{t \in [Q_k, Q_{k+1}]}
 |\log \dist(X_{t},Y_{t})
 -\log \dist(X_{0},Y_{0})|\bone_{A^*_{k+1}} ) \cr
 & \leq
 \sum_{k\geq 0} \bE (\bone _{\{ S_1 \geq Q_k\}}
  \sup_{t \in [Q_k, Q_{k+1}]}
 |\log \dist(X_{t},Y_{t})
 -\log \dist(X_{0},Y_{0})|\bone_{A^*_{k+1}}) \cr
 & \leq \sum_{k \leq  1/(4 c_{138} \dist(X_{0},Y_{0})) }
 2^{-k}  c_{201} k\dist(X_{0},Y_{0})^{\beta_9}
 + \sum_{k >  1/(4 c_{138} \dist(X_{0},Y_{0})) }
 2^{-k}  c_{142}k \cr
 & \leq c_{202} \dist(X_{0},Y_{0})^{\beta_9}
 + c_{203} \dist(X_{0},Y_{0})^{-1}
 \exp( - c_{204} \dist(X_{0},Y_{0})^{-1}) \cr
 & \leq c_{205} \dist(X_{0},Y_{0})^{\beta_9}.
 }$$
The expectation
 $$\bE \left( |\log \dist(X_{S_1}, Y_{S_1}) - \log \dist (X_0, Y_0)|
 \bone_{F^c(T_0, S_1, Z(T_0),\eps^{\beta_1})}\right)$$
is bounded by the estimates on the right hand sides of (3.31),
(3.32), (3.33) and (3.38). This easily implies part (ii) of the
lemma.

(iii) Let $C = \{Z_{T_0} \in K\} $. Note that $C$ signifies the
event discussed in part (n) of Step 4 in part (i) of the proof.
Recall from that step that $\bP(C) \leq c_{206}
\dist(X_0,Y_0)^{\beta_{10}}$ for some $\beta_{10}>1$. If we add
the factor $\bone_C$ to the left hand side of (3.26), the right
hand side of (3.26) is reduced to line (n), and we obtain
 $$\eqalignno{
 \bE& \left( | \log \dist(X_{T_{17}},Y_{T_{17}})
 -\log \dist(X_0,Y_0)| \bone_C\right)\cr
 &\leq \sum_{n\geq 1}
 c_{207} \dist(X_0,Y_0)^{1+\beta_0/2} (1-p_1)^{n+1}
 (c_{208} + (3 \beta_3^{2n} - 1) |\log\dist(X_0,Y_0)|)
  \cr
  &\leq c_{209} \dist(X_0,Y_0)^{\beta_{11}},
 }$$
for some $\beta_{11} >1$. We have the following formula similar to
(3.38),
 $$\eqalignno{
 \bE & \left(|\log \dist(X_{S_1},Y_{S_1})
 -\log \dist(X_{0},Y_{0})|\cdot
 \bone\left( \bigcup_{k\geq 0}
 \{ S_1 \in [Q_k, Q_{k+1}]\}\right)\bone_C \right)
 \cr
 & \leq \bE (\bone _{\{ S_1 \leq Q_{1}\}}
  |\log \dist(X_{S_1},Y_{S_1})
 -\log \dist(X_{0},Y_{0})|\bone_C) \cr
  & \quad +
 \sum_{k\geq 1} \bE (\bone _{\{ S_1 \in [Q_k, Q_{k+1}]\}}
  |\log \dist(X_{S_1},Y_{S_1})
 -\log \dist(X_{0},Y_{0})| \bone_C) \cr
 & \leq \bE (\bone _{\{ S_1 \leq Q_{1}\}}
  |\log \dist(X_{S_1},Y_{S_1})
 -\log \dist(X_{0},Y_{0})|\bone_C) \cr
  & \quad +
 \sum_{k\geq 1} \bE (\bone _{\{ S_1 \geq Q_k\}}
  \sup_{t \in [Q_k, Q_{k+1}]}
 |\log \dist(X_{t},Y_{t})
 -\log \dist(X_{0},Y_{0})|\bone_C) \cr
 & \leq  c_{209} \dist(X_0,Y_0)^{\beta_{11}}
 + \sum_{k \leq  1/(4 c_{210} \dist(X_{0},Y_{0})) }
 2^{-k}  c_{211} k \dist(X_0,Y_0)^{\beta_{10}+1}\cr
 & \quad + \sum_{k >  1/(4 c_{210} \dist(X_{0},Y_{0})) }
 2^{-k}  c_{212}k \dist(X_0,Y_0)^{\beta_{10}+1} \cr
 & \leq  c_{209} \dist(X_0,Y_0)^{\beta_{11}}
 +c_{213} \dist(X_{0},Y_{0})^{\beta_{10}}
 + c_{214} \dist(X_{0},Y_{0})^{-1+\beta_{10}}
 \exp( - c_{215} \dist(X_{0},Y_{0})^{-1}) \cr
 & \leq c_{216} \dist(X_{0},Y_{0})^{\beta_{12}},&(3.39)
 }$$
for some $\beta_{12} >1$. We can bound
 $$\bE \left( |\log\dist(X_{S_1}, Y_{S_1}) - \log \dist (X_0, Y_0)|
 \bone_{\{Z_{T_0}\in K\}}\right)$$
by the sum of the right hand sides of (3.31), (3.32), (3.33) and
(3.39). Part (iii) of the lemma follows.
 \qed

\bigskip
\noindent{\bf 4. Arguments based on excursion theory}.

We start this section with a review of the excursion theory. See,
e.g., [M] for the foundations of the theory in the abstract
setting and [B] for the special case of excursions of Brownian
motion. Although [B] does not discuss reflected Brownian motion,
all results we need from that book readily apply in the present
context. We will use two different but closely related ``exit
systems.'' The first one, presented below, is a simple exit system
representing excursions of a single reflected Brownian motion from
$\prt D$. The second exit system, presented after Lemma 4.2, is
more complex as it encodes the information about two reflected
Brownian motions $X$ and $Y$, and stopping times $S_k$ and $U_k$.

An ``exit system'' for excursions of the reflected Brownian motion
$X$ from $\prt D$ is a pair $(L^*_t, H^x)$ consisting of a
positive continuous additive functional $L^*_t$ and a family of
``excursion laws'' $\{H^x\}_{x\in\prt D}$. We will soon show that
$L^*_t = L^X_t$. Let $\Delta$ denote the ``cemetery'' point
outside $\R^2$ and let ${\cal C}$ be the space of all functions
$f:[0,\infty) \to \R^2\cup\{\Delta\}$ which are continuous and
take values in $\R^2$ on some interval $[0,\zeta)$, and are equal
to $\Delta$ on $[\zeta,\infty)$. For $x\in \prt D$, the excursion
law $H^x$ is a $\sigma$-finite (positive) measure on $\cal C$,
such that the canonical process is strong Markov on
$(t_0,\infty)$, for every $t_0>0$, with the transition
probabilities of Brownian motion killed upon hitting $\prt D$.
Moreover, $H^x$ gives zero mass to paths which do not start from
$x$. We will be concerned only with the ``standard'' excursion
laws; see Definition 3.2 of [B]. For every $x\in \prt D$ there
exists a unique standard excursion law $H^x$ in $D$, up to a
multiplicative constant.

Excursions of $X$ from $\prt D$ will be denoted $e$ or $e_s$,
i.e., if $s< u$, $X_s,X_u\in\prt D$, and $X_t \notin \prt D$ for
$t\in(s,u)$ then $e_s = \{e_s(t) = X_{t+s} ,\,  t\in[0,u-s)\}$ and
$\zeta(e_s) = u -s$. By convention, $e_s(t) = \Delta$ for $t\geq
\zeta$, so $e_t \equiv \Delta$ if $\inf\{s> t: X_s \in \prt D\} =
t$. Let ${\cal E}_u = \{e_s: s \leq u\}$.

Let $\sigma_t = \inf\{s\geq 0: L^*_s \geq t\}$
and let $I$ be the set of left endpoints of all connected
components of $(0, \infty)\setminus \{t\geq 0: X_t\in \partial
D\}$.
 The following is a
special case of the exit system formula of [M],
$$\bE \left[ \sum_{t\in I} V_t \cdot f ( e_t) \right]
= \bE \int_0^\infty V_{\sigma_s}
 H^{X(\sigma_s)}(f) ds = \bE \int_0^\infty V_t H^{X_t}(f) dL^*_t,
 \eqno(4.1)
 $$
where $V_t$ is a predictable process and $f:\, {\cal
C}\to[0,\infty)$ is a universally measurable function which
vanishes on excursions $e_t$ identically equal to $\Delta$. Here
and elsewhere $H^x(f) = \int_{\cal C} f dH^x$.

The normalization of the exit system is somewhat arbitrary, for
example, if $(L^*_t, H^x)$ is an exit system and $c\in(0,\infty)$
is a constant then $(cL^*_t, (1/c)H^x)$ is also an exit system.
One can even make $c$ dependent on $x\in\prt D$. Let $\bP^y_D$
denote the distribution of Brownian motion starting from $y$ and
killed upon exiting $D$. Theorem 7.2 of [B] shows how to choose a
``canonical'' exit system; that theorem is stated for the usual
planar Brownian motion but it is easy to check that both the
statement and the proof apply to the reflected Brownian motion.
According to that result, we can take $L^*_t$ to be the continuous
additive functional whose Revuz measure is a constant multiple of
the arc length measure on $\prt D$ and $H^x$'s to be standard
excursion laws normalized so that
 $$H^x (A) = \lim_{\delta\downarrow 0} {1\over \delta} \bP_D^{x +
 \delta\n(x)} (A),\eqno(4.2)$$
for any event $A$ in a $\sigma$-field generated by the process on
an interval $[t_0,\infty)$, for any $t_0>0$.
 The normalization of the local time is linked to the normalization
of $\omega_x$, given before the statement of Theorem 1.2.

The Revuz measure of $L^X$ is the measure $dx/(2|D|)$ on $\prt D$,
i.e., if the initial distribution of $X$ is the uniform
probability measure $\mu$ in $D$ then $\bE^\mu \int_0^1 \bone_A
(X_s) dL^X_s = \int_A dx/(2|D|)$ for any Borel set $A\subset \prt
D$, see Example 5.2.2 of [FOT].

We will show that $L^*_t=L^X_t$. It is sufficient to verify that
the normalization $L^*_t=L^X_t$ works for the half-space $D_*
=\{(x_1,x_2): x_2>0\}$. Let $K_a= \{(x_1,x_2)\in \R^2: x_2 = a\}$
and $T_a= \inf\{t>0: B_{t}\in K_a\}$. Note that $\bP^{(0,0) + \delta
\n((0,0))}(T_a < T_{\prt D_*}) = \delta/a$ for $\delta<a$, so
$H^{(0,0)}(T_a < T_{\prt D_*}) = 1/a$, assuming that $H^{(0,0)}$
is normalized as in (4.2). The reflected Brownian motion $X_t$ in
the half-plane $D_*$ with $X_0 =(0,0)$, and its local time $L^X_t$
on $\prt D_*$ may be constructed from the planar Brownian motion
$B_t=(B^1_t,B^2_t)$ starting from $(0,0)$ by the following formula
 $$(X_t,L^X_t) =
 \left( \left(B^1_t, \, B^2_t -\min_{0\leq s \leq t} B^2_s \right),
 \,  -\min_{0\leq s \leq t} B^2_s \right).
$$
Note that the $y$-coordinate of an excursion of the reflecting
Brownian motion $X$ from $\partial D^*$ is just an excursion of
1-dimensional Brownian motion away from 0. It is well-known that
such 1-dimensional excursions form a Poisson point process.
The event that $X_t-(0,L^X_t)$
 does not hit $K_1$ before time $\sigma^X_1$ is the same as the
event that there is no excursion of $X$ from $\prt D_*$ such that
it starts at a time $t$ with $L^X_t = a$, $a\in[0,1]$, and the
height of the $y$-coordinate of the excursion exceeds $1+a$.
If we assume that $L^*_t = L^X_t$, then according to the exit
system formula (4.1), the probability of
 this event is equal to the probability that a Poisson random
variable with parameter $\int_0^1 {1 \over 1+a} da$ takes value 0,
i.e., this probability is equal to
 $$\exp\left( - \int_0^1 {1 \over 1+a} da\right) = 1/2.\eqno(4.3)$$
The event ``$X_t-(0,L^X_t)$ does not hit $K_1$ before time
$\sigma^X_1$'' is the same as the event ``$B_t$ does not hit $K_1$
before hitting $K_{-1}$,'' and, obviously, the last event has
probability $1/2$. This agrees with (4.3) so the assumption that
$L^*_t = L^X_t$ is correct. In other words, the normalization of
the local time $L^X_t$ contained implicitly in (1.1) and the
normalization of excursion laws $H^x$ given in (4.2) match so that
$(dL^X_t, H^x)$ is an exit system for $X_t$ from $\prt D$.

Let $T_{\prt D_*}= \inf\{t>0: B_t\in
\prt D_*\}$. Then by Theorem II.1.16 of [Ba],
  $$\bP^{(0,0) + \delta \n((0,0))}(B^1(T_{\prt D_*}) \in dy)
  = {1 \over \delta \pi (1 +(y/\delta)^2)}dy. $$
Hence,
 $$\eqalign{H^{(0,0)} (e(\zeta-) \in dy)
 &= \lim_{\delta\downarrow 0} {1\over \delta}
 \bP^{(0,0) + \delta \n((0,0))}(B^1(T_{\prt D_*}) \in dy)\cr
 &= \lim_{\delta\downarrow 0} {1\over \delta}\cdot
 {1 \over \delta \pi (1 +(y/\delta)^2)}dy
 = {1\over \pi y^2}.}$$
This means that $\omega_x(dy) = H^x(e(\zeta-) \in dy)$, and it is
easy to see that this result extends to all $C^2$-smooth domains
$D$.

 \bigskip
\noindent{\bf Lemma 4.1}. {\sl For some $c_1,c_2,c_3\in(0,\infty)$
and any $x,y\in \ol D$ and $a>0$,

(i) $\bP^{x,y}(L^X_1 > a) \leq  2 e^{-c_1 a}$, and

(ii) $\bP^{x,y}(L^Y_{\sigma^X(1)} > a) \leq c_2 e^{-c_3 a}$.

}

\bigskip

\noindent{\bf Proof}. (i) Since $D$ is a bounded $C^2$-smooth
domain in $\R^2$, it is known (cf. [BH]) that the transition
density function $p(t, x, y)$ of the reflecting Brownian motion in
$D$ satisfies the estimate
 $$ p(t, x, y)\leq c_4 t^{-1} e^{-c_5|x-y|^2/t} \qquad \hbox{for } t\in (0, 1]
\hbox{ and } x, y \in \ol D.
$$
Therefore,
$$ \sup_{x\in \ol D} \bE^x [ L^X_1]
  = c_6 \sup_{x\in \ol D} \int_{\partial D} \int_0^1 p(s, x, y) ds
 dy <\infty .
$$
Take $c_7>0$ so that $\sup_{x\in \ol D} \bE^x [ L^X_1]<1/(2c_7)$.
It follows from  Khasminskii's inequality that
$$ \sup_{x\in \ol D}
\bE^x \left[ e^{c_7 L^X_{t_0}}\right] < 2.
$$
This implies that for any $a>0$,
$$ \bP^{x, y}(L^X_1>a) < 2 e^{-c_7 a}.
$$

(ii) We have proved in part (i) that $\sup_{x\in \ol D} \bE^x
 e^{c_7 L^X_{t_0}} < 2$. This, the additivity of
$L^X$ and the Markov property of $(X, Y)$, imply that $
\bE\exp(c_7 (L^X_j - L^X_{j-1})) < 2 $.
 A routine application of the Markov property at times $t=1,2,
\dots$ shows that
 $$\eqalign{
 \bE^{x,y} \exp ( c_7 L^X_k) &= \bE^{x,y} \exp \left( c_7
 \sum_{1\leq j\leq k} (L^X_j - L^X_{j-1} )\right) \leq  2^k. }$$
It follows that if $c_8 >0$ is sufficiently small, then for
integer $k$ of the form $k = c_8 a$,
 $$\eqalign{
 \bP^{x,y}(L^X_k \geq a)
 &= \bP^{x,y}\left(\exp ( c_7 L^X_k) \geq \exp (c_7 a)\right)\cr
&\leq 2^k \exp(-c_7 a) = 2^{c_8a} \exp(-c_7 a) \leq \exp (-c_7a
/2)
 }$$

Since $y\mapsto \bP^y(L^Y_1 > 1)=\int_{ D} p(\thalf, y, z)
\bP^z(L^Y_{1/2} > 1) dz$ is continuous on $\ol D$, there is
$p_1>0$ such that $\inf_{x\in \ol D} \bP^y(L^Y_1 > 1)\geq p_1$.
Hence
 $$\bP^{x,y}(L^Y_{c_8 a} \leq 1) =
 \bP^{x,y}(L^Y_{k} \leq 1)\leq (1-p_1)^k
 =(1-p_1)^{c_8 a} = e^{-c_9 a}.$$
For $a$ of the form $a=k/c_8$, where $k\geq 1$ is an integer,
we obtain,
 $$\eqalign{\bP^{x,y}(L^X_{\sigma^Y(1)} > a)
 &\leq \bP^{x,y}(L^X_{c_8 a} > a)
 + \bP^{x,y}(\sigma^Y(1) \geq c_8 a) \cr
 &\leq \bP^{x,y}(L^X_{c_8 a} > a)
 + \bP^{x,y}(L^Y_{c_8 a} \leq 1) \cr
 & \leq e^{-c_7 a/2}
 + e^{-c_9 a}
\leq c_{10} e^{-c_{11} a}. }$$
It is elementary to see that by adjusting the values of the
constants $c_{10}$ and $c_{11}$ we can make the formula valid for
all $a\geq 0$. By interchanging the roles of $X_t$ and $Y_t$, we
obtain part (ii) of the lemma.
 \qed

\bigskip

Recall that $|D|$ and $|\prt D|$ denote the area of $D$ and the
length of its boundary.

\bigskip

\noindent{\bf Lemma 4.2}. {\sl Let $\varphi \in C(\partial D)$.
For any $\delta,p>0$ there exists $t_0< \infty$ such that for
every $x \in \ol D$, we have

\noindent (i)
 $$
 \bP^x\left(\sup_{t\geq t_0} \left | {1\over t}\int_0^t \varphi(X_s) dL^X_s
 - {1\over 2|D|}\int_{\prt D} \varphi(y)dy \right | \geq \delta\right) < p.$$

\noindent (ii) In particular,
 $$
 \bP^x\left(\sup_{t\geq t_0} \left |
{L^X_t \over t} - {|\prt D|\over 2|D|} \right | \geq \delta\right) < p.$$
 }

\bigskip
\noindent{\bf Proof}. Let $\mu $ be the uniform probability
distribution on $D$. By Lemma 4.1 (i), $\sup_{x\in\ol D}\bE^x
L^X_1 < \infty$. Since $\varphi(x)$ is bounded, it follows that
$\sup_{x\in \ol D} \bE^x |\int_0^1 \varphi(X_s) dL^X_s| < \infty$
and $\bE^\mu |\int_0^1 \varphi(X_s) dL^X_s| < \infty$. Since $\mu$
is the stationary distribution for $X$, the ergodic theorem shows
that $\bP^\mu$-a.s.,
 $$\lim_{t\to\infty} {1\over t}\int_0^t \varphi(X_s) dL^X_s
 = \lim_{k\to \infty} (1/k) \sum_{1\leq n\leq k}
 \int_{n-1}^{n} \varphi(X_s) dL^X_s
 = \bE^\mu \int_0^1 \varphi(X_s) dL^X_s.$$
The Revuz measure of $L^X$ is $dx/(2|D|)$ on $\prt D$ and $\varphi(x)$
is a continuous function, so it is easy to see that
 $$\bE^\mu \int_0^1 \varphi(X_s) dL^X_s
 ={1\over 2|D|} \int_{\prt D} \varphi(y)dy.$$
Therefore, $\bP^\mu$-a.s.,
$$\lim_{t\to\infty} {1\over t}\int_0^t \varphi(X_s) dL^X_s
={1\over 2|D|}\int_{\prt D} \varphi(y)dy,$$
 and
  $$
 \bP^\mu\left(\sup_{t\geq t_0} \left | {1\over t}\int_0^t \varphi(X_s) dL^X_s
 - {1\over 2|D|}\int_{\prt D} \varphi(y)dy \right | \geq \delta/6\right) < p/6,
 \eqno(4.4)$$
for any given $\delta ,p>0$ and some $t_0<\infty$.

The arguments on p.~6 of [BB] or Theorem 2.4 in [BH] show that
there are constants $c_1, c_2 > 0$ such that
$$ \sup_{x, y\in D} \big| p_t(x, y) - {1\over |D|}
\big| \leq c_1 \, e^{-c_2 t} \quad \hbox{ for } t\geq 1.\eqno(4.5)
$$
Let $t_1<\infty$ be so large that the right hand side of (4.5) is
$p/6$. Then the Markov property applied at $t_1$, (4.4) and (4.5)
imply that for $x\in\ol D$,
  $$\eqalignno{\bP^x&\left(\sup_{t\geq t_0} \left | {1\over  t}
 \int_{t_1}^{t_1+t} \varphi(X_s) dL^X_s
 - {1\over 2|D|}
 \int_{\prt D} \varphi(y)dy \right | \geq \delta/6\right)&(4.6)\cr
&\leq  \bP^\mu\left(\sup_{t\geq t_0} \left | {1\over t}
 \int_{t_1}^{t_1+t} \varphi(X_s) dL^X_s
 - {1\over 2|D|}\int_{\prt D} \varphi(y)dy \right | \geq \delta/6\right)
+{p\over 6} < {p\over 3}.
 }$$
We can increase $t_0$, if necessary, so that for all $x\in\ol D$,
 $$ \bP^x\left(\sup_{t\geq t_0} \left | {1\over  t}
 \int_{t_1}^{t_1+t} \varphi(X_s) dL^X_s
 -  {1\over t_1+ t} \int_{t_1}^{t_1+t} \varphi(X_s) dL^X_s \right |
 \geq \delta/3\right) \leq p/3, \eqno(4.7)$$
and
 $$\bP^x\left(\sup_{t\geq t_0}  \left | {1\over t_1+ t}
 \int_{0}^{t_1} \varphi(X_s) dL^X_s \right | \geq \delta/3\right) < p/3.
 \eqno(4.8)$$
Part (i) of the lemma follows from (4.6)-(4.8) and the triangle
inequality. Part (ii) follows from (i) by taking $\varphi =1$.
 \qed

\bigskip

We will have to analyze excursions of $X$ from $\prt D$ containing
intervals $(S_k,U_k)$. The exit system $(L^X_t, H^x)$ is
inadequate for this purpose so we will now introduce a ``richer''
version of this exit system, capable of keeping track of some
extra information.

Let $\ell(t) = \max\{S_k: k\geq 0, S_k \in[0, t]\}$. Consider the
strong Markov process $(X_t, Y_t, X_{\ell(t)}, Y_{\ell(t)} )$ and
let $e_s(t) = (X_{s+t}, Y_{s+t}, X_{\ell(s+t)}, Y_{\ell(s+t)} )$
for $t\geq 0$ and $s$ such that $X_s\in \prt D$ and $\zeta(e_s)
\df \inf\{t> s: X_t \in \prt D\} -s>0$. For all other $s$, we let
$e_s \equiv \Delta$ (a cemetery state added to $\R^8$). Note a
technical difference with the previous version of the exit
system---here, the excursions are not killed at $\zeta(e_s)$ but
are continued after that time; this version of the exit system is
discussed in Maisonneuve [M].

We will describe an exit system $(L^X_t, H^{(x_1,y_1,x_2,y_2)})$
for the process $(X_t, Y_t, X_{\ell(t)}, Y_{\ell(t)} )$ from the
set $\prt D\times \ol D^3$. For any $(x_1,y_1,x_2,y_2)\in \prt
D\times \ol D^3$, $H^{(x_1,y_1,x_2,y_2)}$ is a $\sigma$-finite
measure defined as follows. The first component of $(X_t, Y_t,
X_{\ell(t)}, Y_{\ell(t)} )$ under $H^{(x_1,y_1,x_2,y_2)}$ has the
same distribution on $[0,\zeta)$ as an excursion under $H^{x_1}$,
defined previously. Under $H^{(x_1,y_1,x_2,y_2)}$, the process
$X_t$ continues after $\zeta$ as a reflected Brownian motion in
$D$, starting from $X_{\zeta-}$, but otherwise independent of
$\{X_t, t\in[0,\zeta)\}$. The other components of $(X_t, Y_t,
X_{\ell(t)}, Y_{\ell(t)} )$ are determined by the first component
as follows. First, we find $B$, the Brownian motion driving $X$,
using the uniqueness of the solution to (1.1). Then we use $B$ and
(1.2) to define a reflected Brownian motion $Y$ in $D$ staring
from $y_1$. We set $S^e_1=U^e_0 = 0$ and for $k\geq 1$,
 $$\eqalign{
 U^e_1 &= \inf\{t> 0:
 \dist(X_t, x_2)  \lor\dist(Y_t, y_2) \geq a_1
 \dist(x_2,y_2) \},\cr
 S^e_k& = \inf\{t> U^e_{k-1}: \dist(X_t, \prt D) \lor
 \dist(Y_t, \prt D) \leq a_2 \dist(X_t, Y_t)^2\},  \cr
 U^e_k &= \inf\{t> S^e_k:
 \dist(X_t, X_{S^e_k})  \lor\dist(Y_t, Y_{S^e_k}) \geq a_1
 \dist(X_{S^e_k}, Y_{S^e_k}) \}.
 }$$
Let $\ell^e(t) = \max\{S^e_k: k\geq 1, S^e_k \in[0, t]\}$. The
last two components of the process under $H^{(x_1,y_1,x_2,y_2)}$
are defined to be $X_{\ell^e(t)}$ and $ Y_{\ell^e(t)} $ if
$\ell^e(t)>0$, and $x_2$ and $y_2$, if $\ell^e(t) = 0$.

Note that the exit system for $(X_t, Y_t, X_{\ell(t)}, Y_{\ell(t)}
)$ from $\prt D\times \ol D^3$ is equivalent, in a sense, to the
exit system of $X_t$ from $\prt D$ because $\ol D^3$ is the state
space for $(Y_t, X_{\ell(t)}, Y_{\ell(t)} )$.
Moreover, since $X$ and $Y$ are strong solutions to the stochastic
Skorokhod equations (1.1)-(1.2) driven by the same Brownian
motion, it follows that
 $(Y_t, X_{\ell(t)}, Y_{\ell(t)} )$ is a deterministic function of $\{X_s,
0\leq s \leq t\}$. We included $(Y_t, X_{\ell(t)}, Y_{\ell(t)} )$
in the process so that we can keep track of the stopping times
$S_k$ and $T_k$ inside the excursions of $X_t$ away from $\prt D$.

In the present context, the exit system formula of [M] changes its
form from that in (4.1) to
 $$\eqalignno{
 \bE \sum_{t\geq 0} V_t \cdot f ( e_t) &= \bE \int_0^\infty V_{\sigma^X_s}
 H^{(X(\sigma^X_s),Y(\sigma^X_s),X(\ell(\sigma^X_s)),Y(\ell(\sigma^X_s)))}
 (f) ds & (4.9)\cr
 &= \bE \int_0^\infty V_{t}
 H^{(X(t),Y(t),X(\ell(t)),Y(\ell(t)))}
 (f) dL^X_t ,
 }$$
where $V_t$ is a predictable process and $f$ is a non-negative
universally measurable function which vanishes on excursions
identically equal to $\Delta$, and those with $\zeta(e)=0$.

Let
 $$T^e = \inf\{ t\geq U^e_1: X_t \in \prt D \hbox{  or  }
 Y_t \in \prt D\},$$
and recall that $Z _{T^e} = X_{T^e}$ if $X_{T^e} \in \prt D$ and
$Z _{T^e} = Y_{T^e}$ otherwise. Recall also that $F(s,u, x, a) $
denotes $ \{\sup_{s\leq t \leq u} \dist(X_t,x) \leq a\}$. Let
$T_X(A) =\inf\{t\geq 0: X_t\in A\}$.

We assume that the constant $a_1$ in the next lemma satisfies
Lemma 3.6.

\bigskip
\noindent{\bf Lemma 4.3}. {\sl (i) For any $c_0>0$ there exist
$c_1<\infty$ and $\eps_0>0$ such that if $x_1 \in \prt D$,
$y_1,x_2, y_2 \in \ol D$, $\dist(x_1, x_2) \lor\dist(y_1, y_2) <
a_1 \dist(x_2, y_2) $, $\dist(x_1, y_1) \geq c_0 \dist(x_2, y_2)
$, $\dist(x_1,y_1)\leq \eps \leq \eps_0$, and $|\pi/2-
\angle(y_1-x_1, \n(x_1))| \leq c_0^{-1} \dist(x_1, y_1)$, then
 $$H^{(x_1,y_1,x_2,y_2)}\left(
 | \log \rho_{U^e_1} - \log \rho_{S^e_2}|
 \mid U^e_1 \leq \zeta \right) \leq c_1\eps
 .$$

(ii) For any $c_0,\beta_1>0$, there exist $c_1<\infty$ and
$\eps_0>0$ such that if $x_1 \in \prt D$, $y_1,x_2, y_2 \in \ol
D$, $\dist(x_1, x_2)  \lor\dist(y_1, y_2) < a_1
 \dist(x_2, y_2) $, $\dist(x_1, y_1) \geq c_0
 \dist(x_2, y_2) $,  $\dist(x_1,y_1)
\leq \eps \leq \eps_0$, $k\geq 1$, and $|\pi/2- \angle(y_1-x_1,
\n(x_1))| \leq c_0^{-1} \dist(x_1, y_1)$, then
 $$
 H^{(x_1,y_1,x_2,y_2)}\big(
 | \log \rho_{U^e_k} - \log \rho_{S^e_{k+1}}|
 \bone_{\{U^e_k \leq \tau^+( \eps^{\beta_1})\land \zeta\}}
 \mid U^e_1 \leq \zeta \big) \leq c_1
 \eps^{1+(k-1)\beta_1}
 .$$

(iii) For any $c_0>0$ there exist $\beta_1>0$, $\beta_2 >1$, and
$\eps_0>0$ such that if $x_1 \in \prt D$, $y_1,x_2, y_2 \in \ol
D$, $\dist(x_1, x_2)  \lor\dist(y_1, y_2) < a_1
 \dist(x_2, y_2) $, $\dist(x_1, y_1) \geq c_0
 \dist(x_2, y_2) $,  $\dist(x_1,y_1)
\leq \eps\leq \eps_0$, and $|\pi/2- \angle(y_1-x_1, \n(x_1))| \leq
c_0^{-1} \dist(x_1, y_1)$, then
 $$H^{(x_1,y_1,x_2,y_2)}\left(
 | \log \rho_{U^e_1} - \log \rho_{S^e_2}|
 \bone_{\{T^e \leq S^e_2\}}
 \bone_{F^c(T^e, S^e_2, Z _{T^e}, \eps^{\beta_1})}
 \mid U^e_1 \leq \zeta \right) \leq \eps^{\beta_2}
 .$$

(iv) For any $c_0>0$ there exist $\beta_1>0$, $\beta_2 >1$, and
$\eps_0>0$ such that if $x_1 \in \prt D$, $y_1,x_2, y_2 \in \ol
D$, $\dist(x_1, x_2)  \lor\dist(y_1, y_2) < a_1
 \dist(x_2, y_2) $, $\dist(x_1, y_1) \geq c_0
 \dist(x_2, y_2) $,  $\dist(x_1,y_1)
\leq \eps\leq \eps_0$, and $|\pi/2- \angle(y_1-x_1, \n(x_1))| \leq
c_0^{-1} \dist(x_1, y_1)$, then
 $$H^{(x_1,y_1,x_2,y_2)}\left(
 | \log \rho_{U^e_1} - \log \rho_{S^e_2}|
 \bone_{\{S^e_2 \leq T^e \}}
 \bone_{F^c(S^e_2, T^e, Z _{T^e}, \eps^{\beta_1})}
 \mid U^e_1 \leq \zeta \right) \leq \eps^{\beta_2}
 .$$

(v) For any $c_0,\beta_0>0$ there exist $\beta_1>0$, $\beta_2 >1$
and $\eps_0>0$ such that if $x_1 \in \prt D$, $y_1,x_2, y_2 \in
\ol D$, $\dist(x_1, x_2)  \lor\dist(y_1, y_2) < a_1
 \dist(x_2, y_2) $, $\dist(x_1, y_1) \geq c_0
 \dist(x_2, y_2) $,  $\dist(x_1,y_1)
\leq \eps\leq \eps_0$, and $|\pi/2- \angle(y_1-x_1, \n(x_1))| \leq
c_0^{-1} \dist(x_1, y_1)$, then
 $$H^{(x_1,y_1,x_2,y_2)}\left(
 | \log \rho_{U^e_1} - \log \rho_{S^e_2}|
 \bone_{\{S^e_2 \leq \tau^+(\eps^{\beta_0})\}}
 \bone_{F^c(T^e, T_X(\prt D), Z _{T^e}, \eps^{\beta_1})}
 \mid U^e_1 \leq \zeta \right) \leq \eps^{\beta_2}
 .$$

(vi) Let $K = K(x_1, \eps, \beta_1)= \{x\in \prt D: \tan
\alpha(x_1, x) \geq \eps^{-\beta_1}\}$. For any $c_0>0$ there
exist $c_1<\infty$, $\beta_1>0$, $\beta_2 >1$ and $\eps_0>0$ such
that if $x_1 \in \prt D$, $y_1,x_2, y_2 \in \ol D$, $\dist(x_1,
x_2) \lor\dist(y_1, y_2) < a_1
 \dist(x_2, y_2) $, $\dist(x_1, y_1) \geq c_0
 \dist(x_2, y_2) $,  $\dist(x_1,y_1)
\leq \eps\leq \eps_0$, and $|\pi/2- \angle(y_1-x_1, \n(x_1))| \leq
c_0^{-1} \dist(x_1, y_1)$, then
 $$H^{(x_1,y_1,x_2,y_2)}\left(
 | \log \rho_{U^e_1} - \log \rho_{S^e_2}|
 \bone_{\{Z_{T^e}\in K\}} \mid U^e_1 \leq \zeta \right)\leq c_1
 \eps^{\beta_2}.$$

(vii) Let $K$ be defined as in (vi). For any $c_0>0$ there exist
$\beta_0, \beta_1>0$, $\beta_2 >1$ and $\eps_0>0$ such that if
$x_1 \in \prt D$, $y_1,x_2, y_2 \in \ol D$, $\dist(x_1, x_2)
\lor\dist(y_1, y_2) < a_1
 \dist(x_2, y_2) $, $\dist(x_1, y_1) \geq c_0
 \dist(x_2, y_2) $,  $\dist(x_1,y_1)
\leq \eps\leq \eps_0$, and $|\pi/2- \angle(y_1-x_1, \n(x_1))| \leq
c_0^{-1} \dist(x_1, y_1)$, then
 $$H^{(x_1,y_1,x_2,y_2)}\left(
 | \log \rho_{U^e_1} - \log \rho_{S^e_2}|
 \bone_{\{S^e_2 \leq \tau^+(\eps^{\beta_0})\}}
 \bone_{\{X( T_X(\prt D))\in K\}}
 \mid U^e_1 \leq \zeta \right) \leq \eps^{\beta_2}
 .$$

 }

\bigskip
\noindent{\bf Proof}. Parts (i), (iii) and (vi) of the lemma
follow from the strong Markov property applied at $U^e_1$ and
Lemma 3.9 (i), (ii) and (iii). Part (vii) follows from (v) (proved
below) and (vi). It remains to prove (ii), (iv) and (v).

(ii) For $k\geq 2$, if $\dist (X_{S^e_{k-1}}, Y_{S^e_{k-1}}) \leq
\eps^{\beta_1}$ then both processes are within distance
$c_2\eps^{2\beta_1}$ of $\prt D$ at time $S^e_{k-1}$. Brownian
motion starting at a point $z$ at most $c_2 \eps^{2\beta_1}$ units
away from $\prt D$ will hit $\prt D$ before hitting $\prt {\cal
B}(z, \eps^{\beta_1})$ with probability no less than $1-
c_3\eps^{\beta_1}$. This and the strong Markov property applied at
$S^e_{j-1}$, $j=2, 3, \dots, k$, imply that
$$
 H^{(x_1,y_1,x_2,y_2)}\left(
 U^e_k \leq \tau^+(\eps^{\beta_1}) \land \zeta
 \mid U^e_1 \leq \zeta \right) \leq c_4
 \eps^{(k-1)\beta_1}.
$$
We combine this with part (i), using the strong Markov property at
time $U^e_k$, to see that (ii) holds.

(iv) If $S^e_2 \leq T^e$ then $\dist(X_{S^e_2}, Y_{S^e_2}) \leq
\eps$ and $\dist(X_{S_e^2}, \prt D) \leq c_1 \eps^{2}$. Suppose
that $\beta_1 \in(0, 2)$. Brownian motion starting at a point $z$
at most $c_1 \eps^{2}$ units away from $\prt D$ will hit $\prt D$
before hitting $\prt {\cal B}(z, \eps^{\beta_1})$ with probability
not less than $1- c_2\eps^{ 2-\beta_1}$. By the strong Markov
property applied at $S^e_2$,
$$
 H^{(x_1,y_1,x_2,y_2)}\left(
  F^c(S^e_2, T_X(\prt D), X_{S^e_2}, \eps^{\beta_1})
 \mid U^e_1 \leq \zeta ,
 S^e_2 \leq T^e,{\cal F}_{S^e_2}\right)
 \leq c_3\eps^{ 2-\beta_1}.
 $$
This, the strong Markov property applied at $S^e_2$ and part (i)
of the lemma imply part (iv), for a suitable choice of $\beta_1$
and $\beta_2$.

(v) We have for some $\beta_2>1$,
 $$\eqalignno{
 &H^{(x_1,y_1,x_2,y_2)}\left(
 | \log \rho_{U^e_1} - \log \rho_{S^e_2}|
  \bone_{\{\tau^+(\eps^{\beta_0}) \geq S^e_2\}}
  \bone_{F^c(T^e, T_X(\prt D), Z _{T^e}, \eps^{\beta_1})}
 \mid U^e_1 \leq \zeta \right)\cr
 &\leq H^{(x_1,y_1,x_2,y_2)}\left(
 | \log \rho_{U^e_1} - \log \rho_{S^e_2}|
  \bone_{F^c(T^e, S^e_2, Z _{T^e}, \eps^{\beta_1}/2)}
 \mid U^e_1 \leq \zeta \right)\cr
 & \quad +H^{(x_1,y_1,x_2,y_2)}\Big(
 | \log \rho_{U^e_1} - \log \rho_{S^e_2}|
 \bone_{\{S^e_2 \leq \tau^+(\eps^{\beta_0}) \land T_X(\prt D) \}}
  \bone_{F(T^e, S^e_2, Z _{T^e}, \eps^{\beta_1}/2)}\cr
  &\qquad \times
  \bone_{F^c(S^e_2, T_X(\prt D), X(S^e_2), \eps^{\beta_1}/2)}
 \mid U^e_1 \leq \zeta \Big)\cr
 & \leq \eps^{\beta_2}  +H^{(x_1,y_1,x_2,y_2)}\Big(
 | \log \rho_{U^e_1} - \log \rho_{S^e_2}|
 \bone_{\{S^e_2 \leq \tau^+(\eps^{\beta_0}) \land T_X(\prt D)
 \}}\cr
 &\qquad \times
  \bone_{F^c(S^e_2, T_X(\prt D), X(S^e_2), \eps^{\beta_1}/2)}
 \mid U^e_1 \leq \zeta \Big),&(4.10)\cr
 }$$
where the last inequality follows from (iii). It will suffice to
bound the second term on the right hand side. If $\dist(X_{S^e_2},
Y_{S^e_2}) \leq \eps^{\beta_0}$, then $\dist(X_{S_e^2}, \prt D)
\leq c_1 \eps^{2\beta_0}$. Choose $\beta_1>0$ so that $2\beta_0 -
\beta_1>0$. Brownian motion starting at a point $z$ at most $c_1
\eps^{2\beta_0}$ units away from $\prt D$ will hit $\prt D$ before
hitting $\prt {\cal B}(z, \eps^{\beta_1}/2)$ with probability not
less than $1- c_2\eps^{ 2\beta_0-\beta_1}$. By the strong Markov
property applied at $S^e_2$,
$$
 H^{(x_1,y_1,x_2,y_2)}\left(
  F^c(S^e_2, T_X(\prt D), X_{S^e_2}, \eps^{\beta_1}/2)
 \mid U^e_1 \leq \zeta ,
 S^e_2 \leq \tau^+(\eps^{\beta_0}),{\cal F}_{S^e_2}\right)
 \leq c_3\eps^{ 2\beta_0-\beta_1}.
 $$
This, the strong Markov property applied at $S^e_2$ and part (i)
of the lemma imply that the second term on the right hand side of
(4.10) is bounded by $c_4\eps^{1+ 2\beta_0-\beta_1}$. This
completes the proof of the lemma.
 \qed

\bigskip
\noindent{\bf Lemma 4.4}. {\sl (i) There exists $c_1<\infty$ such
that for $x_1 \in \prt D$, $y_1,x_2, y_2 \in \ol D$,
 $$H^{(x_1,y_1,x_2,y_2)}
 |\log\cos \alpha(x_1, X_{\zeta-}) | \leq c_1
 .$$

(ii) Let $K = K(x_1, \eps, \beta_1)=\{x\in \prt D: \tan
\alpha(x_1, x) \geq \eps^{-\beta_1}\}$. For any $c_0>0$ there
exist $ \beta_1>0$, $\beta_2>1$ and $\eps_0>0$ such that if $x_1
\in \prt D$, $y_1,x_2, y_2 \in \ol D$, $\dist(x_1, x_2)
\lor\dist(y_1, y_2) < a_1
 \dist(x_2, y_2) $, $\dist(x_1, y_1) \geq c_0
 \dist(x_2, y_2) $,  $\dist(x_1,y_1)
\leq \eps\leq \eps_0$, and $|\pi/2- \angle(y_1-x_1, \n(x_1))| \leq
c_0^{-1} \dist(x_1, y_1)$, then
 $$H^{(x_1,y_1,x_2,y_2)}\left(
  |\log\cos \alpha(x_1, X_{\zeta-}) |
 \bone_{\{X( T_X(\prt D))\in K\}}
 \mid U^e_1 \leq \zeta \right) \leq \eps^{\beta_2}
 .$$

(iii) For any $c_0>0$ there exist $\beta_1>0$, $\beta_2 >1$ and
$\eps_0>0$ such that if $x_1 \in \prt D$, $y_1,x_2, y_2 \in \ol
D$, $\dist(x_1, x_2) \lor\dist(y_1, y_2) < a_1
 \dist(x_2, y_2) $, $\dist(x_1, y_1) \geq c_0
 \dist(x_2, y_2) $,  $\dist(x_1,y_1)
\leq \eps\leq \eps_0$, and $|\pi/2- \angle(y_1-x_1, \n(x_1))| \leq
c_0^{-1} \dist(x_1, y_1)$, then
 $$H^{(x_1,y_1,x_2,y_2)}\left(
 |\log\cos \alpha(x_1, X_{\zeta-}) |
 \bone_{\{T^e\leq S^e_2\}}
 \bone_{F^c(T^e, S^e_2, Z _{T^e}, \eps^{\beta_1})}
 \mid U^e_1 \leq \zeta \right) \leq \eps^{\beta_2}
 .$$

(iv) For any $c_0>0$ there exist $\beta_1>0$, $\beta_2 >1$ and
$\eps_0>0$ such that if $x_1 \in \prt D$, $y_1,x_2, y_2 \in \ol
D$, $\dist(x_1, x_2) \lor\dist(y_1, y_2) < a_1
 \dist(x_2, y_2) $, $\dist(x_1, y_1) \geq c_0
 \dist(x_2, y_2) $,  $\dist(x_1,y_1)
\leq \eps\leq \eps_0$, and $|\pi/2- \angle(y_1-x_1, \n(x_1))| \leq
c_0^{-1} \dist(x_1, y_1)$, then
 $$H^{(x_1,y_1,x_2,y_2)}\left(
 |\log\cos \alpha(x_1, X_{\zeta-}) |
 \bone_{\{ S^e_2\leq T^e\}}
 \bone_{F^c(S^e_2, T^e, Z _{T^e}, \eps^{\beta_1})}
 \mid U^e_1 \leq \zeta \right) \leq \eps^{\beta_2}
 .$$

(v) For any $c_0>0$ there exist $\beta_1>0$, $\beta_2 >1$ and
$\eps_0>0$ such that if $x_1 \in \prt D$, $y_1,x_2, y_2 \in \ol
D$, $\dist(x_1, x_2) \lor\dist(y_1, y_2) < a_1
 \dist(x_2, y_2) $, $\dist(x_1, y_1) \geq c_0
 \dist(x_2, y_2) $,  $\dist(x_1,y_1)
\leq \eps\leq \eps_0$, and $|\pi/2- \angle(y_1-x_1, \n(x_1))| \leq
c_0^{-1} \dist(x_1, y_1)$, then
 $$H^{(x_1,y_1,x_2,y_2)}\left(
 |\log\cos \alpha(x_1, X_{\zeta-}) |
 \bone_{F^c(T^e, T_X(\prt D), Z _{T^e}, \eps^{\beta_1})}
  \mid U^e_1 \leq \zeta \right) \leq \eps^{\beta_2}
 .$$

 }

\bigskip
\noindent{\bf Proof}. (i) Recall from Section 1 that we have
assumed that the boundary of $D$ is $C^4$-smooth and that there
exist at most a finite number of points $x_1, x_2, \dots, x_n$
such that $\nu(x_k)=0$, $k=1,\dots,n$, and the third derivative of
the function representing the boundary does not vanish at any
$x_k$. This implies that there exist $\delta_0,c_2,c_3>0$ such
that if $x\in \prt D$ and $\dist(x, x_k) \leq \delta_0$ then
$|\nu(x)| \geq c_2 \dist (x,x_k)$; moreover, if $x\in \prt D$ and
$\dist(x, x_k) \geq \delta_0$ for every $k=1,\dots,n$, then
$|\nu(x)| \geq c_3$. We make $\delta_0$ smaller, if necessary, so
that we can assume that $\dist(x_j,x_k) \geq 4 \delta_0$ for all
$j\ne k$. It is elementary to see that there exists $c_4>0$ with
the following properties (a)-(c).

(a) For every $x\in \prt D$ such that $\dist(x, x_k) \geq
2\delta_0$ for $k=1,\dots,n$, and every $y\in \prt D$ with
$\dist(x,y) \leq \delta_0$, we have $|\alpha(x,y)| \geq c_4 \dist
(x,y)$.

(b) If $x\in \prt D$ and $\dist(x, x_k) < 2 \delta_0$ for some
$k$, $y\in \prt D$, $\dist(x,y) \leq \delta_0$, and $y$ lies on
the same side of $x_k$ as $x$ then $|\alpha(x,y)| \geq c_4
\dist(x,y) \dist(x,x_k)$.

(c) Suppose that $x=x_k$ for some $k$. If $y\in \prt D$ and
$\dist(x,y) \leq \delta_0$ then $|\alpha(x,y)| \geq c_4
\dist(x,y)^2$.

Make $\delta_0$ smaller, if necessary, so that for any $x,y\in
\prt D$ with $|\alpha(x,y)| \geq \pi/4$, we have $\dist(x,y) \geq
4\delta_0$.

It is standard to prove, using the same methods as in the proof of
Lemma 3.2, that for some $c_5,c_6< \infty$ and all $x\in \prt D$,
we have $\omega_x(dy) \leq c_5 \dist(x,y)^{-2}dy$ if $\dist(x,y)
\leq \delta_0$, and $\omega_x(dy) \leq c_6dy$ if $\dist(x,y) \geq
\delta_0$.

Consider any $z\in \prt D$ and let $z_1, z_2, \dots, z_m$ be all
points in $\prt D$ such that $\alpha(z,z_k) =\pi/2$. The number
$m$ of such points is bounded by a constant $m_0$ depending on $D$
but not on $z$. Let $\Gamma_1$ be the set of points on the same
connected component of $\prt D$ as $z$, within the distance
$\delta_0$ from $z$. The family of points $\{ z_1,\dots, z_m, x_1,
\dots, x_n\}$ divides $\prt D\setminus \Gamma_1$ into Jordan arcs
$\Gamma_k$, $2\leq k\leq n_0$, with $n_0 \leq n+m$. For an arc
$\Gamma_k$, let $y_k^-$ and $y_k^+$ denote its endpoints, and let
$\Gamma_k^- = \{x\in \Gamma_k: \dist(y_k^-,x) \leq \delta_0\}$,
$\Gamma_k^+ = \{x\in \Gamma_k: \dist(y_k^+,x) \leq \delta_0\}$,
and $\Gamma_k^0 = \Gamma_k \setminus (\Gamma_k^- \cup
\Gamma_k^+)$.

Since
 $$H^{(z,y_1,x_2,y_2)}
 |\log\cos \alpha(z, X_{\zeta-}) |
 = H^z |\log\cos\alpha(z,X_{\zeta-})| = \int_{\prt D} |\log\cos
 \alpha(z,y)| \omega_z(dy),
 $$
we will estimate the integral on the right hand side. We will
split the integral into the sum of integrals over $\Gamma_k$'s.
Let $\lambda$ denote the arc length measure on $\prt D$, i.e.,
$\lambda(dx)$ is an alternative and equivalent notation for $dx$.

Since $\nu(x)$ is bounded over $\prt D$, for $x\in \Gamma_1$ we
have $\alpha(x,z) \leq c_7 \dist(x,z)$ and so
$|\log\cos\alpha(x,z)| \leq c_8 \dist(x,z)^2$. Recall that
$\omega_x(dy) \leq c_5 \dist(x,y)^{-2}dy$ if $\dist(x,y) \leq
\delta_0$. This implies that
 $$ \int_{\Gamma_1} |\log\cos\alpha(z,y)| \omega_z(dy)
 \leq 2 \int_0^{\delta_0} c_8 u^2 c_5 u^{-2} du \leq c_9 \delta_0.
 \eqno(4.11)
 $$

Consider an arc $\Gamma_k$ with $k\geq 2$. First assume that the
distance between $y_k^-$ and $y_k^+$ is greater than $\delta_0$.
Since the curvature $\nu(x)$ has the constant sign on $\Gamma_k$,
the function $x\to \alpha(z,x) $ takes its maximum at one or both
endpoints of $\Gamma_k$. Recall that, by convention,
$\alpha(z,x)\leq \pi/2$. This and (a)-(c) show that
 $$\lambda(\{x\in\Gamma_k^-: \pi/2 -\alpha(x,y_k^-)
 \in [2^{-j}, 2^{-j-1}]\}) \leq c_{10} 2^{-j/2}.$$
We have $\omega_z(dx) \leq c_6dx$ for $x\in \Gamma_k$ so
 $$ \int_{\Gamma_k^-} |\log\cos\alpha(z,y)| \omega_z(dy)
 \leq c_{11}\sum_{j\geq 1} j 2^{-j/2} = c_{12},\eqno(4.12)
 $$
and similarly
 $$ \int_{\Gamma_k^+} |\log\cos\alpha(z,y)| \omega_z(dy)
 \leq  c_{12}.
 $$
By (a), for $x\in \Gamma_k^0$, $\pi/2-\alpha(x,z)\geq c_{13}>0$,
where $c_{13}$ does not depend on $z$ or $k$. Hence,
 $$ \int_{\Gamma_k^0} |\log\cos\alpha(z,y)| \omega_z(dy)
 \leq  c_{14}|\prt D| = c_{15}.
 $$

Next suppose that the distance between $y_k^-$ and $y_k^+$ is less
than $\delta_0$. Then one of the endpoints of $\Gamma_k$, say
$y_k^-$, is a point $x_{j_1}$, and $y_k^+$ is a point $z_{j_2}$.
By Lemma 3.7,
 $$\lambda(\{x\in\Gamma_k: \pi/2 -\alpha(x,y_k^-)
 \in [2^{-j}, 2^{-j-1}]\}) \leq c_{16} 2^{-j/2},$$
for $2^{-j} \leq c_{17} \dist(y^-_k,y^+_k)$. Hence,
$$ \int_{\Gamma_k} |\log\cos\alpha(z,y)| \omega_z(dy)
 \leq c_{18}\sum_{2^{-j} \leq c_{17} \dist(y^-_k,y^+_k)}
 j 2^{-j/2} \leq c_{19}. \eqno(4.13)
 $$
Since the number of arcs $n_0$ is bounded by a constant
independent of $z$, we obtain
 $$H^{(z,y_1,x_2,y_2)}
 |\log\cos \alpha(z, X_{\zeta-}) |
 \leq \sum_k \int_{\Gamma_k} |\log\cos
 \alpha(z,y)| \omega_z(dy)\leq c_{20}.
 $$

(ii) It is not hard to prove, using (4.2), that
$H^{(x_1,y_1,x_2,y_2)}(U^e_1 \leq \zeta) \geq c_2 \eps^{-1}$, so
part (i) implies that
 $$H^{(x_1,y_1,x_2,y_2)}\left(
 |\log\cos \alpha(x_1, X_{\zeta-}) |  \mid U^e_1 \leq \zeta
 \right)\leq c_3\eps.\eqno(4.14)
 $$
The proof of part (ii) can be finished by applying the same ideas
as in part (i). The only modification that is needed is to
restrict the range of $j$ in (4.12) to $2^{-j} \leq
\eps^{\beta_1}$, and similarly for (4.13).

(iii) We have
 $$\eqalignno{
 &H^{(x_1,y_1,x_2,y_2)}\left(
 |\log\cos \alpha(x_1, X_{\zeta-}) |
 \bone_{\{T^e \leq S^e_2\}}
 \bone_{F^c(T^e, S^e_2, Z _{T^e}, \eps^{\beta_1})}
 \mid U^e_1 \leq \zeta
 \right)  &(4.15)\cr
 &\leq H^{(x_1,y_1,x_2,y_2)}\left(
 |\log\cos \alpha(x_1, X_{\zeta-}) |
  \bone_{F^c(T^e, S^e_2, X _{T^e}, \eps^{\beta_1})}
 \bone_{\{T^e \leq S^e_2,
 X_{T^e} = Z_{T^e}\}}
 \mid U^e_1 \leq \zeta
 \right) \cr
 &\quad+H^{(x_1,y_1,x_2,y_2)}\Big(
 |\log\cos \alpha(x_1, X_{\zeta-}) |
  \bone_{F^c(T_X(\prt D), S^e_2, Z _{T^e}, \eps^{\beta_1}/2)}\cr
  &\qquad\times
   \bone_{\{
 Y_{T^e} = Z_{T^e}, T_X(\prt D) \leq S^e_2\}}
  \mid U^e_1 \leq \zeta
 \Big)\cr
  &\quad+H^{(x_1,y_1,x_2,y_2)}\Big(
 |\log\cos \alpha(x_1, X_{\zeta-}) |
  \bone_{F^c(T^e, T_X(\prt D), Y _{T^e}, \eps^{\beta_1}/2)}\cr
  &\qquad\times
   \bone_{\{
 Y_{T^e} = Z_{T^e}, T^e \leq T_X(\prt D) \}}
  \mid U^e_1 \leq \zeta
 \Big) .  \cr
  }$$
The first term on the right hand side is bounded by $c_1
\eps^{\beta_3}$ for some $\beta_3>1$ by (4.14), the strong Markov
property applied at $T_e(\prt D)$, and (3.34). A similar bound
holds for the second term, by the strong Markov property applied
at $T_X(\prt D)$. To estimate the last term, we fix a
$\beta_4\in(1, 1+\beta_1)$ and note that,
 $$\eqalignno{
  &H^{(x_1,y_1,x_2,y_2)}\left(
 |\log\cos \alpha(x_1, X_{\zeta-}) |
  \bone_{F^c(T^e, T_X(\prt D), Y _{T^e}, \eps^{\beta_1}/2)}
   \bone_{\{
 Y_{T^e} = Z_{T^e}, T^e \leq T_X(\prt D) \}}
  \mid U^e_1 \leq \zeta
 \right)  \cr
  &\ \leq H^{(x_1,y_1,x_2,y_2)}\left(
 |\log\cos \alpha(x_1, X_{\zeta-}) |
  \bone_{\{\dist(T_X(\prt D), x_1)\leq 2\eps^{\beta_1}\}}
  \mid U^e_1 \leq \zeta
 \right)  \cr
 &\quad + H^{(x_1,y_1,x_2,y_2)}\left(
 |\log\cos \alpha(x_1, X_{\zeta-}) |
  \bone_{\{\dist(X_{U^e_1}, \prt D)\leq \eps^{\beta_4}\}}
  \mid U^e_1 \leq \zeta
 \right)  \cr
  &\quad+H^{(x_1,y_1,x_2,y_2)}\big(
 |\log\cos \alpha(x_1, X_{\zeta-}) |
  \bone_{\{\dist(T_X(\prt D), x_1)> 2\eps^{\beta_1},
  \dist(X_{U^e_1}, \prt D) > \eps^{\beta_4}\}}\cr
  &\qquad\qquad\times
  \bone_{F^c(T^e, T_X(\prt D), Y _{T^e}, \eps^{\beta_1}/2)}
   \bone_{\{
 Y_{T^e} = Z_{T^e}, T^e \leq T_X(\prt D) \}}
  \mid U^e_1 \leq \zeta
 \big).  &(4.16)\cr
 }$$
To bound the first term on the right hand side we use the same
idea that underlies (4.11). The conditioning on $\{U^e_1 \leq
\zeta\}$ transforms the excursion law $H$ into a probability
distribution. The event in question concerns Brownian motion
starting at $X(U^e_1)$ and killed upon hitting the boundary. The
starting point, $X(U^e_1)$, is at most $c_2 \eps$ units away from
$x_1$. For points $y\in \prt D$ with $\dist(x_1, y) \in (2^{-k},
2^{-k+1}]$, we have $|\log \cos \alpha (x_1, y) | \leq c_3
2^{-2k}$, and the probability of hitting the set of such points is
bounded by $c_4 \eps 2^k$, by Lemma 3.2. Let $k_1$ be the smallest
integer such that $2^{-k_1} \leq \eps$, and let $k_2$ be the
smallest integer such that $2^{-k_2} \leq 2\eps^{\beta_1}$. Then
the first term on the right hand side of (4.16) is bounded by
 $$\sum_{k=k_1 }^{k_2} c_3 2^{-2k} c_4 \eps 2^k\leq c_5 \eps^2.$$

We turn to the second term on the right hand side of (4.16). It is
rather easy to show, using (4.14) and the same ideas as in part
(i) of the proof, that for any $r>0$,
 $$H^{(x_1,y_1,x_2,y_2)}( |\log\cos \alpha(x_1, X_{\zeta-}) |
 \mid U^e_1 \leq \zeta, \dist(X_{U^e_1}, \prt D) = r)
 \leq c_6 \eps .
 $$
It is straightforward to see that
 $$H^{(x_1,y_1,x_2,y_2)}(\dist(X_{U^e_1}, \prt D) \leq r
 \mid U^e_1 \leq \zeta)
 \leq c_7 r \eps^{-1} .
 $$
One can use these estimates to find a bound of the form $c_8
\eps^{\beta_5}$ with $\beta_5 >1$, for the second term on the
right hand side of (4.16).

To bound the third term on the right hand side of (4.16), we
condition the process $\{X_t, t\in(U^e_1, \zeta)\}$ under
$H^{(x_1,y_1,x_2,y_2)}$ on its endpoints $X_{U^e_1}$ and
$X_{T_X(\prt D)}$. The result is an $h$-process $R$, in the sense
of Doob, starting at $X_{U^e_1}$ and killed at $X_{T_X(\prt D)}$.
Consider random sets
 $$\eqalign{
 K_1 &= \{x\in D: \dist(x,\prt D) \leq c_9 \eps^{1+\beta_1},
 \dist (x, x_1) \leq \eps^{\beta_1}/2,
 \dist (x,X_{T_X(\prt D)}) \geq \eps^{\beta_1}/2 \},\cr
  K_2 &= \{x\in D: \dist(x,\prt D) \leq \eps,
 \dist (x, x_1) \geq \eps^{\beta_1}/2,
 \dist (x,X_{T_X(\prt D)}) \geq \eps^{\beta_1}/2 \} ,
 }$$
for some $c_9$. Suppose that the events in the indicator functions
in the last term on the right hand side of (4.16) hold. If
$\dist(x_1, Z_{T^e}) \leq \eps^{\beta_1}/2$ then the process $R$
must hit $K_1$, otherwise it must hit $K_2$. The following two
estimates follow from standard properties of harmonic measure. If
$\dist(X_{U^e_1}, \prt D)> \eps^{\beta_4}$ then the probability
that the $h$-process $R$ hits $K_1$ is bounded by $c_{10}
\eps^{1+\beta_1-\beta_4}$, and the probability that it hits $K_2$
is bounded by $c_{11} \eps^{1-\beta_1}$. The conditioning on
$\{U^e_1 \leq \zeta\}$ contributes a factor of $c_{12} \eps$ (see
(4.14)), so we obtain a bound $c_{13} (\eps^{2+\beta_1-\beta_4} +
\eps^{2- \beta_1})\leq c_{14} \eps^{\beta_5}$ for the third term
on the right hand side of (4.16), for some $\beta_5>1$. All terms
on the right hand side of (4.16) have bounds of this form so this
finishes the proof of part (ii) of the lemma.

(iv) We have
 $$\eqalign{
 &H^{(x_1,y_1,x_2,y_2)}\left(
 |\log\cos \alpha(x_1, X_{\zeta-}) |
 \bone_{\{ S^e_2\leq T^e\}}
 \bone_{F^c(S^e_2, T^e, Z _{T^e}, \eps^{\beta_1})}
 \mid U^e_1 \leq \zeta \right)  \cr
 &\leq H^{(x_1,y_1,x_2,y_2)}\left(
 |\log\cos \alpha(x_1, X_{\zeta-}) |
  \bone_{\{ S^e_2\leq T_X(\prt D)\}}
  \bone_{F^c(S^e_2, T_X(\prt D), X_{T_X(\prt D)}, \eps^{\beta_1})}
 \mid U^e_1 \leq \zeta
 \right) .  \cr
  }$$
This can be estimated just like the last term on the right hand
side of (4.15), i.e., (4.16). The crucial point is that the
distance from $X$ to $\prt D$ must be less than $c_1 \eps^2$ at
time $S^e_2$. The estimates of the third term on the right hand
side of (4.16) are based on the fact that the distance from $X$ to
$\prt D$ at time $T^e$ is bounded by $c_2 \eps^{1+\beta_1}$ (see
the definition of $K_1$).

(v) This estimate has been already proved in part (iii) because
the relevant expression appears as the last term on the right hand
side in (4.15).
  \qed

\bigskip
\noindent{\bf Lemma 4.5}. {\sl (i) There exist
$\beta_0,\beta_1>0$, $\beta_2 >1$ and $c_0,\eps_0>0$ such that if
$x_1 \in \prt D$, $y_1,x_2, y_2 \in \ol D$, $\dist(x_1, x_2)
\lor\dist(y_1, y_2) < a_1
 \dist(x_2, y_2) $, $\dist(x_1, y_1) \geq c_0
 \dist(x_2, y_2) $,  $\dist(x_1,y_1)
\leq \eps\leq \eps_0$, and $|\pi/2- \angle(y_1-x_1, \n(x_1))| \leq
c_0^{-1} \dist(x_1, y_1)$, then
 $$\eqalign{
 H^{(x_1,y_1,x_2,y_2)}&\Big( \left|
  \log \rho_{U^e_1} - \log \rho_{S^e_2}
  -|\log\cos \alpha(x_1, X_{\zeta-}) | \right|
  \bone_{\{S^e_2 \leq \tau^+(\eps^{\beta_0})\}} \cr
  &\quad \times
 \bone_{F(T^e, S^e_2, Z _{T^e}, \eps^{\beta_1})}
 \bone_{F(T^e, T_X(\prt D), Z _{T^e}, \eps^{\beta_1})}
 \mid U^e_1 \leq \zeta
 \Big) \leq \eps^{\beta_2}. \cr
 }$$

(ii) There exist $\beta_0>0$, $\beta_1>1$ and $c_0, \eps_0>0$ such
that if $x_1 \in \prt D$, $y_1,x_2, y_2 \in \ol D$, $\dist(x_1,
x_2) \lor\dist(y_1, y_2) < a_1 \dist(x_2, y_2) $, $\dist(x_1, y_1)
\geq c_0 \dist(x_2, y_2) $,  $\dist(x_1,y_1) \leq \eps\leq
\eps_0$, and $|\pi/2- \angle(y_1-x_1, \n(x_1))| \leq c_0^{-1}
\dist(x_1, y_1)$, then
 $$
 H^{(x_1,y_1,x_2,y_2)}\left( \left|
  \log \rho_{U^e_1} - \log \rho_{S^e_2}
  -|\log\cos \alpha(x_1, X_{\zeta-}) | \right|
   \bone_{\{S^e_2 \leq \tau^+(\eps^{\beta_0})\}}
 \mid U^e_1 \leq \zeta
 \right) \leq \eps^{\beta_1}.
 $$

 }

\bigskip
\noindent{\bf Proof}. (i) Let $z_1, z_2, \dots, z_{m_0} \in \prt
D$ be all points with $\alpha(x_1,z_k) =\pi/2$. For integer $k$,
let $M_k\subset \prt D$ be the set of all points $y$ such that
$\tan\alpha(y,x_1 ) \in [2^{-k}, 2^{-k+1}]$. Fix some $\beta_3$
such that $0< \beta_3 < \beta_1 < 1$. Let $k_0$ be the smallest
integer with $2^{-k_0} \leq \eps^{-\beta_3}$ and $K=\bigcup_{k\leq
k_0} M_{k}$. By Lemmas 4.3 (vii) and 4.4 (ii), we have for some
$\beta_0>0$, $\beta_2>1$,
$$\eqalign{
 &H^{(x_1,y_1,x_2,y_2)}\Big( \left|
  \log \rho_{U^e_1} - \log \rho_{S^e_2}
  -|\log\cos \alpha(x_1, X_{\zeta-}) | \right|\cr
  &\qquad \qquad \times
 \bone_{\{S^e_2 \leq \tau^+(\eps^{\beta_0})\}}
  \bone_{\{X(T_X(\prt D)) \in K\}}
 \mid U^e_1 \leq \zeta
 \Big)\cr
 &\ \leq H^{(x_1,y_1,x_2,y_2)}\Big(
  | \log \rho_{U^e_1} - \log \rho_{S^e_2}|
   \bone_{\{S^e_2 \leq \tau^+(\eps^{\beta_0})\}}
  \bone_{\{X(T_X(\prt D)) \in K\}}
 \mid U^e_1 \leq \zeta
 \Big)\cr
 &\quad + H^{(x_1,y_1,x_2,y_2)}\Big(
  |\log\cos \alpha(x_1, X_{\zeta-}) |
  \bone_{\{X(T_X(\prt D)) \in K\}}
 \mid U^e_1 \leq \zeta
 \Big) \cr
 &\ \leq \eps^{\beta_2}
 .}$$
It follows that it will suffice to show that for some $\beta_4>1$,
 $$\eqalign{
 &H^{(x_1,y_1,x_2,y_2)}\Big( \left|
   \log \rho_{U^e_1} - \log \rho_{S^e_2}
  -|\log\cos \alpha(x_1, X_{\zeta-}) | \right| \cr
  &\quad \times
  \bone_{\{X(T_X(\prt D)) \notin K\}}
   \bone_{F(T^e, S^e_2, Z _{T^e}, \eps^{\beta_1})}
 \bone_{F(T^e, T_X(\prt D), Z _{T^e}, \eps^{\beta_1})}
 \mid U^e_1 \leq \zeta
 \Big) \leq \eps^{\beta_4}
 .}$$

Let
 $$\eqalign{
 T_1 &=S^e_2 \land \inf\{t\geq U^e_1: Y_t \in \prt D \hbox{  and  }
 |\angle (\n (Y_t), X_t - Y_t) - \pi/2|
 \leq c_1 \dist(x_1,y_2)^{\beta_1}\}\cr
 &\quad \land \inf\{t\geq U^e_1: X_t \in \prt D \hbox{  and  }
 |\angle (\n (X_t), X_t - Y_t) - \pi/2|
  \leq c_1 \dist(x_1,y_2)^{\beta_1}\}. }$$
for some $c_1$. We will assume that events $F(T^e, S^e_2, Z
_{T^e}, \eps^{\beta_1})$ and $F(T^e, T_X(\prt D), Z _{T^e},
\eps^{\beta_1})$ hold and we will estimate $ \left| \log
\rho_{U^e_1} - \log \rho_{T_1}  -|\log\cos \alpha(x_1, X_{\zeta-})
| \right| $ using this assumption. Let $k_1$ be the smallest
integer with $\dist(x_1,M_{k_1}) \leq \eps^{\beta_1}$. If
$X(T_X(\prt D)) \in \bigcup_{k \geq k_1} M_k$ then $T_1 = T_X(\prt
D) \land T_Y(\prt D)$. Hence, $\dist(X_{U^e_1}, Y_{U^e_1}) =
\dist(X_{T_1}, Y_{T_1})$ and $|\log\cos \alpha(x_1,
X_{\zeta-})|\leq c_2 \eps^{2\beta_1}$. Thus, in this case,
 $$ \left|
   \log \rho_{U^e_1} - \log \rho_{T_1}
  -|\log\cos \alpha(x_1, X_{\zeta-}) | \right|
  \leq c_2 \eps^{2\beta_1}.$$
Now we make an extra assumption that $\beta_1 > 1/2$ and we
conclude that
$$\eqalignno{
 &H^{(x_1,y_1,x_2,y_2)}\Big( \left|
   \log \rho_{U^e_1} - \log \rho_{S^e_2}
  -|\log\cos \alpha(x_1, X_{\zeta-}) | \right|
  \bone_{\{X(T_X(\prt D)) \in \bigcup_{k \geq k_1} M_k\}} \cr
 &\quad \times
 \bone_{F(T^e, S^e_2, Z _{T^e}, \eps^{\beta_1})}
 \bone_{F(T^e, T_X(\prt D), Z _{T^e}, \eps^{\beta_1})}
 \mid U^e_1 \leq \zeta
 \Big) \leq c_2 \eps^{2\beta_1} = c_2 \eps^{\beta_5}
 ,&(4.17)}$$
for some $\beta_5 >1$.

Suppose that $X(T_X(\prt D)) \in  M_k$ with $k_1 \leq k \leq k_0$.
Note that, by assumption, $|\alpha(x_1,Z_t) -\alpha(x_1,
X_{\zeta-})| \leq c_3\eps^{\beta_1}$ for all $t\in [U^e_1, T_1]$
such that one of the processes $X$ or $Y$ is on the boundary at
time $t$. We have $(d/d\gamma) |\log\cos \gamma| =\tan \gamma$ for
$\gamma\in(0,\pi/2)$. This and easy geometry imply that the change
in $\rho_t$ on the interval $[U^e_1, T_1]$ is equal to $\log\cos
\alpha(x_1, X_{\zeta-})$ up to an additive constant bounded by
$c_4 2^{-k} \eps^{\beta_1}$, i.e.,
  $$ \left|
   \log \rho_{U^e_1} - \log \rho_{T_1}
  -|\log\cos \alpha(x_1, X_{\zeta-}) | \right|
  \leq c_4 2^{-k} \eps^{\beta_1}.$$
By the strong Markov property applied at $U^e_1$ and Lemma 3.2,
 $$H^{(x_1,y_1,x_2,y_2)}(X(T_X(\prt D)) \in M_k\mid U^e_1 \leq
 \zeta) \leq c_5 \eps 2^k, $$
so for $k_1 \leq k \leq k_0$,
$$\eqalign{
 &H^{(x_1,y_1,x_2,y_2)}\Big( \left|
   \log \rho_{U^e_1} - \log \rho_{S^e_2}
  -|\log\cos \alpha(x_1, X_{\zeta-}) | \right|
  \bone_{\{X(T_X(\prt D)) \in M_k\}} \cr
 &\quad \times
 \bone_{F(T^e, S^e_2, Z _{T^e}, \eps^{\beta_1})}
 \bone_{F(T^e, T_X(\prt D), Z _{T^e}, \eps^{\beta_1})}
 \mid U^e_1 \leq \zeta
 \Big) \leq c_6 \eps^{1+\beta_1}
 ,}$$
and
$$\eqalign{
 &H^{(x_1,y_1,x_2,y_2)}\Big( \left|
   \log \rho_{U^e_1} - \log \rho_{S^e_2}
  -|\log\cos \alpha(x_1, X_{\zeta-}) | \right|
  \bone_{\{X(T_X(\prt D)) \in \bigcup_{k_1 \leq k \leq k_0}M_k\}} \cr
 &\quad \times
 \bone_{F(T^e, S^e_2, Z _{T^e}, \eps^{\beta_1})}
 \bone_{F(T^e, T_X(\prt D), Z _{T^e}, \eps^{\beta_1})}
 \mid U^e_1 \leq \zeta
 \Big) \leq c_7 \eps^{1+\beta_1} |\log\eps|\leq c_8 \eps^{\beta_6}
 ,}$$
for some $\beta_6 >1$. This and (4.17) imply that for some
$\beta_7>1$,
$$\eqalignno{
 &H^{(x_1,y_1,x_2,y_2)}\Big( \left|
   \log \rho_{U^e_1} - \log \rho_{T_1}
  -|\log\cos \alpha(x_1, X_{\zeta-}) | \right| &(4.18)\cr
  &\quad \times
  \bone_{\{X(T_X(\prt D)) \notin K\}}
   \bone_{F(T^e, S^e_2, Z _{T^e}, \eps^{\beta_1})}
 \bone_{F(T^e, T_X(\prt D), Z _{T^e}, \eps^{\beta_1})}
 \mid U^e_1 \leq \zeta
 \Big) \leq c_9\eps^{\beta_7}
 .}$$

The following definitions assume that $X_{T^k_1} \in \prt D$. If
$Y_{T^k_1} \in \prt D$ then the roles of $X$ and $Y$ should be
interchanged in the definitions of $T^k_2, T^k_3, T^k_4, T^k_5$
and $T^{k+1}_1$. Let $T_1^1 = T_1$ and for $k\geq 1$,
 $$\eqalign{
 T^k_2 &= \inf\{t\geq T_1^k: \dist(Y_t, Y_{T^k_1}) \geq 2\dist(Y_{T^k_1},\prt
 D)\},\cr
 T^k_{3} &= \inf\{t \geq T_1^k: Y_t \in \prt D\}, \cr
 T^k_{4} &= \inf\{t \geq T_1^k: L^X_t - L^X_{T^k_1} \geq c_{10}
 \dist(Y_{T^k_1},\prt D)\},\cr
 T_{5}^k &= T^k_2 \land T^k_{3} \land T^k_{4},\cr
 T_1^{k+1} &= \inf\{t \geq T^k_5: X_t \in \prt D
 \hbox {  or  } Y_t \in \prt D\}.
 }$$
If $T_1^{k+1} \leq \tau^+(\eps^{\beta_0})$ and $\dist(X_{T_1^k},
X_{T^e}) \leq \eps^{\beta_1}$ then $|\pi/2-\angle(\n(Z_t),
X_t-Y_t)| \leq c_{11} \eps^{\beta_1}$ for $t\in [T^k_1,
T^{k+1}_1]$ and $ L^X_{T^{k+1}_1} - L^X_{T^k_1} \leq c_{12}
\dist(X_{T_1^k}, Y_{T_1^k})\eps^{\beta_1}$. The change of
$\dist(X_t, Y_t) $ on the interval $[T_1^k, T_1^{k+1}]$ is bounded
by the product of these numbers, that is $c_{13} \dist(X_{T_1^k},
Y_{T_1^k})\eps^{2\beta_1}$. This implies that the increment of
$|\log \rho_t|$ on the interval $[T_1^k, T_1^{k+1}]$ is bounded by
$c_{14}\eps^{2 \beta_1}$. By Lemma 3.3 (i), the probability that
$S^e_2 \geq T_5^k$ is bounded by $p_1^k$, for some $p_1 < 1$. We
obtain, by the strong Markov property applies at $T_1$,
$$\eqalign{
H^{(x_1,y_1,x_2,y_2)}&\Big( \left|
   \log \rho_{T_1} - \log \rho_{S^e_2} \right|
   \bone_{\{S^e_2 \leq \tau^+(\eps^{\beta_0})\}}
   \bone_{F(T^e, S^e_2, Z _{T^e}, \eps^{\beta_1})}
 \mid U^e_1 \leq \zeta
 \Big) \cr
 &\leq \sum_{k\geq 1} p_1^k c_{14}\eps^{2\beta_1}
 \leq c_{15}\eps^{2\beta_1}= c_{15} \eps^{\beta_8}
 .}$$
The exponent $\beta_8 $ is greater than $1$ provided $\beta_1>
1/2$. Part (i) of the lemma follows from the last estimate and
(4.18).

(ii) We have
 $$\eqalign{H&\left(
 \left| \log \rho_{U^e_1} - \log \rho_{S^e_2} -
 |\log\cos \alpha(x_1, X_{\zeta-}) |\right|
  \bone_{\{S^e_2 \leq \tau^+(\eps^{\beta_0})\}}
 \mid U^e_1 \leq \zeta
 \right) \cr
 & \leq H^{(x_1,y_1,x_2,y_2)}\Big( \left|
  \log \rho_{U^e_1} - \log \rho_{S^e_2}
  -|\log\cos \alpha(x_1, X_{\zeta-}) | \right|
  \bone_{\{S^e_2 \leq \tau^+(\eps^{\beta_0})\}} \cr
  &\quad \quad\times
 \bone_{F(T^e, S^e_2, Z _{T^e}, \eps^{\beta_1})}
 \bone_{F(T^e, T_X(\prt D), Z _{T^e}, \eps^{\beta_1})}
 \mid U^e_1 \leq \zeta
 \Big)\cr
 &\quad+H^{(x_1,y_1,x_2,y_2)}\left(
 | \log \rho_{U^e_1} - \log \rho_{S^e_2}|
 \bone_{\{T^e \leq S^e_2\}}
 \bone_{F^c(T^e, S^e_2, Z _{T^e}, \eps^{\beta_1})}
 \mid U^e_1 \leq \zeta \right)\cr
 &\quad +H^{(x_1,y_1,x_2,y_2)}\left(
 | \log \rho_{U^e_1} - \log \rho_{S^e_2}|
 \bone_{\{S^e_2 \leq T^e \}}
 \bone_{F^c(S^e_2, T^e, Z _{T^e}, \eps^{\beta_1})}
 \mid U^e_1 \leq \zeta \right)\cr
 &\quad +H^{(x_1,y_1,x_2,y_2)}\left(
 | \log \rho_{U^e_1} - \log \rho_{S^e_2}|
 \bone_{\{S^e_2 \leq \tau^+(\eps^{\beta_0})\}}
 \bone_{F^c(T^e, T_X(\prt D), Z _{T^e}, \eps^{\beta_1})}
 \mid U^e_1 \leq \zeta \right)\cr
 &\quad +H^{(x_1,y_1,x_2,y_2)}\left(
 |\log\cos \alpha(x_1, X_{\zeta-}) |
 \bone_{\{T^e\leq S^e_2\}}
 \bone_{F^c(T^e, S^e_2, Z _{T^e}, \eps^{\beta_1})}
 \mid U^e_1 \leq \zeta \right)\cr
 &\quad +H^{(x_1,y_1,x_2,y_2)}\left(
 |\log\cos \alpha(x_1, X_{\zeta-}) |
 \bone_{\{ S^e_2\leq T^e\}}
 \bone_{F^c(S^e_2, T^e, Z _{T^e}, \eps^{\beta_1})}
 \mid U^e_1 \leq \zeta \right)\cr
 &\quad +H^{(x_1,y_1,x_2,y_2)}\left(
 |\log\cos \alpha(x_1, X_{\zeta-}) |
 \bone_{F^c(T^e, T_X(\prt D), Z _{T^e}, \eps^{\beta_1})}
  \mid U^e_1 \leq \zeta \right) .
 }$$
Part (ii) of the lemma follows from the above formula, part (i) of
this lemma, and estimates in Lemma 4.3 (iii)-(v) and Lemma 4.4
(iii)-(v).
 \qed

\bigskip
\noindent{\bf Lemma 4.6}. {\sl For any $\beta_1\in(1,2)$ there
exist $\beta_2>0$, and $\eps_0>0$ such that if $\eps\leq \eps_0$
and $\dist(X_0,Y_0) \leq \eps$ then
 $$\bP(\dist(Y_{\sigma_1^X}, \prt D)\geq \eps^{\beta_1}) \leq
 \eps^{\beta_2}.$$
 }

\bigskip
\noindent{\bf Proof}. By Lemma 4.1 (ii), $\bP(L^Y_{\sigma^X_1} \geq
a )\leq c_1 e^{-c_2 a}$. Hence, for any $\beta_3>0$ and some
$\beta_4>0$,
 $$\bP(L^Y_{\sigma^X_1} \geq \beta_3 |\log \eps| )
 \leq c_1 \exp(-c_2 \beta_3 |\log \eps|) = \eps^{\beta_4}.$$
If the event $A_1\df \{L^Y_{\sigma^X_1} \leq \beta_3 |\log \eps|
\}$ holds then, by Lemma 3.8,
 $$\sup_{t\in[0,\sigma^X_1]}\dist(X_{t}, Y_{t}) \leq \dist(X_0, Y_0)
 \exp( c_4 (1 + \beta_3 |\log \eps|)) \leq c_5 \eps^{1-c_4\beta_3}
 = c_5 \eps^{1-\beta_5}.
 $$
Choose $\beta_3>0$ so small that we can find $\beta_6$ and
$\beta_7$ such that $\beta_5 < \beta_6 < \beta_7 < 1-\beta_5$ and
$\beta_1 = 1-\beta_5+\beta_6$.

Let $T_1 = \inf\{t\geq 0: X_t \in \prt D\}$ and $\{V_t, 0\leq t
\leq \sigma^X_1 - T_1\} \df \{X_{\sigma^X_1-t},  0\leq t \leq
\sigma^X_1 - T_1\} $. If we condition on the values of $X_{T_1}$
and $X_{\sigma_1^X}$, the process $V$ is a reflected Brownian
motion in $D$ starting from $X_{\sigma_1^X}$ and conditioned to
approach $X_{T_1}$ at its lifetime. It is easy to see that
$\bP(\dist(X_{T_1}, X_{\sigma^X_1}) \leq \eps^{\beta_6}) \leq c_6
\eps^{\beta_6}$.

Suppose that the event $A_2=\{\dist(X_{T_1}, X_{\sigma^X_1}) \geq
\eps^{\beta_6}\}$ holds. Conditional on this event, the
probability that $V$ does not spend $\eps^{\beta_7}$ units of
local time on the boundary of $\prt D$ before leaving the disc
${\cal B}(V_0, \eps^{\beta_6})$ is bounded by $\eps^{\beta_7-
\beta_6}$. Let $A_3$ be the event that $V$ spends $\eps^{\beta_7}$
or more units of local time on the boundary of $\prt D$ before
leaving the disc ${\cal B}(V_0, \eps^{\beta_6})$. If $A_3$ holds
then it is easy to see that $\dist(Y_{\sigma_1^X}, \prt D)\leq
\eps^{1-\beta_5+\beta_6}=\eps^{\beta_1}$.We have shown that
$\dist(Y_{\sigma_1^X}, \prt D)\leq \eps^{\beta_1}$ if $A_1\cap A_2
\cap A_3$ holds. Since $\bP((A_1\cap A_2 \cap A_3)^c) \leq
\eps^{\beta_4} + c_6 \eps^{\beta_6} + \eps^{\beta_7- \beta_6}$,
the lemma follows.
 \qed
\bigskip

Recall that $e_s$ denotes an excursion of $X$ from the boundary of
$D$ starting at time $s$ and let $\alpha(e_s) =\alpha
(e_s(0),e_s(\zeta-))$.

\bigskip
\noindent{\bf Lemma 4.7}. {\sl For any $\delta,p>0$ there exist
$t_0, c_0< \infty$ such that for every $x \in \ol D$, we have

(i) $$ \bP^x\left(\sum_{ e_{s} \in {\cal E}_{t_0}}
 |\log\cos \alpha(e_{s}) | \geq c_0\right) < p,$$

(ii)
 $$\bP^x\left(\sup_{u\geq t_0} \left | {1 \over u}
 \sum_{ e_{s} \in {\cal E}_u}
 |\log\cos \alpha(e_{s}) |
 - {1\over 2|D|}\int_{\prt D} \int_{\prt D}
 |\log \cos\alpha(z,y)|\omega_z(dy) dz \right | \geq \delta\right)
 < p.$$
 }

\bigskip
\noindent{\bf Proof}. (i) It suffices to show that $\sum_{ e_{s}
\in {\cal E}_{t_0}} |\log\cos \alpha(e_{s}) |$ has a finite
expectation, bounded by a constant independent of $x$. This
follows from the exit system formula (4.1), Lemma 4.1 (i) and
Lemma 4.4 (i).

(ii) Suppose that $X_0$ has the uniform distribution in $D$. Then,
by the exit system formula, and since the Revuz measure of $L^X_t$
is $dx/(2|D|)$,
 $$\eqalign{
 \bE \sum_{ e_{s} \in {\cal E}_1} |\log\cos \alpha(e_{s}) |
 & = \bE \int_0^1 H^{X(s)}(|\log \cos \alpha(e_s)|) dL^X_s \cr
 & = {1\over 2|D|}\int_{\prt D} \int_{\prt D}
 |\log \cos\alpha(z,y)|\omega_z(dy) dz.
 }$$
By Lemma 4.4 (i) and its proof, the last integral is finite.

Let $V_k =  \sum_{ e_{s} \in {\cal E}_k \setminus {\cal E}_{k-1}}
|\log\cos \alpha(e_{s}) |$. By the ergodic theorem,
 $$\eqalign{
 \lim_{u\to\infty} {1 \over u}
 \sum_{ e_{s} \in {\cal E}_u}
 |\log\cos \alpha(e_{s}) |
 &= \lim_{k\to\infty} (1/k)\sum_{1 \leq n \leq k} V_n\cr
 &= {1\over 2|D|}\int_{\prt D} \int_{\prt D}
 |\log \cos\alpha(z,y)|\omega_z(dy) dz.}
 $$

Recall from (4.5) that the transition density $p_t(x,y)$ of
reflected Brownian motion converges to $1/|D|$ exponentially fast
as $t\to \infty$, uniformly in $(x,y) \in \ol D^2$. This can be
used to finish the proof of part (ii) of the present lemma, using
the same argument as in the proof of Lemma 4.2 (i).
 \qed

\bigskip

Recall that ${\cal T}^c = \bigcup_k (U_{k}, S_{k+1}]$.

\bigskip
\noindent{\bf Lemma 4.8}. {\sl There exist $c_1, \eps_0>0$ such
that if $\dist(X_0,Y_0) \leq \eps\leq \eps_0$ then for any $t\geq
1$,
 $$\bE \int_{{\cal T}^c \cap [0,\sigma^X_t\land \tau^+_ \eps]}
 dL^X_s \leq c_1 t \eps |\log \eps|.$$
}

\bigskip
\noindent{\bf Proof}. It is easy to deduce from Lemma 3.3 that
 $$\bE (L^X_{S_{1}\land \tau^+_ \eps} - L^X_{U_0}) \leq c_2 \eps.
 \eqno(4.19)$$

Next we will estimate $(L^X_{S_{k+1}} - L^X_{U_k}) \bone_{\{U_k<
\tau^+_{\eps}\}}$. Fix some $k\geq 1$ and assume that $U_k<
\tau^+_{\eps}$. Note that $\dist(X_{U_k},\prt D) \lor
\dist(Y_{U_k},\prt D)\leq c_3\dist(X_{U_k},Y_{U_k})$. Let $T_1$ be
the first time after $U_k$ when either $X$ or $Y$ is in $\prt D$.
Let $n_0$ be the greatest integer such that $2^{-n_0}$ is greater
than the diameter of $D$ and let $n_1 $ be the least integer
greater than $|\log \dist(X_{U_k},Y_{U_k})|/\log 2$. By Lemma 3.2,
 $$\bP(\dist(X_{U_k}, X_{T_1}) \geq 2^{-n})
 \leq c_3 \dist(X_{U_k},Y_{U_k}) 2^n,
 $$
for $n_0\leq n \leq n_1$. This obviously implies that
 $$\bP(\dist(X_{U_k}, X_{T_1})
 \in[ 2^{-n}, 2^{-n+1}])
 \leq c_3 \dist(X_{U_k},Y_{U_k}) 2^n,
 $$
for $n_0\leq n \leq n_1$. Simple geometry shows that if
$\dist(X_{U_k}, X_{T_1}) \in[ 2^{-n}, 2^{-n+1}]$ and $Z_{T_1} =
X_{T_1}$ then $\dist(Y_{T_1}, \prt D) \leq c_4
\dist(X_{U_k},Y_{U_k}) 2^{-n}$, and if $Z_{T_1} = Y_{T_1}$ then
$\dist(X_{T_1}, \prt D) \leq c_4 \dist(X_{U_k},Y_{U_k}) 2^{-n}$.
Hence,
 $$\eqalign{
 \bP&(\dist(X_{T_1}, \prt D) \lor \dist(Y_{T_1}, \prt D)
 \in [c_4 \dist(X_{U_k},Y_{U_k}) 2^{-n-1},
 c_4 \dist(X_{U_k},Y_{U_k}) 2^{-n}])\cr
 &\leq c_3 \dist(X_{U_k},Y_{U_k}) 2^{n},
 }$$
for $n_0\leq n \leq n_1$, and
 $$\bE(\dist(X_{T_1}, \prt D) \lor \dist(Y_{T_1}, \prt D))
 \leq c_5 \dist(X_{U_k},Y_{U_k})^2
 |\log \dist(X_{U_k},Y_{U_k})|.
 $$
By Lemma 3.3 (ii), assuming that $\eps_0$ is small,
 $$\bE \left(L^X_{S_{k+1}} - L^X_{U_k} \mid
 U_k< \tau^+_{\eps}\right)
 \leq c_6 \dist(X_{U_k},Y_{U_k})^2
 |\log \dist(X_{U_k},Y_{U_k})|.
 $$

It is elementary to check that
 $$\bE \left(L^X_{U_{k}} - L^X_{S_k} \mid
 S_k< \tau^+_{\eps}\right)
 \geq c_7 \dist(X_{S_k}, Y_{S_k}), $$
and the conditional distribution of $L^X_{U_{k}} - L^X_{S_k}$
given $ \{S_k< \tau^+_{\eps}\}$ is stochastically bounded by an
exponential random variable with mean $c_8 \dist(X_{S_k},
Y_{S_k})$. Note that $\dist(X_{U_k},Y_{U_k}) \leq c_9
\dist(X_{S_k}, Y_{S_k})$. Hence,
 $$N_m \df \sum_{k=1}^m c_{10}\eps|\log\eps|(L^X_{U_{k}} - L^X_{S_k})
 \bone_{\{S_k< \tau^+_{\eps}\}}
 -(L^X_{S_{k+1}} - L^X_{U_k})
 \bone_{\{U_k< \tau^+_{\eps}\}}
 $$
is a submartingale with respect to the filtration ${\cal F}^*_m =
{\cal F}^{X,Y}_{S_{m+1}}$. If
 $$M=\inf\{m: \sum_{k=1}^m (L^X_{U_k}- L^X_{S_k}) \geq t\}$$
and $M_j =  M\land j$ then
 $$ \bE \sum_{k=1}^{M_j}\left( c_{10} \eps|\log\eps|
 (L^X_{U_{k}} - L^X_{S_k})
 \bone_{\{S_k< \tau^+_{\eps}\}}
 -(L^X_{S_{k+1}} - L^X_{U_k})
 \bone_{\{U_k< \tau^+_{\eps}\}}\right) \geq 0,$$
and
 $$\bE \sum_{k=1}^{M_j} (L^X_{S_{k+1}} - L^X_{U_k})
 \bone_{\{U_k< \tau^+_{\eps}\}}
 \leq
 \bE \sum_{k=1}^{M_j} c_{10} \eps|\log\eps|(L^X_{U_{k}} - L^X_{S_k})
 \bone_{\{S_k< \tau^+_{\eps}\}}.
 $$
We let $j\to\infty$ and obtain by the monotone convergence
 $$\bE \sum_{k=1}^{M} (L^X_{S_{k+1}} - L^X_{U_k})
 \bone_{\{U_k< \tau^+_{\eps}\}}
 \leq
 \bE \sum_{k=1}^{M} c_{10} \eps|\log\eps|(L^X_{U_{k}} - L^X_{S_k})
 \bone_{\{S_k< \tau^+_{\eps}\}}
 \leq c_{11} t \eps|\log\eps|.
 $$
Hence,
 $$\bE \int_{{\cal T}^c \cap [U_1,\sigma^X_t\land \tau^+( \eps)]}
 dL^X_s
 \leq \bE \sum_{k=1}^{M} (L^X_{S_{k+1}} - L^X_{U_k})
 \bone_{\{U_k< \tau^+_{\eps}\}}
 \leq c_{11} t \eps|\log\eps|.
 $$
This and (4.19) imply the lemma.
 \qed
\bigskip

Recall that $e_s$ denotes an excursion of $X$ from $\prt D$
starting at time $s$, $\alpha(e_s) =\alpha(e_s(0),e_s(\zeta-))$,
and ${\cal E}_t$ is the family of excursions $e_s$ with $s\leq t$.
See the beginning of Section 3 for the definition of $\ol \rho_t$.

\bigskip
\noindent{\bf Lemma 4.9}. {\sl Let ${\cal E}^*(t)$ be the
restriction of ${\cal E}_t$ to those excursions $e_u$ that satisfy
the condition $\sup_{s\in[u,u+\zeta(e_u))}\dist(X_s, Y_s) \leq
\eps^{\beta_0}$. For any $\beta_0 \in(0,1)$ there exist
$\beta_1\in(3/2,2)$, $\eps_0,\beta_2>0$ and $c_1<\infty$ such that
if $X_0\in \prt D$, $\eps< \eps_0$, $\dist(X_0,Y_0) \leq \eps$ and
$\dist(Y_0,\prt D) \leq \eps^{\beta_1}$ then,
 $$\bE \left | \log \ol \rho_{\sigma^X_1}  -
 \sum_{ e_{s} \in {\cal E}^*(\sigma^X_1)}
 |\log\cos \alpha(e_{s}) | \right| \leq c_1 \eps^{\beta_2}
 .$$
 }

\bigskip
\noindent{\bf Proof}. Recall the ``rich'' version of the exit
system introduced before Lemma 4.3, and the accompanying notation,
i.e., stopping times $S^e_k$ and $U^e_k$. Let $k_1$ be the
smallest (random) integer such that $S_{k_1} \geq \sigma_1^X$. We
will show that the triangle inequality yields
 $$\eqalignno{ &\left | \log  \ol \rho_{\sigma^X_1}  -
 \sum_{ e_{s} \in {\cal E}^*(\sigma^X_1)}
 |\log\cos \alpha(e_{s}) | \right|
 \leq  \sum_{ e_{s} \in {\cal E}^*(\sigma^X_1), s \in {\cal T}}
 \bone_{\{U^e_1 \geq \zeta(e_s)\}}
 |\log\cos \alpha(e_{s}) | &\cr
 &\qquad+  \sum_{ e_{s} \in {\cal E}^*(\sigma^X_1), s \in {\cal T}}
 \bone_{\{U^e_1 \leq \zeta(e_s)\}}
 \left| \log \rho_{U^e_1} - \log \rho_{S^e_2} -
 |\log\cos \alpha(e_{s}) |\right|\cr
 &\qquad+  \sum_{ e_{s} \in {\cal E}^*(\sigma^X_1)}
 \sum_{k\geq 2}
 \bone_{\{U^e_k \leq \zeta(e_s)\}}
 \left| \log \rho_{U^e_k} - \log \rho_{S^e_{k+1}} \right|\cr
 &\qquad +  \sum_{ e_{s} \in {\cal E}^*(\sigma^X_1), s \notin {\cal T}}
 |\log\cos \alpha(e_{s}) |\cr
 &\qquad +|\log \rho_{S _1} - \log \rho_0|
 +|\log \rho _{S_{k_1}} - \log \rho_{\sigma_1^X}|.
 &(4.20)}$$
We will argue that the right hand side properly accounts for all
the terms on the left hand side of the last formula. All the terms
of the sum $\sum_{ e_{s} \in {\cal E}^*(\sigma^X_1)} |\log\cos
\alpha(e_{s}) |$, appearing on the left hand side, are accounted
for on the first, second and fourth lines on the right hand side.
The quantity $\log  \ol \rho_{\sigma^X_1}$ can be represented as
the sum of $\log  \rho_{U_k} - \log  \rho_{S_{k+1}}$, for all $k$
such that $0 \leq U_k \leq S_{k+1} \leq \sigma^X_1$, except that
there are two extra terms corresponding to subintervals at the
very beginning and at the end of $[0,\sigma^X_1]$. The two extra
subintervals are accounted for on the last line of (4.20). The
intervals $[U_k, S_{k+1}]\subset [0,\sigma^X_1]$ are matched with
excursions $e_s$ in the following way. Consider a $U_k$ and find
an excursion $e_s = \{e_s(t), t\in[s, s +\zeta(e_s))\}$ such that
$U_k \in [s, s +\zeta(e_s)]$. Then $U_k$ is one of the times
$U_k^e$ for this excursion. Note that if $ U^e_1
> \zeta(e_s)$ for an excursion $e_s$ then there are no $k$ such
that $U_k \in [s, s +\zeta(e_s)]$, so we restrict the sums on the
second and third lines appropriately. We split the sums according
to whether $k=1$ or $k\geq 2$, and whether $s\in\T$ or not. The
sums on the second and third lines do not contain terms
corresponding to $U_1^e$ with $s\notin\T$. This is because if
$s\notin \T$ then $[s, S^e_2]$ is a subinterval of $[U_k,
S_{k+1}]$ with $S_{k+1} = S^e_2$. Then $U_k = U^{e^*}_j$ for some
$j$ and some excursion $e^*_u $ with $u<s$ but note that we cannot
have $j=1$ and $u\notin \T$. Hence, there is already a term
accounting for the interval $[U_k, S_{k+1}]$.

The following estimate is based on the same ideas as the proof of
Lemma 4.4 (i). If $\dist(x_1,y_1) \leq \eps$ then
 $$\eqalign{ H^{(x_1,y_1,x_2,y_2)}(\bone_{\{U^e_1 \geq \zeta(e_s)\}}
 |\log\cos \alpha(x_1, X_{\zeta-}) |)
 &\leq \int_0^\eps c_2 r^{-2} |\log\cos (c_3 r) |dr\cr
 &\leq \int_0^\eps c_2 r^{-2} c_4 ( c_3 r )^2 dr
 = c_5 \eps.}
 $$
Hence, by the exit system formula (4.9),
 $$\bE \sum_{ e_{s} \in {\cal E}^*(\sigma^X_1), s \in {\cal T}}
 \bone_{\{U^e_1 \geq \zeta(e_s)\}}
 |\log\cos \alpha(e_{s}) |
 \leq c_5 \eps. \eqno(4.21)
 $$

We have $H^{(x_1,y_1,x_2,y_2)} |\log\cos \alpha(x_1, X_{\zeta-}) |
\leq c_6$ for all $(x_1,y_1,x_2,y_2)$ by Lemma 4.4 (i) so, by the
exit system formula (4.9) and Lemma 4.8,
 $$\bE \sum_{ e_{s} \in {\cal E}^*(\sigma^X_1), s \notin {\cal T}}
 |\log\cos \alpha(e_{s}) | \leq
 c_6 \bE \int_{{\cal T}^c \cap [0,\sigma^X_1\land \tau^+(\eps^{\beta_0})]}
 dL^X_s \leq c_7  \eps^{\beta_0} |\log \eps|.\eqno(4.22)$$

We have by Lemma 4.5 (ii), for some $\beta_3 >1$,
 $$\eqalign{
 H^{(x_1,y_1,x_2,y_2)}&\left( \left|
  \log \rho_{U^e_1} - \log \rho_{S^e_2}
  -|\log\cos \alpha(x_1, X_{\zeta-}) | \right|
   \bone_{\{S^e_2 \leq \tau^+(\eps^{\beta_0})\}}
 \mid U^e_1 \leq \zeta
 \right) \cr
 &\leq c_7 \dist(x_1,y_1)^{\beta_3},
 }$$
so, by the strong Markov property applied at $U^e_1$,
 $$\eqalign{
 &H^{(x_1,y_1,x_2,y_2)} \left(
 \left| \log \rho_{U^e_1} - \log \rho_{S^e_2} -
 |\log\cos \alpha(x_1, X_{\zeta-}) |\right|
 \bone_{\{S^e_2 \leq \tau^+(\eps^{\beta_0})\}}
 \bone_{\{ U^e_1 \leq \zeta\}}
 \right) \cr
 &\leq c_7 \dist(x_1,y_1)^{\beta_3} H^{(x_1,y_1,x_2,y_2)} \left(
  U^e_1 \leq \zeta
 \right).
 }$$
This and the exit system formula (4.9) yield,
 $$\eqalignno{
 &\bE\sum_{ e_{s} \in {\cal E}^*(\sigma^X_1), s \in {\cal T}}
 \bone_{\{ U^e_1 \leq \zeta(e_s)\}}
 \left| \log \rho_{U^e_1} - \log \rho_{S^e_2} -
 |\log\cos \alpha(e_{s}) |\right|
  \bone_{\{S^e_2 \leq \tau^+(\eps^{\beta_0})\}}\cr
 & \leq
 \bE\sum_{ e_{s} \in {\cal E}^*(\sigma^X_1), s \in {\cal T}}
 c_7 \dist(X_s, Y_s)^{\beta_3}
 \bone_{\{ U^e_1 \leq \zeta(e_s)\}}.&(4.23)
 }$$
We will now estimate the right hand side of (4.23). Let
 $$M_j =\left(\eps^{\beta_3-1}(L^X_{U_{j}} - L^X_{S_j})- c_8
 \dist(X_{U_j} , Y_{U_j})^{\beta_3}\right) \bone_{\{S_j\leq
 \tau^+(\eps^{\beta_0})\}}.$$
It is not hard to check that $\bE(L^X_{U_{j}} - L^X_{S_j}\mid {\cal
F}_{S_j}) \geq c_9 \dist(X_{S_j} , Y_{S_j})$ and $\dist(X_{U_j} ,
Y_{U_j}) \leq c_{10} \dist(X_{S_j} , Y_{S_j})$. Hence for an
appropriate $c_8>0$, we have $\bE M_j >0$ and the process $N_k =
\sum_{j\leq k} M_j $ is a submartingale. Let $K = \inf\{k:
\sum_{j\leq k}(L^X_{U_{j}} - L^X_{S_j}) \geq 1\}$. It is easy to
check that $L^X_{U_{j}} - L^X_{S_j}$ is stochastically majorized
by an exponential random variable with mean $c_{11} \dist(X_{S_j}
, Y_{S_j}) \leq c_{11} \eps^{\beta_0}$. By the strong Markov
theorem applied at time $\sigma^X_1$, we have $\bE \sum_{j\leq
K}(L^X_{U_{j}} - L^X_{S_j}) \leq 1 + c_{11} \eps^{\beta_0}$. By
the optional stopping theorem we have $\bE N_{K\land n}\geq 0$ for
any fixed $n$, so
 $$\eqalign{
 c_8 \bE \sum_{j\leq K\land n}\dist(X_{U_j} , Y_{U_j})^{\beta_3}
  \bone_{\{S_j\leq \tau^+(\eps^{\beta_0})\}}
 &\leq \eps^{\beta_3-1}
 \bE \sum_{j\leq K\land n}(L^X_{U_{j}} - L^X_{S_j})
 \bone_{\{S_j\leq \tau^+(\eps^{\beta_0})\}}\cr
 &\leq \eps^{\beta_3-1} (1 + c_{11} \eps^{\beta_0}) .
 }$$
Letting $n\to \infty$, we obtain
 $$\bE \sum_{j\leq K}\dist(X_{U_j} , Y_{U_j})^{\beta_3}
  \bone_{\{S_j\leq \tau^+(\eps^{\beta_0})\}}
 \leq c_{12} \eps^{\beta_3-1}.\eqno(4.24)$$
Note that
 $$ \bE\sum_{ e_{s} \in {\cal E}^*(\sigma^X_1), s \in {\cal T}}
 c_7 \dist(X_s, Y_s)^{\beta_3}
 \bone_{\{ U^e_1 \leq \zeta(e_s)\}}
 \leq \bE \sum_{j\leq K}\dist(X_{U_j} , Y_{U_j})^{\beta_3}
  \bone_{\{S_j\leq \tau^+(\eps^{\beta_0})\}},$$
so this, (4.23) and (4.24) imply that
 $$\eqalignno{
 \bE&\sum_{ e_{s} \in {\cal E}^*(\sigma^X_1), s \in {\cal T}}
 \bone_{\{ U^e_1 \leq \zeta(e_s)\}}
 \left| \log \rho_{U^e_1} - \log \rho_{S^e_2} -
 |\log\cos \alpha(e_{s}) |\right|
  \bone_{\{S^e_2 \leq \tau^+(\eps^{\beta_0})\}}\cr
 &\leq c_{13} \eps^{\beta_3-1}.&(4.25)
 }$$

By Lemma 4.3 (ii) and the strong Markov property applied at
$U^e_1$, if $\dist(x_1,y_1) \leq \eps^{\beta_0}$ and $k\geq 2$
then
 $$
 H^{(x_1,y_1,x_2,y_2)}\big(
 | \log \rho_{U^e_k} - \log \rho_{S^e_{k+1}}|
 \bone_{\{U^e_k \leq \tau^+( \eps^{\beta_0/2})\land \zeta\}}
 \bone_{\{ U^e_1 \leq \zeta\}} \big) \leq c_{14}
 \eps^{\beta_0(k-1)/2}
 ,$$
so the exit system formula (4.9) implies
 $$\bE \sum_{ e_{s} \in {\cal E}^*(\sigma^X_1)}
 \sum_{k\geq 2}
 \bone_{\{ U^e_1 \leq \zeta(e_s)\}}
 \left| \log \rho_{U^e_k} - \log \rho_{S^e_{k+1}} \right|
 \leq \sum_{k\geq 2} c_{14} \eps^{\beta_0(k-1)/2}
 \leq c_{15} \eps^{\beta_4},\eqno(4.26)
 $$
for some $\beta_4>0$.

By Lemma 3.9 (i),
 $$ \bE |\log \rho _{S_1} - \log \rho_0|
 \leq c_{16} \eps.\eqno(4.27)
 $$

By Lemma 4.6, for some $\beta_5>0$,
 $$\bP(\dist(Y_{\sigma_1^X}, \prt D)\geq \eps^{\beta_1}) \leq
 \eps^{\beta_2}.\eqno(4.28)$$
By Lemma 3.9 (i),
 $$
 \bE\left(|\log \rho _{S_{k_1}} - \log \rho_{\sigma_1^X}|
 \bone_{\{\dist(Y_{\sigma_1^X}, \prt D)\leq \eps^{\beta_1}\}}
 \right)
 \leq c_{17} \eps^{\beta_0}.\eqno(4.29)
 $$
The lemma follows from (4.20), (4.21), (4.22), (4.25), (4.26),
(4.27), (4.28) and (4.29).
 \qed

\bigskip
\noindent{\bf Proof of Theorem 1.2}. The proof will consist of
three steps. First we are going to define some events. Then we
will estimate their probabilities and choose the values of the
parameters so that the probabilities of the events defined in Step
1 are large. Finally, we will prove that if all of the events
defined in Step 1 hold then $\dist(X_t,Y_t)$ has the asymptotic
behavior asserted in the theorem.

\medskip
\noindent {\it Step 1}. Suppose that $\Lambda(D) >0$ and fix arbitrarily
small $\delta,p >0$. Assume that $\delta < \Lambda(D) /8$.

We will define a number of events and stopping times, depending on
parameters $\eps, k_0, c_1, \dots, c_6$, whose values will be
specified later on. We will assume that all these parameters are
reals in $(0,\infty)$, except that $k_0>0$ is a (large) integer.
The constant $\eps$ will represent the initial distance between
the two Brownian particles, i.e., $\eps=\dist(X_0,Y_0)$. Let
$$\eqalignno{
 A_1 &= \left \{ 1-\delta\leq  (s/\sigma^X(s) )(2|D| /|\prt D|) \leq
 1+\delta,
 \quad \forall s\geq k_0\right\} ,\cr
 A_2 &= \left\{\sup_{t\geq k_0}
\left(  L^Y_{\sigma^X(t)}-2 t \right) \leq 0\right\}
\cap \left\{\sup_{t\leq \sigma^X(k_0)} \dist(X_t,Y_t)
 \leq c_1\right\}, \cr
 A_3&=
 \left\{\sup_{t\geq \sigma^X(k_0)} \left | {1\over t}\int_0^t \nu(X_s) dL^X_s
 - {1\over 2|D|}\int_{\prt D} \nu(y)dy \right | < \delta\right\}, \cr
 \wt A_3&=
 \left\{\sup_{t\geq \sigma^X(k_0)} \left | {1\over t}\int_0^t \nu(Y_s) dL^Y_s
 - {1\over 2|D|}\int_{\prt D} \nu(y)dy \right | < \delta\right\}, \cr
 A_4 &= \left\{\sup_{u\geq \sigma^X(k_0)} \left | {1 \over u}
 \sum_{ e_{s} \in {\cal E}_u}
 |\log\cos \alpha(e_{s}) |
 - {1\over 2|D|}\int_{\prt D} \int_{\prt D}
 |\log \cos\alpha(z,y)|\omega_z(dy) dz \right | < c_2
 \delta\right\}, \cr
 A_5 &=\left\{\left | \log \ol \rho_{\sigma^X(k_0)} -
 \sum_{ e_{s} \in {\cal E}_{\sigma^X(k_0)}}
 |\log\cos \alpha(e_{s}) | \right| \leq \sigma^X(k_0)
 \delta\right\},\cr
 A_6 &= \left\{
 \left|\int_{{\cal T}^c \cap [0,\sigma^X(k_0)]}
 \nu(X_s) dL^X_s \right| \leq \sigma^X(k_0)\delta\right\}, \cr
 \wt A_6 &= \left\{
 \left|\int_{{\cal T}^c \cap [0,\sigma^X(k_0)]}
 \nu(Y_s) dL^Y_s \right| \leq \sigma^X(k_0)\delta\right\}. \cr
 }$$
For each event $A_j$ on the above list, let $A'_j$ denote the
event defined in the same way except that $k_0$ is replaced by
$k_0-1$. The same remark applies to $\wt A_j$ and $\wt A'_j$.

For integer $k\geq 0$, let
 $$\eqalign{
  A_7^k &= \left\{L^Y_{\sigma^X(k_0+k+1)} - L^Y_{\sigma^X(k_0+k)}
 \leq (c_3 \log (k+2))^2\right\},\cr
 T_k &= \sigma^X(k_0+k+1)
 \land \tau^+\Big( \eps
 \exp(- c_2 (k_0+k) (\Lambda(D)-7\delta)
 + c_4(1+ (c_3 \log (k+2))^2)) \Big),\cr
 T^*_k &= \inf\{t> \sigma^X(k_0+k+1): L^Y_t - L^Y_{\sigma^X(k_0+k)}
 \geq (c_3 \log (k+2))^2 \}\cr
 &\qquad \land \tau^+\Big( \eps
 \exp(- c_2 (k_0+k)(\Lambda(D)-7\delta)
 + c_4(1+ (c_3 \log (k+2))^2) )\Big),\cr
 }$$
 $$\eqalign{
 A_8^k &= \Bigg\{\bone_{\{\dist(X_{\sigma^X(k_0+k)}, Y_{\sigma^X(k_0+k)}) \leq
 \eps \exp(- c_2 (k_0+k)(\Lambda(D)-7\delta))\}}\cr
 &\qquad \times
 \left|\int_{{\cal T}^c \cap [\sigma^X(k_0+k),T_k]}
 \nu(X_s) dL^X_s \right| \leq c_2 \delta \Bigg\}, \cr
 \wt A_8^k &= \Bigg\{\bone_{\{\dist(X_{\sigma^X(k_0+k)}, Y_{\sigma^X(k_0+k)}) \leq
 \eps \exp(- c_2 (k_0+k)(\Lambda(D)-7\delta))\}}\cr
 &\qquad \times
 \left|\int_{{\cal T}^c \cap [\sigma^X(k_0+k),T^*_k]}
 \nu(Y_s) dL^Y_s \right| \leq c_2 \delta \Bigg\},\cr
 A_9^k &= \{\dist(X_{\sigma^X(k_0+k)}, Y_{\sigma^X(k_0+k)})>
 \eps \exp(- c_2 (k_0+k)(\Lambda(D)-7\delta))\} \cr
 &\cup
 \{\dist( Y_{\sigma^X(k_0+k+1)},\prt D) \leq
 \eps^{c_5} \exp(- c_5 c_2 (k_0+k)(\Lambda(D)-7\delta))\}.
 }$$
Let ${\cal E}^{k,\eps}(s,t)$ be the subset of ${\cal
E}_{t}\setminus {\cal E}_s$ consisting of these excursions $e_u$
that satisfy
 $$\sup_{v\in[u, u+\zeta(e_u))}\dist(X_v,Y_v) \leq \eps^{c_6}
 \exp(- c_6 c_2 (k_0+k)(\Lambda(D)-7\delta)), $$
and let
$$\eqalign{
 A_{10}^k &=\Bigg\{
 \prod_{0\leq j \leq k}
 \bone_{\{\dist(X_{\sigma^X(k_0+j)}, Y_{\sigma^X(k_0+j)}) \leq
 \eps \exp(- c_2 (k_0+j)(\Lambda(D)-7\delta))\}} \cr
 & \qquad \times
 \left | \log {\ol \rho_{\sigma^X(k_0+k+1)}
 \over \ol \rho _{\sigma^X(k_0+k)}} -
 \sum_{ e_s \in {\cal E}^{k,\eps}
 (\sigma^X(k_0+k), \sigma^X(k_0+k+1))}
 |\log\cos \alpha(e_{s}) | \right| \leq c_2\delta\Bigg\}.
 }$$

\medskip
\noindent {\it Step 2}. In this step, we will choose the
parameters so that all events defined in the previous step have
large probabilities. All the bounds on probabilities will hold
uniformly for all starting points $(X_0,Y_0)\in\ol D^2$ with
$\dist(X_0,Y_0) =\eps$. The starting points will not be reflected
in the notation. We can assume that $\eps>0$ is arbitrarily small
in view of Lemma 3.1.

First we use Lemma 4.2 (ii) to choose $k_0$ so large that
$\bP(A_1)>1- p$. Let $c_2\in(0,1)$ be such that if $A_1$ holds then
for all $k\geq 0$,
 $$ c_2 (k_0+k) \leq \sigma^X(k_0+k). \eqno(4.30)$$

We choose $c_1$ to be the constant $\eps_0$ of Lemma 3.6, assuming
that the constant $c_1$ in that lemma takes the value
$c_2\delta/3$.

Lemma 4.2 is stated for the process $X_t$ but it applies equally
to $Y_t$. It is an easy consequence of part (ii) of that lemma
applied to both $X_t$ and $Y_t$ that if we enlarge $k_0$, if
necessary, then $\bP \left(\sup_{t\geq k_0} ( L^Y_{\sigma^X(t)}-2t
) \leq 0 \right)>1- p$. It follows from Lemma 3.8 that if
$L^Y_{\sigma^X(k_0)} \leq 2 k_0$ and $\eps$ is sufficiently small
then $\sup_{t\leq \sigma^X(k_0)} \dist(X_t,Y_t) \leq c_1$. Hence,
for sufficiently small $\eps$, $\bP(A_2)>1- p$.

Using Lemma 4.2 (i) we can find $t_1 $ so large that
$$\bP\left(\sup_{t\geq t_1} \left | {1\over t}\int_0^t \nu(X_s) dL^X_s
 - {1\over 2|D|}
 \int_{\prt D} \nu(y)dy \right | < \delta\right)> 1-p/2.\eqno(4.31)$$
By part (ii) of the same lemma, we can enlarge $k_0$, if
necessary, so that $\bP(\sigma^X(k_0)> t_1)> 1-p/2$. Hence, for this
value of $k_0$ we have $\bP(A_3) >1- p$. Since (4.31) holds with $X$
replaced by $Y$, our argument shows that $\bP(\wt A_3)>1- p$ for the
same value of $k_0$.

Enlarge $k_0$, if necessary, so that $\bP(\sigma^X(k_0)> t_0)>
1-p/2$, where $t_0$ is the constant in the statement of Lemma 4.7
(ii), assuming that in that lemma $p$ is replaced with $p/2$ and
$\delta$ is replaced with $c_2\delta$. Then it is easy to check,
using Lemma 4.7 (ii), that $\bP(A_4)>1- p$ holds with this choice of
$k_0$.

We will next show that with an appropriate choice of the
parameters, $\bP(A_5)>1- p$. Suppose that $\beta_0>0$ and $\beta_1
\in (3/2,2)$ so that we can apply Lemma 4.9 with these parameters.
Recall the notation ${\cal E}^*(t)$ from that lemma. Let
 $$\eqalign{
 F_k &=\{\dist(X_{\sigma^X_k},Y_{\sigma^X_k})\leq
 \eps_1\} \cap \{\dist(Y_{\sigma^X_k}, \prt D)
 \leq \eps_1^{\beta_1}\},\cr
 G_k^* & = \left\{ \left | \log
 (\ol \rho_{\sigma^X_{k+1}} /\ol \rho_{\sigma^X_{k}})-
 \sum_{ e_{s} \in
 {\cal E}^*(\sigma^X_{k+1})\setminus {\cal E}^*(\sigma^X_{k})}
 |\log\cos \alpha(e_{s}) | \right|
 \leq \sigma^X(k_0) \delta/(2k_0)\right\},\cr
  G_k & = \left\{ \left | \log
 (\ol \rho_{\sigma^X_{k+1}} /\ol \rho_{\sigma^X_{k}})-
 \sum_{ e_{s} \in
 {\cal E}(\sigma^X_{k+1})\setminus {\cal E}(\sigma^X_{k})}
 |\log\cos \alpha(e_{s}) | \right|
 \leq \sigma^X(k_0) \delta/(2k_0)\right\}.\cr}$$
Choose $\eps_1>0$ so small that for $k=0,1, \dots , k_0-1$, using
Lemma 4.9 and the strong Markov property at $\sigma^X_{k}$, $\bP(
G^*_k \mid F_k) \geq 1-p/(4k_0)$. This implies that $\bP(G^*_k) \geq
1 - p/(4k_0) - \bP(F^c_k)$.

By Lemma 4.6, we can find $\eps_2>0$ so small that for
$k=1,\dots,k_0-1$, the conditional probability of  $
\{\dist(Y_{\sigma^X_k}, \prt D) \leq \eps_1^{\beta_1}\}$ given
$\{\dist(X_{\sigma^X_{k-1}},Y_{\sigma^X_{k-1}})\leq \eps_2\}$ is
greater than $1-p/(8k_0)$. By Lemmas 3.8 and 4.1 (ii), we can make
$\eps$ so small that
$\dist(X_{\sigma^X_{k-1}},Y_{\sigma^X_{k-1}})\leq \eps_1\land
\eps_2$ for $k=1,\dots,k_0-1$ with probability greater than
$1-p/(8k_0)$. With this choice of $\eps$ we have $\bP(G^*_k) \geq 1
- p/(2k_0) $ for $k=0,\dots,k_0-1$, so $\bP(\bigcup_{0\leq k \leq
k_0-1} G^*_k) \geq 1-p/2$.

It follows from Lemmas 3.8 and 4.1 (ii) that if $\eps$ is
sufficiently small then $\bP(\tau^+(k_0^{\beta_1}) \geq k_0) \geq 1
- p/2$. If $\tau^+(k_0^{\beta_1}) \geq k_0$ then ${\cal E}^*(k_0)
= {\cal E}(k_0)$, so we obtain $\bP(\bigcup_{0\leq k \leq k_0-1}
G_k) \geq 1-p$. It is easy to check that $ \bigcup_{0\leq k \leq
k_0-1} G_k \subset A_5$, so with our choice of parameters $\eps$
and $k_0$, $\bP(A_5)>1- p$.

Recall that $\nu^* = \sup_{x\in \prt D} |\nu(x)|<\infty$. Lemma
4.8 implies that for some $\eps_1>0$, $C_1<\infty$, and $\eps \leq
\eps_0\leq \eps_1$,
 $$\bE \int_{{\cal T}^c \cap [0,\sigma^X(k_0)\land \tau^+( \eps_0)]}
 |\nu(X_s)| dL^X_s  \leq \nu^* C_1 k_0 \eps_0 |\log \eps_0|.$$
It follows that,
 $$\bP\left(
 \int_{{\cal T}^c \cap [0,\sigma^X(k_0)\land \tau^+( \eps_0)]}
 |\nu(X_s)| dL^X_s  \geq c_2 k_0 \delta\right)
 \leq  \nu^* C_1 k_0 \eps_0 |\log \eps_0|/(c_2 k_0 \delta).\eqno(4.32)$$
According to (4.30), if $A_1$ holds then $c_2 k_0 \delta \leq
\sigma^X(k_0)\delta$. Suppose that $\eps_0$ is so small that we
have $\nu^* C_1 k_0 \eps_0 |\log \eps_0|/(c_2 k_0 \delta) < p/2$.
We have shown that $\bP(\sup_{t\leq \sigma^X(k_0)} \dist(X_t,Y_t)
\leq c_1) \geq \bP(A_2) \geq 1-p$ if $\eps$ is small. The same
argument applied with $\eps_0$ in place of $c_1$ shows that we can
choose $\eps$ so small that $\bP(\tau^+(\eps_0)> \sigma^X(k_0))>
1- p/2$. This and (4.32) imply that $\bP(A_1\cap A_6)>1- 2p$.

We make $k_0$ larger, if necessary, so that by Lemma 4.2 (ii) we
have $\bP(\sigma^Y(2k_0)\leq \sigma^X(k_0)) \leq p/6$. By the same
lemma, we can choose $C_2$ so that $\bP(C_2 k_0 \delta \leq
\sigma^X(k_0)\delta) \geq 1-p/6$. Then we can make $\eps$ so small
that the same argument that leads to (4.32) gives
 $$\bP\left(
 \int_{{\cal T}^c \cap [0,\sigma^Y(2k_0)\land \tau^+( \eps_0)]}
 |\nu(Y_s)| dL^Y_s  \geq C_2 k_0 \delta\right)
 \leq  p/6.$$
Recall that $\bP(\tau^+(\eps_0)> \sigma^X(k_0))> 1- p/2$. Combining
all these estimates, we obtain $\bP(A_1\cap A_6 \cap \wt A_6)>1-
3p$.

Recall the definition of events $A'_j$ and $\wt A'_j$. By
enlarging $k_0$, if necessary, and making $\eps$ smaller, we
obtain the same estimates for events $A'_j$ and $\wt A'_j$ as for
$A_j$ and $\wt A_j$, for example, $\bP(A'_1\cap A'_6 \cap \wt
A'_6)>1- 3p$.

By Lemma 4.1 (ii), for some $C_3,C_{4}\in(0,\infty)$,
$$\bP(A^k_7) \leq C_3 \exp(-C_4 (c_3\log(k+2))^2).$$
We choose $c_3$ so large that $\sum_{k=0}^\infty \bP((A_7^k)^c) \leq
p$.

Let $c_4$ be the constant called $c_1$ in Lemma 3.8.

We make $\eps$ smaller, if necessary, so that
 $$\gamma
 \df \eps
 \exp(- c_2 (k_0+k) (\Lambda(D)-7\delta)
 + c_4(1+ (c_3 \log (k+2))^2))
 $$
is smaller than the constant $\eps_0$ in Lemma 4.8, for $k\geq 0$.
Let $C_5$ be the constant $c_1$ of Lemma 4.8. Using the strong
Markov property at $\sigma^X(k_0+k)$, and applying Lemma 4.8, we
see that for $k\geq 0$,
 $$\eqalign{
 \bP((A^k_8)^c) &\leq {1 \over c_2\delta}
 \bE \left|\int_{{\cal T}^c \cap [\sigma^X(k_0+k),T_k]}
 \nu(X_s) dL^X_s \right|\cr
 &\leq {1 \over c_2\delta} \nu^*
 \bE \int_{{\cal T}^c \cap [\sigma^X(k_0+k),T_k]}
 dL^X_s \cr
 &\leq {1 \over c_2\delta} \nu^* C_5 \gamma |\log \gamma|.}
$$
We make $\eps$ smaller, if necessary, so that $\sum_{k=0}^\infty
\bP((A_8^k)^c) \leq p$.

We apply Lemma 4.8 to $Y_t$ in place of $X_t$ to see that for
$k\geq 0$,
 $$\eqalign{
 \bP((\wt A^k_8)^c) &\leq {1 \over c_2\delta}
 \bE \left|\int_{{\cal T}^c \cap [\sigma^X(k_0+k),T^*_k]}
 \nu(Y_s) dL^Y_s \right|\cr
 &\leq {1 \over c_2\delta} \nu^*
 \bE \left|\int_{{\cal T}^c \cap [\sigma^X(k_0+k),T^*_k]}
 dL^Y_s \right|\cr
 &\leq {1 \over c_2\delta} \nu^* C_5 (c_3\log (k+2))^2
 \gamma |\log \gamma|.}
$$
We make $\eps$ smaller, if necessary, so that $\sum_{k=0}^\infty
\bP((\wt A_8^k)^c) \leq p$.

We choose $c_5\in(3/2,2)$ and $C_6$ such that, by Lemma 4.6,
 $$\bP((A_9^k)^c)\leq  \eps ^{C_6}
 \exp(- C_6 c_2 (k_0+k)(\Lambda(D)-7\delta)).
 $$
By making $\eps$ smaller, if necessary, we obtain
$\sum_{k=-1}^\infty \bP((A_9^k)^c)  \leq p/2$. Note that the
summation index starts from $k=-1$, not $k=0$ (obviously, we can
assume that $k_0 >1$).

Choose a $c_6\in(0,1)$. By Lemma 4.9, for some $C_7,C_8$, and
$k\geq 0$,
 $$\bP((A_{10}^k)^c \mid A_9^{k-1}) \leq {1 \over c_2\delta} C_7
 \eps^{C_8}
 \exp(- C_8 c_2 (k_0+k) (\Lambda(D)-7\delta)).
 $$
We make $\eps$ smaller, if necessary, so that $\sum_{k=0}^\infty
\bP((A_9^k)^c \cup (A_{10}^k)^c) \leq p$.

We decrease the value of $\eps$ once again so that the argument of
$\tau^+$ in the definition of $T_k$ is less than the constant
$\eps_0$ in the statement of Lemma 3.6 for $k\geq k_0$, assuming
that the constant $c_1$ in Lemma 3.6 takes the value
$c_2\delta/3$.

By our choice of the parameters we arrive at the following bound,
 $$\eqalignno{
 \bP\Big(& A_1 \cap A_2 \cap A_3 \cap \wt A_3 \cap A_4 \cap
 A_5 \cap A_6 \cap \wt A_6  \cr
 & \cap A'_1 \cap A'_2 \cap A'_3 \cap \wt A'_3 \cap A'_4 \cap
 A'_5 \cap A'_6 \cap \wt A'_6  \cr
 & \cap \bigcap_{k=0}^\infty (A_7^k \cap A_8^k \cap \wt A_8^k
 \cap A_9^k \cap A_{10}^k) \Big) > 1- 21 p. &(4.33)
 }$$

\medskip
\noindent {\it Step 3}. We will assume that all parameters have
the values chosen in the previous step and that all events that
appear in (4.33) hold. Given this assumption, we will prove that
$(\log \rho_t) /t \in[-\Lambda(D) - 8 \delta, -\Lambda(D) + 8
\delta]$ for all $t\geq \sigma^X(k_0)$.

Recall that $\dist(X_0,Y_0) = \eps>0$. First we will deal with the
case $t=\sigma^X(k_0)$. Since $A_3, \wt A_3, A_6$ and $\wt A_6$
hold,
 $$\left |\int_{{\cal T} \cap [0,\sigma^X(k_0)]} \nu(X_s) dL^X_s
 - \sigma^X(k_0) {1\over 2|D|}\int_{\prt D} \nu(y)dy \right|
 \leq 2 \sigma^X(k_0) \delta,\eqno(4.34)$$
 and
 $$\left |\int_{{\cal T} \cap [0,\sigma^X(k_0)]} \nu(Y_s) dL^Y_s
 - \sigma^X(k_0) {1\over 2|D|}\int_{\prt D} \nu(y)dy \right|
 \leq 2 \sigma^X(k_0) \delta.\eqno(4.35)$$
Since $A_1$ and $A_2$ hold, we can use Lemma 3.6 and (4.30) to
conclude that
 $$\eqalign{
  \left |\log \wt \rho_{\sigma^X(k_0) }
 - (1/2)\int_{{\cal T} \cap [0,\sigma^X(k_0)]}
 (\nu(X_s) dL^X_s +\nu(Y_s) dL^Y_s) \right|
 &\leq (c_2\delta/3)(k_0 + L^Y_{\sigma^X(k_0})\cr
 &\leq \sigma^X(k_0) \delta.
 }$$
This and (4.34)-(4.35) yield,
 $$ \left |\log \wt \rho_{\sigma^X(k_0) }
 - \sigma^X(k_0) {1\over 2|D|}\int_{\prt D} \nu(y)dy \right|
 \leq 5 \sigma^X(k_0) \delta.\eqno(4.36)$$
Since $A_4$ and $A_5$ are assumed to hold,
 $$\left |\log \ol \rho_{\sigma^X(k_0)}
 - \sigma^X(k_0) {1\over 2|D|}\int_{\prt D} \int_{\prt D}
 |\log \cos\alpha(z,y)|\omega_z(dy) dz \right| \leq 2 \sigma^X(k_0) \delta.$$
This and (4.36) imply
 $$\eqalign{ &\Big |\log  \rho_{\sigma^X(k_0)}
 - \sigma^X(k_0) {1\over 2|D|}\int_{\prt D} \nu(y)dy
 - \sigma^X(k_0) {1\over 2|D|}\int_{\prt D} \int_{\prt D}
 |\log \cos\alpha(z,y)|\omega_z(dy) dz \Big| \cr
 &\leq 7 \sigma^X(k_0) \delta.}$$
Hence we have  $\log \rho_{\sigma^X(k_0)} /\sigma^X(k_0)
\in[-\Lambda(D) - 7\delta, -\Lambda(D) + 7\delta]$. Using events
$A'_j$ and $\wt A'_j$ in place of $A_j$ and $\wt A_j$, we can also
prove that $\log \rho_{\sigma^X(k_0-1)} /\sigma^X(k_0-1)
\in[-\Lambda(D) - 7\delta, -\Lambda(D) + 7\delta]$.

Suppose that
 $$\log \rho_{\sigma^X(k_0 + j)} /\sigma^X(k_0+j)
 \in[-\Lambda(D) - 7\delta, -\Lambda(D) + 7\delta]
 $$
for some $k\geq 0$ and all $j\leq k$. We will show that the same
holds for $j=k+1$.

The event $A^k_7$ holds so
$$L^X_{\sigma^X(k_0+k+1)} - L^X_{\sigma^X(k_0+k)} +
L^Y_{\sigma^X(k_0+k+1)} - L^Y_{\sigma^X(k_0+k)} \leq 1 + (c_3 \log
(k+2))^2. \eqno(4.37) $$
 By the induction assumption,
 $$\dist(X_{\sigma^X(k_0+k)}, Y_{\sigma^X(k_0+k)}) \leq \eps
 \exp(-\sigma^X(k_0+k)(\Lambda(D)-7\delta)).\eqno(4.38)$$
This, (4.37) and Lemma 3.8 imply that
 $$\sup_{t \in [\sigma^X(k_0+k), \sigma^X(k_0+k+1)]}
 \dist(X_{t}, Y_{t}) \leq \eps
 \exp(-\sigma^X(k_0+k)(\Lambda(D)-7\delta)
 + c_4 (1 + (c_3 \log (k+2))^2).$$
This and the assumption that $A^k_7$ holds show that $T_k =
\sigma^X(k_0+k+1)\leq T_k^*$.

By (4.38) and (4.30),
 $$\dist(X_{\sigma^X(k_0+k)}, Y_{\sigma^X(k_0+k)}) \leq \eps
 \exp(-c_2 (k_0+k)(\Lambda(D)-7\delta)),$$
so the indicator functions in the definitions of events $A^k_8$
and $\wt A^k_8$ take values 1. Since these events are assumed to
hold, we obtain
 $$ \left|\int_{{\cal T}^c \cap [\sigma^X(k_0+k),\sigma^X(k_0+k+1)]}
 \nu(X_s) dL^X_s \right| \leq c_2\delta,
$$
and
 $$ \left|\int_{{\cal T}^c \cap [\sigma^X(k_0+k),\sigma^X(k_0+k+1)]}
 \nu(Y_s) dL^Y_s \right| \leq c_2\delta. $$
By the induction assumption, these estimates hold for all
$n=0,1,\dots, k$ in place of $k$, so
 $$ \left|\int_{{\cal T}^c \cap [0,\sigma^X(k_0+k+1)]}
 \nu(X_s) dL^X_s \right| \leq c_2(k+1)\delta,
$$
and
 $$ \left|\int_{{\cal T}^c \cap [0,\sigma^X(k_0+k+1)]}
 \nu(Y_s) dL^Y_s \right| \leq c_2(k+1)\delta. $$
This, (4.30), and the inequalities in $A_3$ and $\wt A_3$ imply
that
 $$\left |\int_{{\cal T} \cap [0,\sigma^X(k_0+k+1)]} \nu(X_s) dL^X_s
 - \sigma^X(k_0+k+1){1\over 2|D|} \int_{\prt D} \nu(y)dy \right|
 \leq 2 \sigma^X(k_0+k+1) \delta,\eqno(4.39)$$
and
 $$\left |\int_{{\cal T} \cap [0,\sigma^X(k_0+k+1)]} \nu(Y_s) dL^Y_s
 - \sigma^X(k_0+k+1) {1\over 2|D|}\int_{\prt D} \nu(y)dy \right|
 \leq 2 \sigma^X(k_0+k+1) \delta.\eqno(4.40)$$
Since $T_k = \sigma^X(k_0+k+1)\leq T_k^*$ and $A_2$ holds, Lemma
3.6 and (4.39)-(4.40) imply that,
 $$\eqalignno{
 &\left |\log \wt \rho_{\sigma^X(k_0+k+1)}
 - \sigma^X(k_0+k+1) {1\over 2|D|}\int_{\prt D} \nu(y)dy \right|&(4.41)\cr
 &\qquad\leq (c_2\delta/3)(k_0 + L^Y_{\sigma^X(k_0})
 + 4\sigma^X(k_0+k+1) \delta
 \cr
 &\qquad\leq 5 \sigma^X(k_0+k+1) \delta.
 }$$

In view of the assumption that $A_4, A_5$ and $A_{10}^j$ hold for
all $j=0,\dots, k$, and using (4.30), we have
 $$\eqalign{&\left |\log \ol \rho_{\sigma^X(k_0+k+1)}
 - \sigma^X(k_0+k+1) {1\over 2|D|}\int_{\prt D} \int_{\prt D}
 |\log \cos\alpha(z,y)|\omega_z(dy) dz \right|\cr
 &\leq 2 \sigma^X(k_0+k+1) \delta.}$$
This combined with (4.41) shows that
 $$ \eqalign{\Bigg |\log & \rho_{\sigma^X(k_0+k+1)}
  - \sigma^X(k_0+k+1) {1\over 2|D|}\int_{\prt D} \nu(y)dy \cr
  &- \left.\sigma^X(k_0+k+1) {1\over 2|D|}\int_{\prt D} \int_{\prt D}
 |\log \cos\alpha(z,y)|\omega_z(dy) dz \right| \leq 7 \sigma^X(k_0+k+1) \delta.}$$
In other words, $\log \rho_{\sigma^X(k_0+k+1)} /\sigma^X(k_0+k+1)
\in[-\Lambda(D) - 7\delta, -\Lambda(D) + 7\delta]$. This completes
the induction step.

We have proved that if the events in (4.33) hold then $\log
\rho_{\sigma^X(k_0+k)} /\sigma^X(k_0+k) \in[-\Lambda(D) - 7\delta,
-\Lambda(D) + 7\delta]$ for all integer $k\geq0$. We will extend
this claim to all real $t$ greater than some $t_1<\infty$. By
Lemma 4.2 (ii),
$\lim_{k\to\infty}\sigma^X(k_0+k+1)/\sigma^X(k_0+k) = 1$, a.s.
Lemma 3.8 and (4.37) imply that for $t\in [\sigma^X(k_0+k),
\sigma^X(k_0+k+1)]$ we have
 $$\eqalign{
 \log \rho_t -\log
 \rho_{\sigma^X(k_0+k)} &\leq C_9(1 + (c_3 \log(k+2))^2),\cr
 \log \rho_{\sigma^X(k_0+k+1)} -\log \rho_t
 &\leq C_9(1 + (c_3 \log(k+2))^2).\cr
 }$$
These observations easily imply that for some $t_1<\infty$ and all
real $t\geq t_1$ we have $(\log \rho_t) /t \in[-\Lambda(D) -
8\delta, -\Lambda(D) + 8\delta]$. Recall from Step 2 that this
holds with probability greater than $1-21p$. Since $p$ and
$\delta$ are arbitrarily small, the proof is complete. \qed

\vskip1truein \centerline{REFERENCES}
\bigskip

\item{[A]} H.~Aikawa, Potential-theoretic characterizations of
nonsmooth domains. {\it Bull. London Math. Soc. \bf 36} (2004),
469--482.

\item{[BB]} R.~Ba\~nuelos and K.~Burdzy, On the "hot spots"
conjecture of J. Rauch {\it J. Func. Anal. \bf 164}, (1999) 1--33.

\item{[Ba]} R.F. Bass, {\it Probabilistic Techniques in Analysis},
Springer, New York, 1995.

\item{[BBC]} R.~Bass, K.~Burdzy and Z.-Q.~Chen (2002) Uniqueness
for reflecting Brownian motion in lip domains.
{\it Ann. Inst. H. Poincar\'e Probab. Statist.} (to appear)

\item{[BH]} R.F.~Bass and P.~Hsu, Some potential theory for
reflecting Brownian motion in H\"older and Lipschitz domains {\it
Ann. Probab. \bf 19} (1991) 486--508.

\item{[B]} K. Burdzy, {\it Multidimensional Brownian excursions
and potential theory}, Longman Sci. Tech., Harlow, 1987.

\item{[BC]} K. Burdzy\ and\ Z.-Q. Chen, Coalescence of synchronous
couplings, {\it Probab. Theory Related Fields \bf 123} (2002),
no.~4, 553--578.

\item{[CLJ1]}M. Cranston\ and\ Y. Le Jan, On the noncoalescence of
a two point Brownian motion reflecting on a circle, {\it Ann. Inst. H.
Poincar\'e Probab. Statist. \bf 25} (1989), no.~2, 99--107.

\item{[CLJ2]} M. Cranston\ and\ Y. Le Jan, Noncoalescence for the
Skorohod equation in a convex domain of ${\bf R}\sp 2$, {\it
Probab. Theory Related Fields \bf 87} (1990), 241--252.

\item{[FOT]} M.~Fukushima, Y.~Oshima and M.~Takeda,
 {\it Dirichlet Forms and Symmetric Markov Processes}.
 Walter de Gruyter, Berlin, 1994

\item{[H]} E.~Hsu (2002) Multiplicative functional for the heat
equation on manifolds with boundary, {\it Michigan Math. J. \bf
50} (2002), 351--367.

\item{[KS]} I.~Karatzas and S.E.~Shreve {\it Brownian Motion and
Stochastic Calculus}, Second edition, \break
 Springer, New York, 1991.

\item{[LS]} P.~L. Lions and A.~S. Sznitman,  Stochastic differential
equations with reflecting boundary conditions. {\it Comm. Pure Appl.
Math.} {\bf 37} (1984), 511-537.

\item{[M]} B. Maisonneuve, Exit systems, {\it Ann. Probability \bf 3}
(1975), no.~3, 399--411.

\item{[T]} M.~Tsuji, {\it Potential Theory in Modern Function
Theory}, Maruzen, Tokyo, 1959.

\item{[Z]} Z.~Zhao, Uniform boundedness of conditional gauge
          and Schr\"{o}dinger equations
 {\it Commun. Math. Phys. \bf 93} (1984), 19--31.

\bigskip

K.B. and Z.C. Address: Department of Mathematics, Box 354350,
University of Washington, Seattle, WA 98115-4350 \hfill\break
(burdzy@math.washington.edu, zchen@math.washington.edu)

P.J. Address: Department of Mathematics, Yale University, PO Box
208283, New Haven, CT 06520-8283 \hfill\break
(jones@math.yale.edu)

\medskip

\bye